\documentclass[11pt]{article}

\usepackage[all]{xy}

\newtheorem{tm}{Theorem}[subsection]
\newtheorem{lm}[tm]{Lemma}
\newtheorem{pr}[tm]{Proposition}
\newtheorem{rmk}[tm]{Remark}
\newtheorem{cor}[tm]{Corollary}

\newtheorem{fact}[tm]{Fact}
\newtheorem{??}[tm]{Question}
\newtheorem{defi}[tm]{Definition}

\newtheorem{ass}[tm]{Assumption}

\font\tenmsb=msbm10
\font\sevenmsb=msbm7
\font\fivemsb=msbm5

\newfam\msbfam
\textfont\msbfam=\tenmsb
\scriptfont\msbfam=\sevenmsb
\scriptscriptfont\msbfam=\fivemsb
\def\Bbb#1{{\fam\msbfam #1}}

\font\teneufm=eufm10
\font\seveneufm=eufm7
\font\fiveeufm=eufm5
\newfam\eufmfam
\textfont\eufmfam=\teneufm
\scriptfont\eufmfam=\seveneufm
\scriptscriptfont\eufmfam=\fiveeufm
\def\frak#1{{\fam\eufmfam\relax#1}}

\def\lorw{\longrightarrow}
\newcommand\n{\noindent}

\newcommand\ci{\cite}

\newcommand\rat{{\Bbb Q}}
\newcommand\comp{{\Bbb C}}
\newcommand\real{{\Bbb R}}
\newcommand\zed{{\Bbb Z}}
\newcommand\nat{{\Bbb N}}
\newcommand\pn[1]{{\Bbb P}^{#1}}

\newcommand\blacksquare{{\hspace*{\fill} $\fbox{}$}}

\newcommand\f{\varphi}

\newcommand\e{\epsilon}

\newcommand{\fxn}[3]{{#1}_* \rat_{#2} [{#3}]  }

\newcommand{\phix}[2]{ \,^p\!{\cal H}^{#1}(#2)} 

\newcommand{\lij}[3]{ {#1}^{#2}_{#3} }

\newcommand{\dsdix}[2]{ \bigoplus_{#1} \phix{#1}{#2} [-{#1} ]  }

\newcommand{\dsdjix}[3]{ \bigoplus_{#1} \phix{#2}{#3} [-{#2} ]  }

\newcommand{\ihixc}[3]{ {\Bbb H}^{#1} ( {#2}, {#3} ) }

\newcommand{\im}{ \hbox{\rm Im} }

\newcommand{\ke}{ \hbox{\rm Ker} }

\newcommand{\p}{{\Bbb P}}

\newcommand{\pd}{ {\Bbb P}^{\vee} }

\newcommand{\lijx}[4]{ {#1}^{#2}_{#3} ( {#4} ) }

\newcommand{\csix}[2]{ {\cal H }^{#1} ( {#2} ) }

\newcommand{\ptd}[1]{ \,^{p}\!\tau_{ \leq {#1} } }

\newcommand{\ptu}[1]{ \,^{p}\!\tau_{ \geq {#1} } }

\newcommand{\td}[1]{ \tau_{ \leq {#1} } }

\newcommand{\tu}[1]{ \tau_{ \geq {#1} } }

\newcommand{\fihixcj}[4]{ {\Bbb H}^{#1}_{\leq #2} ( {#3}, {#4} ) }

\newcommand{\gihixcj}[4]{ {\Bbb H}^{#1}_{#2} ( {#3}, {#4} ) }

\title{The Hodge theory of algebraic maps}
\author{
Mark Andrea A.  de Cataldo\thanks{
Partially supported by N.S.F. Grant DMS 0202321, NSA Grant 
MDA904-02-1-0100.}\, 
and Luca Migliorini\thanks{ Partially supported by MIUR project {\em 
Propriet\`a Geometriche delle Variet\`a Reali e Complesse} and by GNSAGA}
}
\date{May, 2004}

\begin{document}\maketitle

\begin{abstract}
We give a geometric proof 
of the Decomposition Theorem of Beilinson,
Bernstein, Deligne and Gabber for the direct image
of the intersection cohomology complex
under a proper map of complex algebraic varieties.
The method rests on new Hodge-theoretic results
on the cohomology of projective varieties which extend naturally
the classical theory and provide new applications.
\end{abstract}

\tableofcontents

\section{Introduction}
 Let $f:X \to Y$ be a proper 
map of complex algebraic varieties,
$n =\dim{X}.$ For ease of exposition only, assume that
$X$ is  nonsingular and that $X$ and $Y$ are projective.
 Let $\eta$ and $A$ be ample 
line bundles on $X$ and $Y$ respectively, and set
$L:=f^{*}A.$ 

\medskip
If $f$ is a smooth 
family, then the classical Hard Lefschetz Theorem for
$\eta$ applied to the fibers of $f$ 
gives isomorphisms for every $i \geq 0:$ 
\begin{equation}
    \label{101}
    \eta^{i}: R^{n - \dim{Y} -i} f_{*}\rat_{X} \simeq
R^{n - \dim{Y} +i} f_{*}\rat_{X}
\end{equation}
which 
give rise to a direct sum
decomposition for the direct image complex
\begin{equation}
\label{202}
Rf_*\rat_X \simeq \bigoplus_{i}{ R^if_*\rat_X[-i]}
\end{equation}
in the derived category of the category of sheaves on $Y$
(cf. \ci{dess}). This important fact
implies, for example, the $E_{2}-$degeneration of the Leray spectral 
sequence
for $f.$
The sheaves $R^{i}f_{*}\rat_{X}$ are semisimple local systems,
i.e. they split as a direct sum of local systems
with no nontrivial local subsystems.
Note  that the 
category of finite dimensional local systems is  abelian, noetherian
and artinian.

\medskip
At first sight, nothing similar happens for an arbitrary
map $f:X \to Y.$ The isomorphisms (\ref{101}) and (\ref{202}) fail 
in 
general, 
  the Leray spectral sequence 
  may not degenerate at $E_2$ and the abelian  category of sheaves on 
$Y$
  is neither noetherian, nor artinian.
 
  \n
The Leray spectral sequence is 
associated with the ``filtration'' of $Rf_* \rat_X$ by the truncated 
complexes 
$\tau_{\leq i}Rf_* \rat_X.$  The $i-$th direct image $R^{i}\!f_{*}\rat_{X}$
appears, up to a shift, as the cone of 
the natural map 
$\tau_{\leq i-1}Rf_* \rat_X \to \tau_{\leq i} Rf_* \rat_X,$
i.e. as the  $i-$th cohomology sheaf of the complex
 $Rf_* \rat_X.$

\medskip
One of the main 
ideas  leading to the theory of perverse 
sheaves in \ci{bbd} is that all the  facts mentioned in the case of a smooth 
family hold for an arbitrary
map, provided that   they are  re-formulated with respect to a 
notion of truncation different from the one leading
to the cohomology sheaves,  that is with respect
to the so-called perverse truncation 
$\ptd{i},$ 
and that  we replace the sheaves $R^if_* \rat_{X}$
with the shifted cones $\phix{i}{Rf_{*}\rat_{X}} $ of the mappings 
$\ptd{i-1}Rf_{*}\rat_{X} 
\lorw
\ptd{i} Rf_{*}\rat_{X}.$ 
These cones are called the perverse cohomology  of $Rf_* 
\rat_X$ and  are perverse sheaves.
Despite their name, perverse sheaves are 
complexes in  the derived category of the category of sheaves on 
$Y$  which are characterized by conditions on 
their cohomology sheaves. Just like local 
systems, the category of perverse sheaves is abelian, noetherian
and artinian. 
Its simple objects are the  intersection 
cohomology complexes  of simple local systems on strata.
Whenever $Y$ is 
nonsingular and the stratification is trivial, perverse sheaves are, 
up to 
a shift, just local systems.

\medskip
That these notions are the correct  generalization
to arbitrary proper morphisms of the situation
considered above for   smooth morphisms, 
is shown by the beautiful {\em Relative Hard Lefschetz Theorem}
and {\em Decomposition Theorem},
proved in \ci{bbd} by Beilinson, Bernstein and Deligne
using algebraic geometry in positive characteristic. They generalize 
the  isomorphisms (\ref{101}) and (\ref{202})  for a smooth family
to the case of an arbitrary projective  map from an algebraic manifold:
the map induced by the line bundle  $\eta$ in perverse cohomology
\begin{equation}
    \label{303}
\eta^{i}: \phix{-i}{\fxn{Rf}{X}{n}} \lorw 
\phix{i}{\fxn{Rf}{X}{n}}
\end{equation}
is an isomorphism for every $i \geq 0$ and we have a direct sum 
decomposition
\begin{equation}
\label{404}
  \fxn{Rf}{X}{n} \simeq \dsdix{i}{\fxn{Rf}{X}{n}} .
\end{equation}
 As a consequence,  the so-called
perverse Leray spectral sequence $ \ihixc{l}{Y}
{\phix{m}{\fxn{Rf}{X}{n}}} $ $  \Longrightarrow H^{n+l+m}(X,\rat)$
is $E_{2}-$degenerate. This fact alone has striking
computational and theoretical consequences. For example,
the intersection cohomology groups of a variety $Y$
inject in the ordinary singular cohomology groups
of any resolution $X$ of the singularities of $Y.$
The semisimplicity statement for the local 
systems $R^i\!f_*\rat_X$ 
has a far-reaching generalization in the
{\em  Semisimplicity Theorem}, also proved in \ci{bbd}:
there is a canonical isomorphism of perverse sheaves
\begin{equation}
 \label{505}
\phix{i}{\fxn{Rf}{X}{n}} \simeq \bigoplus_{l}{
IC_{\overline{S_{l}}} (L_{i,l}) }
\end{equation}
where  the  $IC_{\overline{S_{l}}} (L_{i,l})$ 
are the  Goresky-MacPherson intersection cohomology complexes
on  $Y$ associated with
certain semisimple local systems $L_{i,l}$ on the strata of
a  finite algebraic
stratification $Y= \amalg_{l=0}^{\dim{Y}}{
S_{l}}$ for the map $f.$

\n
Analogous results hold for a possibly singular $X,$ provided one
replaces $Rf_{*}\rat_{X}[n]$ by the intersection cohomology
complex $IC_{X}.$

\medskip
These three theorems are  cornerstones of the topology of algebraic 
maps.
They have found many applications to algebraic geometry and 
to representation theory
and, in our opinion, should be regarded as expressing  fundamental 
properties of complex algebraic geometry.

\bigskip
In our previous paper \ci{demigsemi} we proved that if $f$ is 
semismall, then 
$L$ behaves Hodge-theoretically like an ample line bundle:
the Hard 
Lefschetz Theorem holds for $L$ acting on rational cohomology,
i.e. $L^{r}: H^{n-r}(X) \simeq H^{n+r}(X),$
for every
$r,$  and the  primitive 
subspaces $\ke{\,L^{r+1}} \subseteq H^{n-r}(X)$
are  polarized by means of the
 intersection form on $X.$
Associated with a stratification of the map $f$ there 
is a series of intersection forms  describing
how the fiber of a point in a given 
stratum intersects in $X$ the pre-image of the stratum.
In the case of a semismall map there is only one intersection form 
for each stratum component.
The discovery of the  above polarizations,
joined with an argument of mixed Hodge 
structures showing that  for every $y \in Y$ the group $H_{n}(f^{-1}(y))$
injects in $H^{n}(X),$
allowed us to prove that the intersection  forms are  
definite and thus  nondegenerate. 
This generalizes the well-known result of Grauert
for the contraction of curves on surfaces.
By means of an induction on the strata,
the  statement of the 
semisimplicity  theorem  (\ref{505}) for semismall maps $f$ 
was proved to be equivalent to the fact that the intersection forms 
are nondegenerate. The statements (\ref{303}) and (\ref{404})
are trivial for semismall maps since
$Rf_{*}\rat_{X}[n] \simeq \phix{0}{Rf_{*}\rat_{X}[n]}.$
Our  result about the intersection forms being  definite
can  be seen as a 
``Decomposition Theorem with signs for semismall maps,'' i.e. 
as a polarized version of 
this theorem. 

\bigskip

In this paper, 
 in the spirit of our  paper \ci{demigsemi},
 we give a  geometric proof 
  of  the Relative Hard Lefschetz, Decomposition and
 Semisimplicity isomorphisms
 (\ref{303}),  (\ref{404}) and (\ref{505}). We complement
 these results by uncovering 
 a series of Hodge-theoretic
properties of  the singular rational cohomology groups  $H^*(X)$
and of the natural maps 
 $H_{n-*}(f^{-1}(y)) \to H^{n+*}(X),$ $y \in Y.$

\medskip
We now discuss our results. By standard reductions, most statements
remain valid in the context of proper maps of algebraic varieties
(cf. \ref{tacss}).

\n
The perverse truncation, which  is defined locally over $Y$ by means of 
the topological operations of push-forward and truncation
with respect to a stratification for $f,$
gives rise to an increasing filtration
$H^l_{\leq b}(X)
\subseteq  H^{l}(X)$ and to the corresponding
graded perverse cohomology groups
 $H^l_{b}(X)=H^l_{\leq b}(X)/H^l_{\leq b-1}(X).$ See \ref{tpf}.

\medskip
Our first result is
the {\em Hard Lefschetz Theorem for Perverse Cohomology Groups} 
\ref{tm3}.
While  the Hard Lefschetz Theorem
for  the pull-back line bundle $L$  acting on ordinary cohomology 
fails, due to
the  lack of positivity of $L$ on the fibers of the map $f,$
the analogous result, with a natural shift in cohomological
degree, holds for the
perverse cohomology groups. 

\n
{\em It is  as if the perverse 
filtration
were calibrated precisely for the purpose of correcting
the failure of the Hard Lefschetz Theorem for $L.$}
In fact,
this result implies that the perverse filtration coincides with
the canonical weight filtration associated with the Jordan form
of the nilpotent operator $L$ acting on $H^*(X)$ (cf. \ref{tiafiltr}).
We find this aspect of our approach quite intriguing.

\smallskip
Since the weight filtration above, being
characterized in terms of the $(1,1)-$operator $L,$
is automatically
Hodge-theoretic, we get the
{\em Hodge Structure Theorem} \ref{uf} stating
that the perverse filtration, and hence the perverse cohomology 
groups,
are endowed with canonical Hodge structures.

\smallskip

The Hard Lefschetz Theorem for Perverse Cohomology Groups
implies
the {\em $(\eta,L)-$De\-com\-po\-si\-tion Theorem} \ref{etaldecompo},
i.e. a  Lefschetz-type direct sum 
decomposition of  
the perverse cohomology groups  of $X$
into ``$(\eta,L)-$primitive'' 
 Hodge sub-structures. The decomposition
 is orthogonal with respect to
certain polarizing  bilinear forms $S^{\eta L}$
coming from the Poincar\'e pairing
on $X$ modified  by $L$ and $\eta.$

\smallskip

The {\em Generalized Hodge-Riemann Bilinear Relations Theorem} \ref{tmboh}
state that the forms $S^{\eta L}$ polarize up to sign the
$(\eta,L)-$primitive spaces.

\n
In our approach, it is crucial
to describe (cf.  \ref{tapwat})
the subspace of cohomology classes of 
$H^n(X)$ which are limits
for $\epsilon \to 0^+$ of cohomology classes primitive 
with respect to the ample line  bundles of the form $L + \epsilon 
\eta .$

\smallskip

The {\em Generalized Grauert Contractibility Criterion} \ref{nhrbr}
and the {\em Refined Intersection Form Theorem} \ref{rcffv} establish
some of the Hodge-theoretic properties of the homology groups
$H_{*}(f^{-1}(y))$ and of
the  refined intersection forms defined on them (cf. \ref{sperem}).

\medskip

To our knowledge, this rich structure 
on $H^*(X)$ and on $H^{*}(f^{-1}(y))$
has not been spelled-out before and it should have 
significant geometric applications. We discuss
two examples in  \ref{example3fold}
and \ref{examplefamily}.

\n
We propose two applications, the {\em Contractibility Criterion} \ref{solido}
and the {\em Signature  for Semismall Maps Theorem } \ref{sigsemi}.

\medskip
In our approach, not only are 
 these structures and results  
complementary to the Decomposition Theorem, but they are also 
instrumental
in proving it. 
In view of the inductive approach we develop
in \ref{redtopt},
the Decomposition Theorem implies that
 the refined intersection forms
(more precisely, one of the graded parts) are 
nondegenerate.  We establish the converse statement. 
We prove  directly
that these forms  are nondegenerate and
 show how this nondegeneration  implies the 
Decomposition Theorem. In the critical case of 
cohomological degree $n$ and  perversity zero, we
show that the graded class map
$H_{n,0}(f^{-1}(y)) \to H^{n}_{0}(X)$ is an injection of pure Hodge 
structures
and that the refined  intersection form on the fiber $f^{-1}(y)$
underlies a polarization.
Again, our results, can be seen as a ``Decomposition Theorem with 
signs.''

\smallskip
These results, coupled with a  series of simple reductions, 
give a proof of (\ref{404}), 
(\ref{505}) 
for proper maps of complex algebraic varieties and
of (\ref{303}) for projective maps of complex algebraic varieties
(cf. \ref{tacss}).
 
\smallskip
The {\em Purity Theorem}
\ref{tmp1}  states that the direct sum decomposition
for the hypercohomology of (\ref{505})
is by Hodge sub-structures. In particular,
the intersection cohomology groups
of projective varieties carry a pure Hodge structure
which is canonical in the sense of \ref{pic}.a.

\smallskip
The {\em Hodge-Lefschetz Theorem 
for Intersection Cohomology} \ref{pic} is a generalization to the 
intersection cohomology 
of a projective variety of the classical Hodge theory
for the singular cohomology 
of projective manifolds (cf. \ref{thtav}).

\bigskip
What follows compares the results of the present paper
with some of the literature.

\n
Theorem
 \ref{tm1ac}
is proved by Beilinson, Bernstein, Deligne and Gabber in
\ci{bbd}.
The result  is first  proved in positive characteristic
using  the formalism of perverse sheaves
in conjunction with
the purity results proved in \ci{weil2} concerning
the eigenvalues of 
the Frobenius operator acting on complexes  of sheaves
on a variety defined 
over a finite field.
The result is then ``lifted'' to characteristic zero.

\n
The deep elegance of this approach does not seem to
explain the geometry of the result over the complex numbers
and does not give 
a proof of the Hodge-theoretic results in \ref{tpcss}.

\n
In  the series of remarkable papers
 \ci{samhp}, 
\ci{samhm} and \ci{sa}, M. Saito 
has   developed   a D-modules transcendental  approach
via his own mixed Hodge modules.

\n
C. Sabbah \ci{sabbah} has recently extended M. Saito's results
to the case of semisimple local systems on $X$ by developing
his own theory of polarizable twistor D-modules.
See also the related work by T. Mochizuki
\ci{mo1}, \ci{mo2}.

\n
While Saito's results cover and pre-date some of the results of this paper,
namely the Relative Hard Lefschetz Theorem \ref{tm1}.a,
the Decomposition Theorem \ref{tm1}.b and part of the Purity
of Intersection Cohomology \ref{pic}.a,
and give other results as well,
the proofs do not seem to explain the underlying geometry and do not
describe explicitly Hodge structures and polarizations.

\medskip 
We show that the properties of
the   refined intersection form
are responsible for the topological splitting
of $R\fxn{f}{X}{n}$ and we establish these properties
using Hodge theory.

\bigskip
  The paper is not self-contained
as it relies, for instance, on the theory of $t-$structures. 
However, at several stages,  we 
need  results in  a form that seems to be less
general but sharper than what we could find
in the literature.
For this reason, we offer two 
 rather long sections of preliminaries.
 We also  hope
that having 
collected results on   the theory of stratifications, 
constructible sheaves and perverse sheaves can in any case be useful 
to 
the reader.  The statements proved in this paper are 
collected in section \ref{stats}. 
The proofs of the main results are strongly intertwined
and    we give
a detailed account of the steps of the proof
in 
\ref{ssotp}, trying 
to emphasize the main ideas.
Due to the presence of a rich array of structures, many verifications 
of compatibility are necessary in the course of our proofs.
  We have decided to include careful proofs of the ones that
did not seem to be 
just routine.

\bigskip
\n
{\bf Acknowledgments.} The authors 
would like to thank
Victor Ginzburg, Mark Go\-re\-sky, Agnes Szilard, Arpad Toth and
Dror Varolin for useful remarks. They would like to thank
Robert MacPherson for suggesting improvements and making useful remarks.
The first-named author would like to thank the following institutions for
their kind hospitality while portions of this work were being carried
out: the Korean Institute for Advanced Study, Seoul,
the Institute for Advanced Study, Princeton,
the Mathematical Sciences Research Institute, Berkeley,
the Max Planck Institut
f\"ur Mathematik, Bonn and the Department of Mathematics
of the University of Bologna.
The second-named author  would like to thank
the Department of Mathematics
of SUNY at Stony Brook and the Department of Mathematics of the University of 
Rome, ``La Sapienza.''

\n
The first-named author dedicates this paper to 
his family and to the memory
of Meeyoung Kim.

\n
The second-named author dedicates this paper to his father 
for his eightieth birthday, 
and to E., F. and G.

\section{Statements}
\label{stats}
We state and prove our results  for maps of
projective varieties $f: X \to Y$  
 with $X$ nonsingular
in \ref{tpcss} and \ref{pp}. The
Hodge-theoretic results are strongest and more meaningful
in this context.
Most results remain valid (cf. \ref{tacss}) for proper algebraic maps
of algebraic varieties
via standard reductions to the nonsingular projective case.

\subsection{The projective case}
\label{tpcss}
The basic set-up of this paper is as follows:
\begin{itemize}
    \item
let $f:X \to Y$ be a map of projective varieties,
$X$ nonsingular of dimension $n,$ $\eta$ be an ample line bundle on 
$X,$
$A$ be an ample line bundle on $Y$  and $L:= f^{*}A.$
\end{itemize}

\medskip
The results that follow are discussed in the  two examples \ref{example3fold}
and \ref{examplefamily}.

\medskip

We denote $Rf_{*}$ simply by $f_{*}.$
The line bundle $\eta$ defines a map $\rat_{X} \to \rat_{X}[2].$
By pushing-forward in the derived sense, we get
a map $f_{*}\rat_{X} \to f_{*}\rat_{X}[2]$ and maps
$\phix{i}{\fxn{f}{X}{n}} \to \phix{i+2}{\fxn{f}{X}{n}}.$ We denote 
all these maps simply by $\eta.$

\n
The following extends a great deal of classical Hodge Theory
to the case of maps.  It was proved by Beilinson, Bernstein
and Deligne in \ci{bbd}.

\begin{tm}
\label{tm1}
$\,$

\n
(a) ({\bf The Relative Hard Lefschetz Theorem})
For every $i\geq 0,$ the map induced  by $\eta$ in perverse 
cohomology is an isomorphism:
$$
\eta^i: \phix{-i}{ \fxn{f}{X}{n}   } \simeq \phix{i}{ \fxn{f}{X}{n}  
}.
$$
In particular, having set, for $i\geq 0,$ ${\cal P}^{-i}_{\eta}: 
= \ke \; \eta^{i+1} \subseteq \phix{-i}{ \fxn{f}{X}{n} },$
we have  equalities
$$
\phix{-i}{ \fxn{f}{X}{n} } = \bigoplus_{j\geq 0} \eta^{j} {\cal 
P}_{\eta}^{-i -2j},
\qquad
\phix{i}{ \fxn{f}{X}{n} } = \bigoplus_{j\geq 0} \eta^{i+j} {\cal 
P}_{\eta}^{-i -2j}.
$$
(b) ({\bf The Decomposition Theorem})
There is an isomorphism in $D(Y)$:
$$
\f:  \, \dsdix{i}{ \fxn{f}{X}{n} }   \,     \simeq  \,     
\fxn{f}{X}{n}.
$$
(c) ({\bf The Semisimplicity Theorem}) 
The $\phix{i}{ \fxn{f}{X}{n} }$ are semisimple (cf. 
\ref{iccso}).
More precisely, given any stratification for $f$ (cf. 
\ref{wsoav}) $Y= \amalg_{l}{S_{l}},$ $ 0 \leq l \leq  \dim{Y},$
there is a canonical isomorphism in $Perv(Y):$
$$
\phix{i}{\fxn{f}{X}{n}} \simeq
\bigoplus_{l=0}^{\dim{Y}}{ 
IC_{\overline{S_{l}}}(L_{i,l} )}
$$
where the local systems  
$L_{i,l}:= \alpha_{l}^{*} {\cal H}^{-l} (  
\phix{i}{\fxn{f}{X}{n}}  )
$ on $S_{l}$
are semisimple.
\end{tm}

\begin{rmk}
\label{dsgzero}
{\rm The complexes $\phix{i}{\fxn{f}{X}{n}}=0$ if $|i| > r(f),$
where $r(f)$ is the defect of semismallness \ref{defr} of $f.$
It can be shown that this vanishing is sharp (cf. \ci{decmightamv1}).}
\end{rmk}

\begin{rmk}
    \label{paraphi}
{\rm The isomorphism $\f$ of Theorem \ref{tm1}.b is not unique.
It is possible to make some distinguished choices
(cf. \ci{shockwave}). These choices play no role in the present
paper.
}
\end{rmk}
 
\medskip
The symbol $\stackrel{\f}{\simeq}$  indicates 
that a certain isomorphism is realized via $\f$.
The Decomposition Theorem implies that (cf. \ref{psic}), setting
$$
H^{n+l}_{\leq b}(X)\,  :=\, \im{ \,\left\{ 
\ihixc{l}{Y}{\ptd{b}  \fxn{f}{X}{n}} 
\to \ihixc{l}{Y}{ \fxn{f}{X}{n}}  \right\}} \,  \subseteq \, H^{n+l}(X)
$$
and
$$
H^{n+l}_b(X)\, : = \,
H^{n+l}_{\leq b}(X)/ H^{n+l}_{\leq b-1}(X),
$$
we get a canonical identification
$$
H^{n+l}_{b}(X) \,  =  \, {\Bbb H}^{l-b}(Y, 
   \phix{b}{\fxn{f}{X}{n}} )
$$   
and  isomorphisms:
$$
H^{n+l}_{\leq b}(X) \, \stackrel{\f}{\simeq} \,
 \bigoplus_{i \leq b}{\Bbb H}^{l-i}(Y, 
   \phix{i}{\fxn{f}{X}{n}} ).
$$
The cup product with $\eta$ verifies 
$\eta \, H^l_{\leq a}(X) \subseteq H^{l+2}_{\leq a+2}(X)$ 
(cf. \ref{cupetap})
and induces maps,  still denoted 
$\eta: H^l_{a}(X) \to H^{l+2}_{a+2}(X)$.
The cup product with $L$ is  compatible with the direct sum 
decomposition induced by any isomorphism $\f$ (cf. \ref{sepullback})
and induces maps $L: H^l_{a}(X) \to H^{l+2}_{a}(X).$

\begin{tm}
\label{tm3} ({\bf  The Hard Lefschetz Theorem for Perverse Cohomology 
Grou\-ps})
Let $k \geq 0,$ $b, \, j\, \in \zed.$ Then
 the following cup product maps are isomorphisms:

$$
\eta^{k} \, : \, H^{j}_{-k}(X) \, \simeq \, H^{j+2k}_{k}(X),
 \qquad
\qquad
L^{k}\, :\,  H^{n+b-k}_{b}  (X) \,
\simeq \,
H^{n+b+k}_{b}(X).
$$
\end{tm}
The previous result allows to describe the perverse filtration purely
in terms of the nilpotent linear map $L$ acting via cup-product
on the cohomology of $X.$ For the precise statement, involving the 
notion
of weight filtration associated with $L,$ see \ref{fs}, \ref{asoh} and 
Proposition \ref{tiafiltr}. Since $L$ is of type $(1,1),$ we get the 
following.

\begin{tm}
\label{uf}
({\bf The  Hodge Structure Theorem})
For $l\geq 0$ and $b\in \zed,$ the subspaces
$$
H^{l}_{\leq b}(X) \, \subseteq\,  H^{l}(X)
$$
are  pure Hodge sub-structures. The quotient spaces
$$
H^l_b(X)\, = \, H^{l}_{\leq b}(X) /H^{l}_{\leq b-1}(X)
$$
inherit a  pure Hodge structure 
of  weight $l.$
\end{tm}

Note that the Hodge structure thus constructed on 
$H^l_b(X)$ is compatible with restriction to zero-loci of
sections of  $\eta$ and $L$
in the following sense.
Let $i:X^r_s \to X$ be the inclusion of a
complete intersection of $r$ general sections of 
$\eta$ and $s$ general sections of $L.$
The restriction $i^*:H^l(X) \to H^l(X^r_s)$ maps $H^l_{\leq a}(X)$ to 
$H^l_{\leq a+r}(X^r_s)$ (cf. \ref{filtrfac}) 
and  induces a morphism of pure Hodge structures
$H^l_{a}(X) \to H^l_{a+r}(X^r_s).$

\medskip

Define 
$P^{-j}_{-i}:= \ke{\, \eta^{i+1}} \cap \ke{\, L^{j+1}} 
\subseteq H^{n-i-j}_{-i}(X),$ 
$i,\,j \geq 0$ and   $P^{-j}_{-i} :=0$
otherwise.
In the same way in which the classical Hard Lefschetz implies
the Primitive Lefschetz Decomposition for the cohomology of $X,$
Theorem \ref{tm3}  implies the 
double direct sum decomposition of

\begin{cor}
\label{etaldecompo}
({\bf The $(\eta,L)-$Decomposition})
Let $i, \, j \in \zed.$ 
There is  a Lefschetz-type direct sum decomposition into 
pure Hodge sub-structures  of weight $(n-i-j),$
called the  $(\eta,L)-$decomposition:
$$
H^{n-i-j}_{-i}(X) = \bigoplus_{l,\,m\, \in \zed }{
\eta^{-i+l} \, L^{-j+m} \, P^{j-2m}_{i-2l}.}
$$
\end{cor} 

Using representatives in $H^{n-i-j}_{\leq -i}(X),$
$i,j \geq 0,$ define bilinear forms 
$S^{\eta L}_{ij}$ on $H^{n-i-j}_{-j}(X)$ by 
modifying the Poincar\'e pairing
$$
S^{\eta L}_{ij} ([\alpha], [\beta]   ) \, := \,
\int_X{ \eta^i \wedge  L^j \wedge  \alpha \wedge \beta}.
$$
Using \ref{tm3}, one can define forms 
$S^{\eta L}_{ij}$ for all $i$ and $j$.
These forms are well defined and nondegenerate  (cf. \ref{fs}).
 
\begin{tm}
\label{tmboh}
({\bf  The Generalized Hodge-Riemann Bilinear Relations})

\n
The $(\eta, L)-$decomposition is orthogonal with respect
to $S^{\eta L}_{ij}.$ The forms $S^{\eta L}_{ij}$
are, up to a sign made precise in \ref{uptosign}, 
a polarization 
of each $(\eta,L)-$direct summand.
%$\eta^{-i+l} \, L^{-j+m} \, P^{j-2m}_{i-2l}.$
\end{tm}

The following two results are key in proving the Semisimplicity 
Theorem and in our opinion are geometrically significant.
For what follows see  \ref{sperem}  and Lemma \ref{collectinfo}.

\n
The homology groups $H^{BM}_{*}(f^{-1}(y))=H_{*}(f^{-1}(y)),$
$y \in Y,$ 
are filtered.
The natural cycle class map
$cl: H^{BM}_{n-*}(f^{-1}(y)) \to H^{n+*}(X)$ is strict.

\begin{tm}
\label{nhrbr} 
({\bf The Generalized Grauert Contractibility Criterion})
Let $b \in \zed$, $y \in Y$.
 The  natural class  maps
$$
cl_{b}: H^{BM}_{n-b,b}(f^{-1}(y)) \lorw H^{n+b}_b(X)
$$
is injective and identifies
$H^{BM}_{n-b,b}(f^{-1}(y))\subseteq \ke{\, L}\subseteq
H^{n+b}_b(X)$ with a pure  Hodge substructure, compatibly with
the $(\eta,L)-$decomposition. Each $(\eta,L)-$direct summand
of $H^{BM}_{n-b,b}(f^{-1}(y))$ is polarized up to sign
by $S^{\eta L}_{-b,0}.$

\n
In particular, the restriction of $S^{\eta L}_{-b,0}$ to
$H^{BM}_{n-b,b}( f^{-1} (y ))$ 
is nondegenerate.
\end{tm}

By intersecting in $X$ cycles supported on $f^{-1}(y),$
we get the refined intersection form (cf. \ref{sperem})
$H^{BM}_{n-*}(f^{-1}(y)) \to H^{n+*}(f^{-1}(y))$ 
which is strict as well.

\begin{tm}
\label{rcffv}
({\bf The Refined Intersection Form Theorem})
Let $b \in \zed$, $y \in Y$.
The graded refined intersection form
$$
H^{BM}_{n-b,a}(f^{-1}(y)) \lorw H^{n+b}_{a}(f^{-1}(y))
$$
is zero if $a\neq b$ and it is an isomorphism if $a=b.$
 \end{tm}

\begin{rmk} 
    \label{alk}
{\rm
Let $b \in \zed$, $y \in Y$. If $y$ lies on a positive-dimensional 
stratum $S_{l}$ for the map $f$ (cf. \ref{wsoav}), then
the refined intersection form is identically zero (cf. \ref{1111}).
In this case, the geometrically interesting map
is the one of the Splitting Criterion
\ref{splitp}:
$H^{BM}_{n-l-b,b}( f^{-1}(y) ) \to H^{n-l+b}_{b}(f^{-1}(y)).$
This map is the refined intersection form for $f^{-1}(y)$
in $X_{l}$ for the map $X_{l} \to Y_{l}$, where $Y_{l}$ is a 
codimension$-l$
general complete 
intersection  on $Y$ transversal to $S_{l}$ at $y$ and $X_{l}:=
f^{-1}(Y_{l})$ (cf. \ref{dtpb}).
}
\end{rmk}
 
\medskip
The following two corollaries give  examples of
 the kind of geometric applications stemming from 
the results of this paper. See also Example \ref{example3fold}.

\begin{cor}
\label{solido} ({\bf Contractibility Criterion})
Let $f:X \to Y$ be a projective and surjective  map of quasi projective
varieties, $X$ nonsingular. Let $y \in Y$ and assume
that $f$ is semismall over $Y\setminus y, $ e.g. $f_{| X \setminus
f^{-1}(y) }$
is an isomorphism.

\n
For every $b \geq 0$ the natural mixed Hodge structure
$
H^{n+b}(f^{-1}(y))
$
is  pure of weight $n+b.$
The associated  $\eta-$direct summands are
polarized up to sign by the forms
$S^{\eta L}_{-b 0}.$
\end{cor}
{\em Proof.}
We may compactify the morphism  without changing the situation
around and over  $y.$
In particular, the perverse cohomology complexes
will not change around $y$ by this process.
It follows that we may assume   that $X$ is projective.
 By the Decomposition Theorem,
the Semisimplicity Theorem and  semismallness (cf. \ref{defr} and 
\ref{sstri})
$$
\fxn{f}{X}{n}\,  \simeq \bigoplus_{l>0}{
\,IC_{\overline{S}_l} (L_{0,l}) } \, \bigoplus\,  T^0 \,
\, \bigoplus \oplus_{i \neq 0} \phix{i}{\fxn{f}{X}{n}}[-i],
$$
where
$T^0$ and $\phix{i}{\fxn{f}{X}{n}}$ are  skyscraper sheaves at $y.$
It follows that the natural map $H^{n+b}(X) \to H^{n+b}(f^{-1}(y))$
is surjective for every $b \geq 0,$
whence the purity statement.
Since  $\csix{j}{IC_{\overline{S}_l} (L_{0,l})}=0$
for $j \geq 0$ (cf. \ref{iccso}),
$H^n(f^{-1}(y)) =T^0$
and $H^{n+b}(f^{-1}(y))=\phix{b}{\fxn{f}{X}{n}}$ for $b >0.$
It follows that, for $b\geq 0,$ we have $H^{n+b}(f^{-1}(y))=
H^{n+b}_{\leq b}(f^{-1}(y))= H^{n+b}_{b}(f^{-1}(y)).$
The result follows from \ref{rcffv}  and \ref{nhrbr}. 
\blacksquare

\begin{rmk}
\label{nonsufficit}
{\rm If, for example,  $n=4$ and $f^{-1}(y)= 
\cup{S_j}$ is a configuration of surfaces, then
the matrix $[S_j] \cdot [S_k]$ is positive definite.
The example of $\pn{2} \subseteq \pn{4}$ shows that
the necessary Contractibility conditions expressed
by Corollary \ref{solido} are clearly not sufficient.}
\end{rmk}

In  \ci{demigsemi} we study
projective semismall (cf. \ref{defr})   maps
$f: X \to Y$ with $X$ a projective manifold.
We show that the intersection form $I$ 
associated with a  component of a codimension $2h$ 
relevant  stratum  has a precise signature: $(-1)^{h}I >0.$
The methods of \ci{demigsemi}  do not cover the case
of semismall projective maps from a quasi-projective
manifold $X.$ One may compactify the morphism,
but the condition of semismallness could be destroyed at the boundary.
The methods of this paper by-pass this problem
and we have 
the following
\begin{cor}
    \label{sigsemi}
    ({\bf Signature Theorem for Semismall Maps})
    Let $f: X \to Y$ be a projective and surjective
    semismall map, with $X$ nonsingular and quasi-projective.
    Then the intersection form  $I$ associated with a
    component
    of a codimension $2h$ relevant stratum is $(-1)^{h}I>0.$
    \end{cor}
{\em Proof.}
As in the proof of \ref{solido}, we may also assume that
$X$ is projective.
By slicing with hyperplane sections on $Y$ (cf. 
section \ref{sios}
and \ci{demigsemi}), we are reduced to the case when $\dim{X}=h,$ i.e.
 to the zero-dimensional stratum $S_{0}. $ Let $y\in S_{0}.$
 We have a decomposition similar to 
the one in the proof of \ref{solido}, where $T^0\simeq
H_{n}(f^{-1}(y))$ is polarized up to sign
by the intersection form on $X,$ i.e. by $I.$
\blacksquare 

\subsection{Purity and  Hodge-Lefschetz  for Intersection Cohomology}
\label{pp}
Let $H^{j}_{i,l}(X): = {\Bbb H}^{j-i -n}( Y , IC_{{\overline{S}}_{l}} 
(L_{i,l}   ) ) \subseteq H^{j}_{i}(X),$ 
$ 0 \leq l \leq \dim{Y},$ $i \in  \zed,$ $j \geq 0.$
See \ref{tm1}.c.

\n
In general, the pure-dimensional stratum $S_{l}$ is not connected
and one may write
\begin{equation}
    \label{dirsums}
 H^{j}_{i,l}(X) \, = \, \bigoplus_{S}{ H^{j}_{i, l ,S}(X)  }
    \end{equation}
where the  direct sum is over the connected components $S$
of $S_{l}.$

\medskip

A priori, it is  not clear that the l.h.s. of 
(\ref{dirsums}) is a Hodge sub-structure
of $H^{j}_{i}(X)$ and that (\ref{dirsums}) itself is a 
decomposition into Hodge sub-structures. We prove these
facts
in the following

\begin{tm}
    \label{tmp1} ({\bf Purity Theorem})
Let $j \geq 0$ and $ i \in \zed.$ 
The canonical direct sum decomposition
\begin{equation}
    \label{maineq}
H^{j}_{i}(X) \, = \, \bigoplus_{l}{   H^{j}_{i,l}(X) }
\, = \, \bigoplus_{l,S}{   H^{j}_{i,l,S}(X) }
\end{equation}
is by Hodge sub-structures,  it is $S^{\eta L}$-orthogonal,
and it 
is compatible with the $(\eta, L)-$De\-com\-po\-si\-ti\-on 
and its polarization.
\end{tm}

\begin{rmk}
    \label{tuttoht}
    {\rm
The spaces $H^{n+b+d}_{b,l,S} (X)$ behave like the cohomology
of a collection of projective manifolds.
Fix $b \in \zed,$ $l \geq 0$ and $S$ a connected component
of $S_{l},$ and let
$d \in \zed$ vary.
We get a structure which satisfies Hodge-Lefschetz properties
for $L$ analogous to the ones stated in Theorem \ref{pic}.
Similarly,  if we fix $d,$ $l$ and $S,$ and we let $b \in \zed$ vary,
then we get similar properties with respect to $\eta.$
See \ref{example3fold}, \ref{examplefamily}.
    }
    \end{rmk}

    The following is the intersection cohomology counterpart
    of the classical Hodge Theory of projective manifolds.
    
\begin{tm}
    \label{pic}
    ({\bf Hodge-Lefschetz for Intersection Cohomology})
$\;$

\n
(a) ({\bf Purity of Intersection Cohomology})
For every $j \geq 0,$ the intersection cohomology
group
$I\!H^{j}(Y,\rat)$ carries
a  weight$-j$ pure Hodge structure. 

\n
This structure is characterized by 
the property that, given any
 projective resolution $f:X \to Y,$ the resulting natural inclusion
 $$
 I\!H^{j}(Y,\rat)\,  \lorw \, H^{j}_{0}(X) 
 $$
is a map of  weight$-j$ pure Hodge structures.

\n
Moreover,
given a diagram of projective resolutions 
$$
\xymatrix{
    X \ar[r]^{g''} \ar[d]^{g'} \ar[rd]^h & X''\ar[d]^{f''} \\
      X' \ar[r]^{f'} & Y}
$$
the images of the natural inclusions
$$
I\!H^{j}(Y,\rat) \stackrel{f'\circ g'}\lorw  H^{j}_{0}(X)  
\stackrel{f''\circ g''}\longleftarrow I\!H^{j}(Y,\rat).
$$
coincide.

\n
(b) ({\bf Weak Lefschetz Theorem for Intersection Cohomology})
Let $r: Y_{1} \subseteq Y$ be the inclusion
of a subvariety hyperplane section of $Y$
transversal to all the strata of a stratification of $Y.$ Then
$$
r^{*} \, : \, I\!H^{j}(Y,\rat) \lorw I\!H^{j}(Y_{1},\rat)
$$
is an isomorphism for $j \leq \dim{Y} -2$ and injective for
$j= \dim{Y}-1.$

\n
(c) ({\bf Hard Lefschetz Theorem for Intersection Cohomology)}
The cup product map 
$$
A^{j}\,:\, I\!H^{\dim{Y}-j} (Y,\rat) \lorw  I\!H^{\dim{Y}+j}(Y,\rat)
$$
is an isomorphism for every $j \geq 0.$ Setting $P^{n-r}_{A}:=
\ke {\, A^{r+1} } \subseteq I\!H^{n-r}(Y,\rat),$ $r \geq 0$ there is
a primitive Lefschetz decomposition
\begin{equation}
    \label{ldic}
I\!H^{\dim{Y}-j}(Y,\rat) \, = \, \bigoplus_{r\geq 0}{ A^{r} 
\,P_{A}^{\dim{Y}-j-2r}}.
\end{equation}

\n
(d) ({\bf Hodge-Riemann Bilinear Relations for Intersection 
Cohomology})
Given a projective resolution $f:X \to Y$ of the singularities
of $Y,$ the forms $S_{0,-j}^{\eta L}$ polarize, up to sign,
the spaces $P_{A}^{n-j}.$
Moreover, the decomposition (\ref{ldic}) is $S^{\eta 
L}_{0,-j}$-orthogonal.
\end{tm}

\begin{rmk}
\label{comp}
{\rm
\medskip
The Hodge-Lefschetz Theorem \ref{pic} is due to several authors.
The Weak Lefschetz Theorem is due to Goresky-MacPherson \ci{g-m}.
The Hard Lefschetz Theorem for Intersection Cohomology is 
proved in \ci{bbd} using  algebraic geometry in positive characteristic.
We prove it as a corollary
to the Hard Lefschetz Theorem on Perverse Cohomology Groups.
The fact that the  Intersection Cohomology  of a projective variety
carries a canonical pure  Hodge structure is proved by M. Saito
in \ci{samhp}. The  relation with the cohomology of  resolutions is 
stated, but not proved in \ci{sa}.  
Polarizations associated with  Intersection Cohomology
appear implicitly in Saito's work.
The methods of the aforementioned papers are 
completely different from the ones of this paper.
}
\end{rmk}

\subsection{The algebraic case}
\label{tacss}
In this section we point out  that a series of simple reductions
using Hironaka's resolution of singularities
and Chow's Lemma
allows to prove Theorem
\ref{tm1} for proper  maps of algebraic varieties.
Since no new idea is necessary
for this purpose, we omit the proofs. For details see
\ci{decmightamv1}, p.71-74.

\begin{tm}
\label{tm1psc}
Let $f: X \to Y$ be a proper map of projective varieties.
Then all the  results in \ref{tpcss} hold
if we replace $\rat_{X}[n]$ with $IC_{X}.$
\end{tm}

\begin{rmk}
    \label{nontutto}
    {\rm In the singular case above, the refined intersection form 
must be replaced by the maps
${\Bbb H}^{*}(y, \alpha_{y}^{!} f_{*} IC_{X}) \to
{\Bbb H}^{*}(y, \alpha_{y}^{*} f_{*} IC_{X}),$
where $\alpha_{y}: y \to Y$ is a point in $Y.$}
    \end{rmk}
    
\begin{tm}
    \label{tm1ac}
  Let $f: X \to Y$ be a proper map of algebraic varieties.
  Then the Decomposition Theorem \ref{tm1}.b and the Semisimplicity
  Theorem \ref{tm1}.c hold for $f$ if we replace $\rat_X[n]$ by $IC_X.$
  
  \n
If, in addition,  $f$ is projective and $\eta$ is $f-ample,$
then the Relative Hard-Lefschetz Theorem \ref{tm1}.a holds as well.
\end{tm}

\
The results hold, with obvious modifications also
left to the reader, for  the push-forward 
of any complex
$K \simeq \dsdix{i}{K}$
such that each
$\phix{i}{K} \simeq \oplus{
IC_{\overline{Z}}(L_{Z})},$ where  the $Z$ are 
nonsingular locally closed  subvarieties of $X$ and the $L_{Z}$
are  self-dual local systems arising 
as direct summands of some $\phix{b}{\fxn{g}{Z'}{\dim{Z'}}}$
for some algebraic proper  $g: Z' \to Z.$

\subsection {Example: Resolution of singularities of a threefold}
\label{example3fold}
Let $f:X \to Y$  be a generically finite
and surjective map
from a nonsingular three-dimensional projective variety $X.$
For ease of exposition, we assume that $f$ admits a stratification
$Y=U \amalg C \amalg y$ where $y$ is a point in the closure of
 the smooth curve $C$ with the properties : a) $\dim{f^{-1}(y)}=2,$
b) $\dim{f^{-1}(c)} =1,$ $c \in C.$ The map $f$ is locally 
topologically trivial, when restricted to the strata.

Note that $f$ is semismall over $Y \setminus y.$
Let ${\cal L}_{U}$ be the local system on $U$ associated
with the topological covering $f^{-1}(U) \to U$ and ${\cal L}_{C}$
be the local system on $C$ associated with 
$H^{2}(f^{-1}(c)),$ $c \in C.$ Let $D:= f^{-1}(y),$
$D_{j} \subseteq D$ be the irreducible surface components.
We have that $H_{4}(D)$ is spanned precisely by the fundamental
classes of $\{D_{j} \}.$ Denote by $[D_{j}] \in H^{2}(X)$
the image of $\{D_{j} \}$ via the natural class map $cl:
H_{4}(D) \to H^{2}(X).$  Let $X^{1} \subseteq X$ be a general
$\eta-$hyperplane section of $X$  and $X_{1}$ be the zero locus
of a general section of $L=f^{*}A.$
The map $f_{|}: X^{1} \to Y$ is generically finite onto
a surface and contracts to $y$ the irreducible curves
$E_{j}: = D_{j} \cap X^{1}.$
The map $f_{|}: X_{1} \to Y$ has analogous properties,
but contracts the irreducible curves inside $f^{-1}(c),$
$c \in C \cap f(X_{1}).$

\smallskip
The defect of semismallness $r(f)=,1$ so that 
$ \phix{i}{ \fxn{f}{X}{3} }=0,$   $i\neq -1,$ $0,$ $1$ and 
$$
\begin{tabular}{||c|c|c||}  \hline
$\phix{-1}{ \fxn{f}{X}{3} }$ & $\phix{0}{ \fxn{f}{X}{3} }$ &
$\phix{1}{ \fxn{f}{X}{3} }$ \\ \hline
$H_{4}(D)_{y}$ & $IC_{Y}( {\cal L}_{U} ) \oplus IC_{\overline{C}}( {\cal 
L}_{C} ) \oplus H_{3}(D)_{y}$ & $ H^{4}(D)_{y}$ \\ \hline 
\end{tabular}
$$
If $f^{-1}(c)$ is irreducible, then $IC_{\overline{C}}( {\cal 
L}_{C} ) \simeq \nu_{*} \rat_{\widetilde{C}} [1],$
where $\nu: \widetilde{C} \to \overline{C}$ is the normalization.

\n
The Relative Hard Lefschetz Theorem, i.e.
$\eta: \phix{-1}{ \fxn{f}{X}{3} } \simeq \phix{1}{ \fxn{f}{X}{3} },$
becomes the statement that  the composition
$H_{4}(D) \to   H^{2}(X) \stackrel{\eta}\to H^{4}(X) \to H^{4}(D)$
is an isomorphism so that the ensuing bilinear form on
$H_{4}(D)$ 
\begin{equation}
    \label{chebello}
 \langle \{D_{j} \}, \{ D_{k} \}  \rangle \, = \, 
 \int_{X}{ \eta \wedge [D_{j}] \wedge [D_{k}]   }   \, = \, 
 \int_{X^{1}}{ [E_{j}] \wedge [E_{k}]   }
\end{equation}
is nondegenerate.
We have $\phix{-i}{ \fxn{f}{X}{3} } = {\cal P}^{-i}_{\eta},$ $i=0,\,1.$

\n
The Decomposition Theorem reads:
$$
\fxn{f}{X}{n} \, \simeq \, 
\left( H_{4}(D)_{y} [1] \right)  \oplus
\left( 
IC_{Y}( {\cal L}_{U} ) \oplus IC_{\overline{C}}({\cal L}_{C}) 
\oplus H_{3}(D)_{y} [0]
\right)
\oplus
\left(
H^{4}(D)_{y}[-1]
\right).
$$
The Semisimplicity Theorem implies the elementary fact that
${\cal L}_{U}$ and
${\cal L}_{C}$ are semisimple.

\n
The Hodge Structure Theorem gives the following table of Hodge 
(sub-)structures:
$$
\begin{tabular}{||c||c|c|c|c|c|c|c||}  \hline\
coh.deg. & $0$ & $1$ & $2$ & $3$ & $4$ & $5$ & $6$  \\ \hline\hline
$H^{*}_{\leq -1}(X)$ & $0$ & $0$ & $\im{\, \{ H_{4}(D) \to H^{2} \}  }$ 
& $0$ & $0$ & $0$ & $0$ \\ \hline 
$H^{*}_{\leq 0}(X)$ & $H^0$ & $H^1$ & $H^2$ 
& $H^3$ & $\ke{\, \{ H^{4} \to H^{4}(D) \}  }$ & $H^5$ & $H^6$  \\ \hline
$H^{*}_{\leq 1}(X)$ & $H^0$ & $H^1$ & $H^2$ 
& $H^3$ & $H^4$ & $H^5$ & $H^6$  \\ \hline
\end{tabular}
$$
The forms $S^{\eta L}_{ij}$ are defined via $\int_{X}{ \eta^i \wedge
L^j \wedge a \wedge b } = 
\int_{ X^{i}_{j} }{ a_{|X^{i}_{j} } \wedge b_{|X^{i}_{j}} }.$

\n
The $(\eta,L)-$Decomposition Theorem and the Generalized Hodge-Riemann
Bilinear Relations give the following table of
$S^{\eta L}_{ij}-$orthogonal decompositions of the pure
Hodge structures
$H^{3-i-j}_{-i}(X).$ Each term $P^{-j}_{-i}$ is a Hodge sub-structure
polarized by $S^{\eta L}_{ij}.$
$$
\begin{tabular}{||c||c|c|c|c|c|c|c||} \hline
coh.deg. & $0$ & $1$ & $2$ & $3$ & $4$ & $5$ & $6$  \\ \hline\hline
$H^{*}_{ -1}(X)$ & $0$ & $0$ & $P^{0}_{-1} \simeq H_{4}(D)$ 
& $0$ & $0$ & $0$ & $0$ \\ \hline 
$H^{*}_{ 0}(X)$ & $P^{-3}_{0}$ & $P^{-2}_{0} $ &
$ P^{-1}_{0} \oplus L P^{-3}_{0}$ &
$   P^{0}_{0} \oplus L P^{-2}_{0}      $ 
& $L P^{-1}_{0} \oplus L^{2} P^{-3}_{0}$ 
& $ L^{2} P^{-2}_{0}$ &   $L^{3} P^{-3}_{0}$  \\ \hline
$H^{*}_{ 1}(X)$ & $0$ & $0$ & $0$ 
& $0$ & $\eta P^{0}_{-1} \simeq H^{4}(D)$ & $0$ & $0$  \\ \hline
\end{tabular}
$$
In general, the display above has several rows. Each row
presents a horizontal shifted symmetry with respect to the action of $L.$
There is 
an additional
diagonal symmetry of the display
due to $\eta.$ 

\smallskip
\n
{\em Each  horizontal and diagonal row behaves like the
cohomology of a projective manifold under the action of an ample
line bundle}.

\n
These symmetries are explained by the Hard Lefschetz Theorem for 
Perverse Cohomology Groups which gives the following display of 
isomorphisms
$$
\xymatrix{
         &           &       H_4(D) \ar@/_2pc/ 
@{-->}[ddrr]_>>>>>>>>{c_1(\eta)}    &           
&                                 &               &          \\
H^0  \ar@/_2pc/[rrrrrr]^{c_1(L)^3}  & H^1 
\ar@/^3pc/[rrrr]^{c_1(L)^2}   &  H^2/H_4(D)  \ar@/^1pc/[rr]^{c_1(L)} 
&   H^3  &  \ke \{ H^4 \to H^4(D) \}    &   H^5      &   H^6   \\
         &           &                  &           &     
H^4(D)                      &                &   
}
$$
The Purity Theorem implies that
the various pieces decompose further
according to strata into polarized
Hodge sub-structures in a fashion
compatible
with the $(\eta, L)-$Decomposition (cf. \ref{tuttoht}).

\n
By the Generalized  Grauert Contractibility Criterion,
the fibers of $f$ contribute to the pieces
$P^{0}_{-1},$ $P^{-1}_{0}$ and $P^{0}_{0}$ as we now explain.

\n
a) The class map $H_{4,-1}(D) = H_{4}(D) \to H^{2}_{-1}(X)$ is injective.
The image is the Hodge structure $P^{0}_{-1}$ polarized
by $S^{\eta L}_{10}$ which in turn coincides with
the negative-definite (\ref{chebello}). This is Grauert's
Criterion for the contraction of the curves $E_{j}$ on $X^{1}.$

\n
b) The group $H_{2,1}(f^{-1}(c)) =\{0 \}$ and the injectivity
statement \ref{nhrbr} is trivial for $c \in C.$
The map $ H_{2,0}(f^{-1}(c)) =H_{2}(f^{-1}(c)) \to H^{4}_{0}(X)$
lands in $LP^{-1}_{0}.$  This map is not the zero map:
its image is isomorphic to the invariants
$H_{2}(f^{-1}(c))^{\pi_{1}(C,c)}.$ The fact that
$S^{\eta L}_{0,-1}$ is a polarization merely reflects
the Grauert criterion for the contraction of the curves in
$f^{-1}(c),$ $c \in C \cap f(X_{1}),$ given by the map $f_{|X_{1}}:
X_{1} \to Y.$

\n
c) The class map $H_{3,0}(D) = H_{3}(D) \to H^{3}(X)$ is injective.
The image is a split  Hodge sub-structure of $P^{0}_{0}.$  Both
Hodge structures are polarized by $S^{\eta L}_{00},$
i.e. by $\int_{X}.$  In fact, by the
 Contractibility Criterion  \ref{solido}: if $D$ is to appear
in such a morphism $f,$  then its a-priori mixed Hodge structures
$H^{i}(D)$ must be pure for every  $i\geq 3.$

\subsection {Example: Families of Varieties}
\label{examplefamily}
Let $f: X \to C$ be a surjective map of projective manifolds,
$\dim{X}=m+1,$ $\dim{C} =1.$ 

\n
There is a stratification of $f$ given by $C = U\amalg S,$ where
$S$ is a finite set. Denote by $g: f^{-1}(U)=:U' \to U$ 
and  by $\beta: U \to C \leftarrow S: \alpha$ the resulting 
maps.
If ${\cal L}$ is a local system on $U,$ then $IC_{C}({\cal L})=
(\beta^{0}_{*} {\cal L} ) [1],$ where we denote by
$\beta^{0}_{*}$ the sheaf-theoretic direct image
(i.e. not the derived functor).
We have $r(f)=m$ and 
$$
\phix{j}{\fxn{f}{X}{m+1}} \, \simeq \,
( R^{m+j} g_{*} \rat_{U'} ) [1] \oplus K^{j}, \;\; j \in 
[-m, m], \qquad 0 \; \; \hbox{otherwise,}
$$
where the $K^{j}$ are sheaves supported on $S.$

\n
The Relative Hard Lefschetz Theorem translates into
$$
\eta^{j}\, :\, R^{m-j} g_{*} \rat_{U'} \, \simeq \, 
 R^{m+j} g_{*} \rat_{U'}, \qquad \qquad
 \eta^{j}\, :\,  K^{-j}\,  \simeq\,  K^{j} .
 $$
 The first isomorphism is the classical Hard Lefschetz Theorem
 for the fibers
 of the smooth map $g.$

\n
The  Decomposition Theorem reads
$$
\fxn{f}{X}{m+1} \, \simeq  \,
\bigoplus_{j}{ 
\left(
( \beta^{0}_{*} R^{m+j} g_{*}\rat_{U'} ) [1] [-j] \oplus \alpha_{*}
K^{j}[-j]
\right).
}
$$
The Semisimplicity Theorem gives the well-known
semisimplicity of the local systems
$R^{m+j} g_{*} \rat_{U'}.$ We omit drawing tables
as in \ref{example3fold}.  We point out that the stalks
$K^{j}_{p},$ $p \in S,$ are split Hodge sub-structures
of $P^{0}_{-j},$ polarized up to sign by the forms
$S^{\eta L}_{j0}.$ 
By taking cohomology sheaves we find isomorphisms
$R^{i} f_{*} \rat_{X} \simeq \beta^{0}_{*} \beta^{*}
R^{i} g_{*} \rat_{U'} \oplus \alpha_{*} K^{i-m}.$ It follows that
the natural adjunction map $R^{i} f_{*} \rat_{X}
\to \beta^{0}_{*} \beta^{*}
R^{i} g_{*} \rat_{U'}$ is surjective. By taking stalks at
$p \in S$ we get that
$$
H^{i}( f^{-1}(p) )  \lorw H^{i} ( g^{-1}(u) )^{\zed}, \qquad u \in U
$$
is surjective, i.e. that the classes in $H^{i}(g^{-1}(u) )$
which are invariant under the local monodromy
around $p$ come from $H^{i}(f^{-1}(p))$ or, equivalently, 
from  $H^{i}( f^{-1}(U_{p}) ),$ $U_{p}$ a small
Euclidean neighborhood of $p.$
This statement is known as the {\em Local Invariant Cycle
Theorem.} 
We note that, compared with the sharp versions of this theorem, due 
to 
various authors, see for instance
\ci{steenbrink}, \ci{clemens}, \ci{elzein} and \ci{navarro}, this 
proof works only for 
projective  (as opposed to  K\"ahler) families over a quasi 
projective base (as opposed to over the disk).

\subsection{The structure of the proof}
\label{ssotp}
The  set-up is as in \ref{tpcss}.
The proof of the results in \ref{tpcss}
is by a double induction on the defect of semismallness and 
on the dimension of the target of the map $f:X \to Y.$

\n
The Purity Theorem \ref{tmp1} is  proved in \ref{rtoho}.
The Hodge-Lefschetz Theorem \ref{pic} is proved in \ref{poftmpic}. 

\medskip
The starting point of the proof by induction is the following
\begin{fact}
\label{fatto}
{\rm
If  $\dim{f(X)}=0,$ then $L$ and the perverse filtration
are  trivial 
and all the results of \ref{tpcss} are either trivial
or hold  by classical Hodge 
theory.
See Theorem \ref{chl}. 
}
\end{fact}
The inductive hypothesis takes the following form
\begin{ass}
\label{basicass}
{\rm  Let $R\geq 0$ and $m>0.$
Assume that the results of \ref{tpcss}
%Theorem \ref{tm1} (the Relative Hard Lefschetz,
%the Decomposition and the Semisimplicity theorems),
%Theorem \ref{tm3} (the Hard Lefschetz Theorem for
%Perverse Cohomology), Theorem
%\ref{uf}  (the Hodge Structure Theorem),
%Corollary \ref{etaldecompo} (the $(\eta, L)-$Decom\-po\-si\-ti\-on
%Theorem)
%and Theorem \ref{tmboh} (the Polarization Theorem)
hold for every 
projective map $g:Z' \to Z$ of projective varieties with $Z'$ 
nonsingular
such that either $r(g) < R$, or $\dim f(Z) < m $ {\em and} $r(g) \leq 
R.$
}
\end{ass}

\smallskip
We  prove that if Assumption \ref{basicass} holds, then 
 the results of
\ref{tpcss}   hold 
 for every map $f:X \to Y$ as in \ref{tpcss}
 with $r(f)\leq R$ and $\dim f(X) \leq m.$

\n
In view of Fact \ref{fatto}, {\em all} the results
in \ref{tpcss}
  follow  by induction.

\medskip
What follows is an outline of the structure   of the proof
of the results of this paper.

\begin{rmk}
    \label{ovviamento}
    {\rm
    We assume \ref{basicass} and  prove
    the results in \ref{tpcss}  
    for the map $f.$ 
   Once a result has been established for $f$ we   use
   it in the proof of the results that follow.
    }
    \end{rmk}

\begin{itemize}
    
    \item
    {\bf Step 1.} We prove the Relative Hard Lefschetz Theorem
    \ref{tm1}.a  in Proposition \ref{rhl}. 
    Either the defect of semismallness is zero and there is nothing 
    to prove, or we consider the universal hyperplane section
    morphism $g$ (cf. \ref{hyperplane}).
    In this case, the defect of semismallness
    $r(g) < r(f)$ 
    and the inductive hypothesis \ref{basicass} apply. The rest of the proof
    follows a  classical path: one is reduced to the case $i=1$ and concludes
    by using  the inductive  semisimplicity statement for $g.$
    As is well-known (cf. \ci{shockwave}),
    the Decomposition Theorem \ref{tm1}.b is a formal
    homological consequence
    of the Relative Hard Lefschetz Theorem \ref{tm1}.a. 
    The Weak-Lefschetz-type result
    \ref{wltphc} implies  the Semisimplicty
    Theorem \ref{tm1}.c  for
    $i \neq 0.$ The critical case  $i=0$
     is proved  at the end of the inductive procedure. See Step 6.
    
    \item
    {\bf Step 2.} We  prove the Hard Lefschetz Theorem 
    for Perverse Cohomology Groups \ref{tm3}  for  in \ref{phltpcq}.
    The statement for $\eta$ follows  from
    the Relative Hard Lefschetz after taking hypercohomology.
    Using  hyperplane sections $X^1$ on $X$ 
    (cf. \ref{wlmech})
    and $Y_1$ on $Y$ (cf. \ref{pwl}), the  statement for $L$ is reduced to 
    checking that $L : H^{n-1}_{0}(X) \simeq H^{n+1}_{0}(X).$
     At this stage we do not know
yet that $H^{n-1}_0(X)$ has a Hodge structure, so that care is
needed when using Hodge-theoretic statements.
    We use the map $X_1= f^{-1}(Y_1) \to Y_1$ 
and the inductive Generalized Hodge-Riemann Bilinear Relations
  \ref{tmboh} to conclude.
    
   \item
   {\bf Step 3.} The Hodge Structure Theorem \ref{uf}  is proved
   in \ref{phltpcq}.
   The Hard Lefschetz Theorem for Perverse Cohomology Groups  implies that
   the filtration $W^{tot}$ (cf. \ref{asoh}) 
   on $\oplus_{j}{H^{j}(X)}$
   coincides with the weight filtration $W^{L}$ for the nilpotent operator
   $L.$  Since the latter is clearly Hodge-Theoretic, so is the 
   former, whence Theorem \ref{uf}.
   The $(\eta,L)-$Decomposition \ref{etaldecompo} is then a formal algebraic
   consequence of what has  already  been proved. 
   
   \item
   {\bf Step 4.} The Generalized Hodge-Riemann Bilinear Relations
  \ref{tmboh} 
   are proved 
   in the cases  $P^{j}_{i} \neq P^{0}_{0}$
   in 
   \ref{eccettopoo} 
   using the inductive hypothesis.
   The crucial case $P^{0}_{0}$ is proved in 
   \ref{tapwat}. 
Let $\e >0$ and    $L_{\e} := L + \e\, \eta: H^{n}(X) \to H^{n+2}(X),$
   $\Lambda_{\e}:= \ke{\, L_{\e}}  \subseteq H^{n}(X).$
   We consider the pure Hodge structure
   $\Lambda  = \lim_{\e \to 0}{ \Lambda_{\e}} \subseteq
   \ke{\, L} \subseteq H^{n}(X).$
   Since, by classical Hodge Theory (cf. \ref{chl}),
   every $\Lambda_{\e}$ is polarized (up to sign) by the
   intersection form, the space $\Lambda$ is semipolarized
   by the same form, i.e.  the relevant bilinear form 
   has a non-trivial radical, but 
   it is positive semidefinite. 
   At this point something quite remarkable happens:
   $\Lambda \subseteq H^{n}_{\leq 0}(X),$
   i.e. $\Lambda = \Lambda_{\leq 0}$ and the radical
   of the semipolarization is $\Lambda_{\leq  -1} =\Lambda
   \cap H^{n}_{\leq -1}(X).$ It follows that $\Lambda_{0}:=
   \Lambda_{\leq 0}/ \Lambda_{\leq -1} \subseteq
   H^{n}_{0}(X)$ is polarized by 
   $S^{\eta L}_{00}.$ Finally, since $P_{0}^{0} \subseteq 
   \Lambda_{0}$  is a Hodge sub-structure, it is automatically
   polarized (cf.
  \ref{important}).
   
   \item
   {\bf Step 5.} The Generalized Grauert Contractibility Criterion
   \ref{nhrbr}  is proved in \ref{iadsoicc} in the context
   of proving the semisimplicity of  $\phix{0}{\fxn{f}{X}{n}},$
   i.e. the remaining case of Theorem \ref{tm1}.c.
   By slicing the strata, we reduce the proof to the key case
   when the perversity index is zero and the stratum is 
   zero-dimensional. To deal with this case
   we use Deligne's Theory of Mixed Hodge structures to
   infer the injectivity part of the statement. The
   relevant graded piece of the homology of the fiber is then
   a Hodge sub-structure of the corresponding perverse  cohomology
   group of 
   $X,$ compatibly with the 
   $(\eta,L)-$Decomposition. 
   The nondegeneration statement follows from
    the Generalized Hodge-Riemann Bilinear Relations 
   and the elementary \ref{important}.

    \item
    {\bf Step 6.} We prove that $\phix{0}{\fxn{f}{X}{n}}$ is 
    semisimple
    in $\S$\ref{sios}. 
    We first  prove Theorem \ref{itissum} which 
    states that that the complex in question is
    a direct sum of intersection cohomology complexes
    of local systems on strata. The proof uses  
    the Splitting Criterion \ref{splitp}: the  criterion is met
   by virtue of the Generalized Grauert Contractibility Criterion.
    The Refined Intersection
   Form Theorem \ref{rcffv} follows from Theorem \ref{itissum}
   and the Splitting Criterion  \ref{splitp}.
    We prove Theorem \ref{sss}, i.e. 
    that the local systems above are semisimple
    by exhibiting them as quotients
    of local systems associated with certain auxiliary  smooth proper maps.
    
    \end{itemize}
    
    \medskip
  The Purity  Theorem \ref{tmp1}
   is proved in \ref{rtoho} by a similar induction,
   using   the results of \ref{tpcss} as well as
    Deligne's theory of Mixed Hodge structures.

\section{Notation and preliminary results}
\label{rcst}
We work over $\comp$  and denote
 rational singular cohomology groups by 
$H^{*}(-).$

\subsection{The topology and Hodge Theory of algebraic varieties}
\label{thtav}

In this section we collect classical results concerning the topology
and the  Hodge theory of projective manifolds. 

Let $l \in \zed$, $H$ be a finitely generated abelian group, 
$H_{\rat}:= H 
\otimes_{\zed} \rat$, $H_{\real}= H \otimes_{\zed}\real,$
$H_{\comp}= H \otimes_{\zed} \comp$.

\n
A pure Hodge structure  of weight $l$ on $H$,  $H_{\rat}$ or
$H_{\real},$
is a direct sum decomposition $H_{\comp}= \oplus_{p+q=l} H^{pq}$
such that $H^{pq}= \overline{ H^{qp} }.$
The Hodge  filtration  is the decreasing filtration 
$F^{p}(H_{\comp}):=
\oplus_{p'\geq p}H^{p'q'}.$ A morphism of Hodge structures
$f: H \to H'$ is a group homomorphism such that
$f \otimes Id_{\comp}$ is compatible with the Hodge filtration,
i.e. such that it is a filtered map. Such maps are automatically 
strict.
The category of Hodge structures of weight $l$ is abelian.

\n
Let $C$ be the Weil operator, i.e. $C: H_{\comp} \simeq H_{\comp}$
is such that $C(x)=i^{p-q} x$, for every $x\in H^{pq}$.
It is a real operator. Replacing $i^{p-q}$ by 
$z^{p}{\overline{z}}^{q}$
 we get a real  action $\rho$ of $\comp^{*}$ on $H_{\real}.$
A polarization of the real pure Hodge structure
$H_{\real}$ is a  real bilinear form $\Psi$ on $ H_{\real}$
which is  invariant under the action given by $\rho$ restricted to 
$S^{1} \subseteq \comp$
and such that the bilinear form
$\widetilde{\Psi} (x,y):= \Psi(x, Cy)$ is symmetric and positive definite.
If $\Psi$ 
is a polarization, then $\Psi$ is symmetric
if $l$ is even, and antisymmetric if $l$ is odd.
In any case, $ \Psi$ is nondegenerate.
In addition, for every $0 \neq x \in H^{pq}$, 
$(-1)^{l}i^{p-q}\Psi(x,\overline{x}) >0,$
where $\Psi$ also denotes 
the $\comp-$bilinear extension of $\Psi$ to $H_{\comp}.$
The following remark is used several times in this paper.
\begin{rmk}
 \label{important}
 {\rm
 If $H' \subseteq H$ is a Hodge sub-structure,
then $H_{\real}$ is fixed by  $C.$
It follows  that
 $\Psi_{|H'_{\real}}$ is a polarization, hence it is nondegenerate.
}
\end{rmk}

\smallskip
Let $X$ be a nonsingular projective
variety of dimension $n,$ $\eta$ be an ample line bundle
on $X.$ For every $r \geq 0$ define $P^{n-r} = \ke{ \,\eta^{r+1} }
\subseteq H^{n-r}(X, \rat).$ Classical Hodge Theory
states that, for every $l,$  $H^{l}(X, \zed)$ is a pure Hodge 
structure
of weight $l,$  $P^{n-r}$ is a rational pure Hodge structure 
of weight $n-r$ polarized by a modification of the Poincar\'e pairing 
on $X.$

\begin{tm}
\label{chl}
(a) ({\bf The Hard Lefschetz Theorem})
For every $r \geq 0$ one has
$$
\eta^{r}\,: \, H^{n-r}(X, \rat) \, \simeq \, H^{n+r}(X, \rat).
$$
\n
(b) ({\bf The Primitive Lefschetz Decomposition})
For every $r \geq 0$ there is the direct sum decomposition
$$
H^{n-r}(X, \rat)\, = \, \bigoplus_{j \geq 0} \eta^{j}P^{n-r-2j}
$$
where each summand
is a  pure Hodge sub-structure of weight $n-r$   and
all summands are mutually orthogonal with respect to the
bilinear form $\int_{X}{ \eta^{r} \wedge - \wedge -}.$

\n
(c) ({\bf The Hodge-Riemann Bilinear Relations})
For every $0\leq l \leq  n,$ the bilinear form
$(-1)^{\frac{l(l+1)}{2}}\int_{X}{\eta^{n-l} \wedge - \wedge -}$
is a polarization of the pure weight$-l$ Hodge structure
$P^{l} \subseteq H^{l}(X, \real).$ In particular,
$$
(-1)^{\frac{l(l-1)}{2}} i^{p-q} \int_{X}{\eta^{n-l}\wedge
\alpha \wedge \overline{\alpha} } \; >\; 0, \qquad \forall \;\, 0 \neq \alpha
\in P^{l} \cap H^{pq}(X, \comp).
$$
\end{tm}

A local system
$\cal L$ on an algebraic variety $Y$ is said to be
{\em semisimple} if every local subsystem
${\cal L}'$ of $\cal L$ admits a complement, i.e. a  local 
subsystem 
${\cal L}''$ of $\cal L$ such that ${\cal L} \simeq {\cal L}' \oplus 
{\cal L}''.$

\begin{rmk}
\label{sszos}
{\rm 
If $Y$ is normal and $Y'\subseteq Y$
is a Zariski-dense open subset, then $\cal L$ is semisimple
if and only ${\cal L}_{|Y'}$ is semisimple. In fact,
the natural map 
$\pi_{1}(Y',y') \lorw \pi_{1}(Y,y')$ is surjective
for any $y' \in Y'.$
}
\end{rmk}

\begin{tm}
\label{dss}
({\bf Decomposition Theorem for proper smooth maps})
Let $f: X^{n} \to Y^{m}$ be a smooth proper 
map of smooth algebraic
varieties of the indicated dimensions and $\eta$ be an ample
line bundle on $X.$ 
Then 
$$
\eta^{i}\, : \, R^{n-m-i} f_{*} \rat_{X}\,  \simeq\,  R^{n-m+i} f_{*}
\rat_{X}, \; \forall i \geq 0, \qquad Rf_{*}\rat_{X} \, \simeq \,
\bigoplus_{i\geq 0} R^{i}f_{*}\rat_{X} [-i]
$$
and
the local systems $R^{j}f_{*}\rat_{X}$
are semisimple on $Y.$
\end{tm}
{\em Proof.} See \ci{dess} and  \ci{ho2}, Th\'eor\`eme 4.2.6.
\blacksquare

\bigskip
\begin{tm}
\label{mhssc} ({\bf Mixed Hodge structure on cohomology})
Let $X$ be an algebraic variety. For each $j$ there is an increasing
weight filtration
$$
\{0 \} \, =\,  W_{-1} \, \subseteq W_{0} \, \subseteq \, \ldots 
\, \subseteq W_{2j} \, = \, H^{j}(X, \rat)
$$
and a decreasing Hodge filtration
$$
H^{j}(X, \comp) \, = \, F^{0} \, \supseteq \, F^{1} \, \supseteq 
\ldots \, \supseteq\,  F^{m}\, \supseteq \, F^{m+1} \, = \, \{0 \}
$$
such that the filtration induced by $F^{\bullet}$ on the 
complexified graded pieces
of the weight filtration endows 
every graded piece $W_{l} / W_{l-1}$ with a pure Hodge structure of 
weight $l.$

\n
This structure is functorial for maps of algebraic
varieties and the induced maps strictly preserve
both filtrations.
\end{tm}

We shall need the following two properties
of this structure. See \ci{ho3}, 8.2 for more.

\begin{tm}
\label{propofmhs}
Let $Z \subseteq U \subseteq X$ be embeddings
where $U$ is a Zariski-dense  op\-en subvariety of
the nonsingular variety $X$ and $Z$ is a closed subscheme
of $X.$ Then
$$
\im{ \left( \, H^{j}(X) \lorw H^{j}(Z) \right) } \, =
\, \im{ \left(\, H^{j}(U) \lorw 
H^{j}(Z) \right) }, \quad 
\forall \, j \geq 0.
$$
Let $g: T \to Z$ be a proper algebraic map
of proper schemes,   $T$ nonsingular. Then
$$
\ke{ \left( \, g^{*}: H^{j}(Z)  \lorw  H^{j}(T )\right) }  \, = \, W_{j-1}
\left( H^j(Z) \right), \quad
\forall \, j \geq 0.
$$
 \end{tm}

\subsection{Whitney stratifications of algebraic maps}
\label{wsoav}
It is known that every algebraic variety $Y$
of dimension $d$ admits a
Whitney stratification ${\frak Y}$
where the  strata  are locally closed algebraic subsets
with a {\em finite} number of irreducible nonsingular components.
See \ci{borel}, I.1, I.4,
\ci{g-m}, I, and the references contained therein.
In particular, $Y$ admits a filtration
$Y= Y_{d} \supseteq Y_{d-1} \supseteq Y_{d-2} \supseteq 
\ldots \supseteq Y_{1} \supseteq Y_{0} \supseteq Y_{-1}=
\emptyset$ by closed algebraic subsets subject to the following 
properties.

\medskip
\n
(1) $S_{l} := Y_{l} \setminus Y_{l-1}$ is either
empty or a locally closed algebraic subset of pure
dimension $l;$ the connected components
of $S_{l}$ 
are a finite number of  nonsingular algebraic varieties.
We  have Zariski-dense open sets 
$U_{l}: =Y \setminus Y_{l-1} = \amalg_{l'\geq l}{S_{l'}},$
such that
$U_{l} = U_{l+1} \amalg S_{l}.$
Note that $U_{d}$ is a nonsingular Zariski-dense open subset of $Y$ 
and that $U_0=Y.$

\medskip
\n
(2) (Local normal triviality) 
Let
$y \in S_{l}$, $\overline{N}$ be a normal slice  through $S_{l}$ at $y$, 
$\frak L$
be the link of $S$ at $y$, $N: = \overline{N} \setminus {\frak L}$ be the 
(open) normal slice. The spaces
$\overline{N}$, ${\frak L}$ and $N$  inherit  Whitney stratifications.  
$\overline{N}$ ($N$, resp.)
is homeomorphic in a stratum-preserving manner to the cone $c({\frak L})$
($c({\frak L}) \setminus {\frak L}$, resp.) over the 
link $\frak L$ with vertex identified to $y$. 
The cone is stratified using the cone structure and the 
given stratification of the link.
The point $y$ admits an open euclidean  neighborhood $W$  in $Z$ which 
is 
homeomorphic in a stratum-preserving manner to 
$\comp^{l } \times N$.

%One can shrink $W$ in the two directions
%of the product. This gives rise to the notion
%of {\em standard neighborhoods}, with respect to
%${\frak Y}$, for the points of $Y$.

\begin{defi}
    \label{defofstr}
{\rm ({\bf Stratification of $Y$}) In this paper,  
the term {\em stratification
of $Y$}
indicates a finite, algebraic  Whitney stratification of $Y.$
The resulting open and closed embeddings are denoted
 by $ S_l \stackrel{\alpha_l}\lorw U_l 
\stackrel{\beta_l}\longleftarrow U_{l+1}.$ 
}
\end{defi}

\begin{rmk}
\label{cutstr}
{\rm
Let $Y \subseteq \pn{N}$ be a a quasi-projective variety,
 ${\frak Y}$ a stratification
of $Y.$ 
Bertini Theorem implies that, for every $l>0$ for which
$S_{l}$ is not empty, the 
normal slice $N$ through a point $y \in S_{l}$ can be chosen to 
be the trace,
in a suitable euclidean neighborhood of $y$ in $Y,$ of a
complete intersection of $l$ hyperplanes  of $\pn{N}$
passing through $y$, transversal to all strata of ${\frak Y}$.
}
\end{rmk}
The Thom Isotopy Lemmas, adapted to the algebraic setting,
yield the following result. See \ci{g-m}, I.7.

\begin{tm}
\label{tila}
Let $f: X \to Y$ be an algebraic map of algebraic varieties.
There exist finite algebraic  Whitney stratifications
${\frak X}$ of $X$ and ${\frak Y}$ of $Y$ such that,
given any connected component $S$ of a ${\frak Y}$ stratum
$S_{l}$ on $Y:$

\n
1) $f^{-1} (S)$ is an union of connected components of strata
of ${\frak X}$ each of which is mapping submersively to $S$;
in particular, every fiber $f^{-1}(y)$ is stratified  by
its intersection with the strata of ${\frak X}.$

\n
2) $\forall y \in S$ there exists an euclidean open  neighborhood
$U$ of $y$ in $S$ and a stratum-preserving homeomorphism
$h: U \times f^{-1} (y) \simeq f^{-1}(U)$ such that
$f \circ h$ is the projection to $U$.
\end{tm}

\begin{defi}
\label{defofstrforf}
{\rm ({\bf Stratification for $f$})
A pair of stratifications ${\frak X}$ and $\frak Y$ as in Theorem
\ref{tila} is called a {\em stratification for f.}
}
\end{defi}

If $f$ is an open immersion,
 then a stratification ${\frak Y}$ induces one on $X$.
If $f$ is a closed immersion, one can choose
a finite Whitney stratification 
${\frak X}$ so that every stratum of it is the intersection
of $X$ with strata of ${\frak Y}$ of the same dimension. 
In either case, one obtains a stratification for $f.$

\subsection{The category $D(Y)$}
\label{stabi}
Let $Y$ be an algebraic variety   and
$D^{b}(Y)$ be the bounded derived category 
of sheaves of 
rational vector spaces on $Y$.
We refer to \ci{borel}, $\S$V and to \ci{iv} for an account of the 
formalism of
derived categories and Poincar\'e-Verdier Duality.

\begin{defi}
{\rm ({\bf Cohomologically-constructible})
 Let ${\frak Y}$ be  stratification of $Y.$
We say that {\em $K \in  Ob(D^{b}(Y))$ 
is   ${\frak Y}$-cohomologically-constructible}
(in short, ${\frak Y}-$cc)
if 
$\forall \,j \, \in \zed$ and  $\forall  \,l$, 
the sheaves $\csix{j}{K}_{| S_{l} }$ are locally constant
and the stalks are finite dimensional.
}
\end{defi}

Let $D(Y)$ be the full  sub-category of $D^{b}(Y)$ 
consisting
of those complexes $K$ which are ${\frak Y}$-cc with respect to 
{\em some} stratification ${\frak Y}$ of $Y.$

\n
The category $D(Y)$  is triangulated and it is
preserved by the truncation, Verdier Duality and $Rhom$ functors. 
The dualizing complex of $Y$ is denoted by  $\omega_{Y} \in Ob (D(Y))$
and the Verdier dual of a complex $K$ by ${\cal D}(K) = Rhom(K, \omega_{Y}).$
In fact,  ${\frak Y}-$cc complexes are stable under these 
constructions. By {\em triangle}  we mean a distinguished triangle.

\smallskip
The four functors $(Rf_{!},f^{!}, f^{*}, Rf_{*}),$ are denoted 
here simply
by $(f_{!}, f^{!}, f^{*}, f_{*} ).$

One has the following properties; see \ci{borel}, V.10.13 and 16:
if  $F$ is  ${\frak Y}-$cc and $G$ is  ${\frak X}-$cc,
then $f^{*}F,$ $f^{!}F$ are ${\frak X}-$cc
and
 $f_{*}G$ and $f_{!}G$ are ${\frak Y}-$cc.
In particular, 
$f_{!}f^{!}F,$ $f_{*} f^{*}F$ are ${\frak Y}-$cc  and
$f^{!} f_{!} G,$ $f^{*} f_{*}G$ are ${\frak X}-$cc.

\smallskip
The following facts are used in the sequel, often without
explicit mention.

\smallskip
 \n
 The pairs
 $(f_{!}, f^{!})$ and $(f^{*}, f_{*})$  are pairs of adjoint functors
 so that there are natural transformations
 $Id \to f_{*} f^{*},$ $f^{*}f_{*} \to Id,$ $f_{!}f^{!} \to Id$
 and
 $Id \to f^{!} f_{!}.$

 \n
 Let $\alpha: Z \to Y $ be the embedding
 of a closed algebraic subset,  $\beta:
 U \to Y$ be the embedding of the open complement
 and $K \in Ob(D(Y)).$
 There are natural isomorphism:
 $\alpha_{!} \simeq \alpha_{*},$  $\beta^{!} \simeq \beta^{*}$
 and
 dual   triangles:
 $$
 \alpha_{!} \alpha^{!}K \lorw K \lorw \beta_{*}\beta^{*} K
 \stackrel{[1]}\lorw, \quad
 \beta_{!} \beta^{!} K \lorw K \lorw \alpha_{*} \alpha^{*}K
 \stackrel{[1]}\lorw,
 $$
 whose the associated long exact sequences in hypercohomology
 are the ones of the pairs  
 $\ihixc{l}{Y,U}{K}$ and $\ihixc{l}{Y,Z}{K},$ respectively.

\n
There are  canonical isomorphisms $f^{!} \omega_{Y} \simeq 
\omega_{X},$
${\cal D}{\cal D} \simeq Id_{D(Y)},$ $ {\cal D} f_{*} \simeq
f_{!} {\cal D}$ and ${\cal D}f^{*} \simeq f^{!} {\cal D}.$

\n
Let  $f$  be proper and consider a Cartesian diagram 
$$
\xymatrix{
X' \ar[r]^{u'} \ar[d]^{f'} & X \ar[d]^f \\
Y' \ar[r]^u & Y.
}$$
%\begin{array}{ccc}
%X' & \stackrel{u'}\lorw  & X \\
%\downarrow f' & & \downarrow f \\
%Y' & \stackrel{u}\lorw & Y.
%\end{array}
%$$
Since $f_{!} \simeq f_{*},$  the 
{\em  Change of Coefficients Formula} reads:
$$
K \stackrel{\Bbb L}\otimes f_{*}K'
\,  \simeq \,
 f_{*} (f^{*}K \stackrel{\Bbb L}\otimes 
K'), \qquad  \forall \, K \in Ob(D(Y)), \; \forall \,K' \in Ob(D(X)) 
$$
and 
the {\em Base Change Theorem
for Proper Maps} reads: 
$$ u^{*}f_{*} \simeq {f'}_{*}{u'}^{*}, \qquad \qquad
 u^{!} f_{*} \simeq f'_{*} {u'}^{!} .
$$

\subsection{The intersection form on the fibers of the map $f:X \to 
Y$.}
\label{sperem}
In this section we introduce the intersection form
on the fibers of an algebraic map $f:X \to Y.$ We could not find a reference
serving the needs of the present paper.
Lemma \ref{collectinfo} 
and  Theorem \ref{rcffv} describe important  properties of the intersection
form.

\medskip
Let $Z$ be an algebraic set, $c: Z \to pt$ be the constant map.
We have $\omega_{Z}\simeq c^{!} \rat_{pt}.$ Define
the Borel-Moore homology groups with rational coefficients
of $Z$ as 
$$
H^{BM}_{l}(Z)\,:=\, \ihixc{-l}{Z}{\omega_{Z}}.
$$
We have that $H^{BM}_{l}(Z) \simeq H^{l}_{c}(Z)^{\vee}.$ If
$Z$ is compact, then $H^{BM}_{l}(Z) \simeq H_{l}(Z).$
Let $i: Z \to W$ be a  map of algebraic sets. If $i$ is proper,
then the natural
adjunction map, the identification
$i_{*} \simeq i_{!}$ and the isomorphism,
$\omega_{Z} \simeq i^{!} \omega_{W}$, give the map $i_{*}\omega_{Z}
\to \omega_{W}.$ The resulting maps in hypercohomology
$i_{*}:H^{BM}_{l}(Z) \lorw H^{BM}_{l}(W)$ are the usual
proper-push-forward maps.
If $i$ is an open immersion, then using the natural adjunction map
and the identification
$i^{*} \simeq i^{!},$ we get a map $\omega_{W} \to i_{*} \omega_{Z}$
whose counterparts is hypercohomology are the restriction to 
an open set
maps $H^{BM}_{l}(W) \lorw H^{BM}_{l}(Z).$

Let $y \in Y$  and $i: f^{-1}(y) \to X.$ 
Using the isomorphism $\omega_{X}[-n] \simeq \rat_{X}[n],$
we get  a  natural  sequence of maps
$$
i_{!} \omega_{f^{-1}(y)} [-n] \lorw \omega_{X}[-n]
\simeq \rat_{X}[n]
\lorw i_{*} \rat_{f^{-1}(y)} [n]
$$
where the first and third map are each other's dual.
Taking  degree $l$ hypercohomology we get maps
$$
H^{BM}_{n-l}(f^{-1}(y)) \lorw H^{BM}_{n-l}(X)
\simeq H^{n+l}(X) \lorw H^{n+l}(f^{-1}(y)).
$$
The resulting pairing 
\begin{equation}
    \label{rintpa}
H^{BM}_{n-l}(f^{-1}(y)) \times H_{n+l}(f^{-1}(y)) \lorw \rat
\end{equation}
is called the {\em refined intersection form on
$f^{-1}(y) \subseteq X.$} Note that we may replace $X$ by any 
euclidean open neighborhood  of $f^{-1}(y).$ Geometrically,
this form 
corresponds to intersecting locally finite cycles supported on 
$f^{-1}(y)$ with finite cycles of complementary
dimension in $X$ supported on $f^{-1}(y).$ 

\begin{rmk}
    \label{1111}
    {\rm  
If $y$ lies on a positive dimensional stratum $S_{l},$ $l >0,$ then
the refined intersection form is trivial: in fact, 
by the local triviality of the stratification, a cycle supported
on $f^{-1}(y)$ can be moved to a homologous one supported
on a nearby fiber $f^{-1}(y'),$ $y' \neq y,$ $y' \in S_{l}.$
See also \ref{alk}.
     }
    \end{rmk}

We now assume that $f$ is proper.
We have  the canonical identification
$H^{BM}_{l}(f^{-1}(y)) \simeq H_{l}(f^{-1}(y))$ and
the usual  base-change
identifications:
%$$
%\begin{array}{ccccccc}
%f^{-1}(y) &  \stackrel{i}\lorw & X &&&& \\
%\downarrow \Phi & \fbox{} & \downarrow f  & \qquad & 
%\alpha^{*} f_{*} \simeq \Phi_{*} i^{*}, & \quad &
%\alpha^{!} f_{*}
%\simeq \Phi_{*} i^{!}
%\\
%y & \stackrel{\alpha}\lorw & Y &&&&
%\end{array}
$$
\xymatrix{
f^{-1}(y) \ar[r]^i \ar[d]^{\Phi} &  X \ar[d] ^f \\
y  \ar[r]^{\alpha} & Y & &
\alpha^{*} f_{*} \simeq \Phi_{*} i^{*},  \quad 
\alpha^{!} f_{*}
\simeq \Phi_{*} i^{!}}
$$
giving rise to a self-dual diagram of adjunction maps:
$$
\alpha_{!} \alpha^{! }f_{*} \omega_{X}[-n]
\lorw f_{*}  \omega_{X}[-n] \simeq 
f_{*}\rat_{X}[n] \lorw 
\alpha_{*} \alpha^{*} f_{*}\rat_{X}[n]
$$
which, after taking hypercohomology,
give the refined intersection form (\ref{rintpa})
on $f^{-1}(y).$

\begin{rmk}
    \label{thepoint}
    {\rm
We shall consider  the map
$\alpha_!\alpha^!f_* \rat_X[n] \to  f_* \rat_X[n]$
in connection with the Splitting Criterion 
 \ref{splitp}. On the other hand, 
the map that arises geometrically is
$\alpha_{!} \alpha^{! }f_{*} \omega_{X}[-n]
\to f_{*}  \omega_{X}[-n].$  Using
 that $\alpha_{!}\alpha^{!} \to Id$ is a natural 
transformation of additive functors and the isomorphism
$\omega_{X}[-n] \simeq \rat_{X}[n],$ one has a commutative diagram:
$$
\xymatrix{\alpha_!\alpha^!f_* \omega_X[-n] \ar[r] \ar[d]^{\simeq} & 
f_* \omega_X[-n]  \ar[d]^{\simeq}\\
\alpha_!\alpha^!f_* \rat_X[n] \ar[r] & f_* \rat_X[n].
}
$$
%$$
%\begin{array}{ccc}
%\alpha_!\alpha^!f_* \omega_X[-n]& \lorw & f_* \omega_X[-n]  \\
%\downarrow \simeq  &  & \downarrow \simeq   \\
%\alpha_!\alpha^!f_* \rat_X[n]& \lorw & f_* \rat_X[n]. \\
%\end{array}
%$$
The two horizontal maps
are  thus equivalent and one  can check the hypotheses of the Splitting
Criterion 
\ref{splitp} on the top one.
}
\end{rmk}

\subsection{The local structure of a ${\frak Y}-$cc
complex along a stratum }
\label{locstr} 
The references are \ci{g-m}, 1.4,   \ci{borel}, V.3. and Lemma V.10.14.
Let $Y$ be a projective variety, ${\frak Y}$ be a stratification,
$y \in S \subseteq S_{l}$ be a point in a connected component
$S$ of a stratum $S_{l}$, $N$ be a normal slice  through $S$ at $y.$
Let $W$ be a standard open neighborhood of $y$  in $Y$,  
homeomorphic in a stratum-preserving manner to 
$\comp^{l } \times N$. Let $\pi:W \to N$ be the corresponding map,
$\dot{N}:= N \setminus y$, $\dot{W}: =W \setminus (S\cap W).$
We have a commutative diagram of cartesian  squares,
where $a,$ $\alpha$, $i_{y},$  and $i_{N}$
are closed immersions, $b$
and $\beta$ are open immersions, $c,$ $\pi$
and $\dot{\pi}$ are trivial
topological $\comp^{l}-$bundles, $c\circ \alpha_{y} =Id_{y}$ and
$\pi \circ i_{N} = Id_{N}$:
$$
\xymatrix{
S\cap W  \ar[r]^{\alpha} \ar[d]_c & W \ar[d]_{\pi} &  \dot{W}  \ar[d]^{\dot{\pi}} \ar[l]^{\beta}   \\
y \ar@/_/[u]_{\alpha_y} \ar[r]^{a} \ar[ur]_{i_Y} & N \ar@/_/[u]_{i_N} &  \dot{N} \ar[l]^{b} 
}
$$
%%%%%$$
%\begin{array}{ccccc}
%y  & \stackrel{a}\lorw  & N &  & \\
%\downarrow \alpha_{y} & \searrow i_{y} & \downarrow i_{N} &  & 
%\\
%S\cap W & \stackrel{\alpha}\lorw & W
%& \stackrel{\beta}\longleftarrow &   \dot{W}
%\\
%\downarrow c &  & \downarrow \pi &  & \downarrow \dot{\pi}
%\\
%y & \stackrel{a}\lorw & N & \stackrel{b}\longleftarrow & \dot{N}
%\end{array}
%$$
The following rules apply:
$\alpha_{*} \simeq \alpha_{!},$  $\alpha^{*} 
\alpha_{*} \simeq Id \simeq \alpha^{!}\alpha_{*},$ $\pi^{*} \simeq 
\pi^{!} [-2l],$ 
$c^{*}\simeq c^{!}[-2l],$ $\beta_{*} \dot{\pi}^{*} \simeq 
\pi^{*} b_{*}.$
If $K$ has locally constant cohomology sheaves on $S$, then
$\alpha_{y}^{!}K \simeq \alpha_{y}^{*} K [-2l].$

\medskip
Let $K$ be ${\frak Y}-$cc.
On $W$, we  have that: 
$\;
K \simeq \pi^{*} \pi_{*}K \simeq \pi^{*} K_{|N}.
$
That is, $K$ has, locally at a point of any stratum,  a product 
structure
along the stratum. See \ci{borel}, Lemma V.10.14.

\n
The sheaves ${\cal H}^i(\alpha^*K)$
and  ${\cal H}^i(\alpha^!K)$ on $S \cap W$
are constant
with representative stalks
${\cal H}^i(\alpha^*K)_{y} \simeq  {\Bbb H}^i(N,K_{|N})
\simeq {\cal H}^i(a^*K_{|N}),$
and
 ${\cal H}^i(\alpha^!K)_{y} \simeq  {\Bbb H}^i(N,\dot{N};K_{|N})
 \simeq {\cal 
H}^i(a^!K_{|N}).$

\n
\begin{rmk}
\label{resshperv}
{\rm
Since $\pi^{*}$ is fully faithful,
if $K$ is self-dual, then
$K_{|N}[-l]$ is self-dual. If $K$ is perverse 
(cf. \ref{perverse}) on $Y$, then, using the characterization
of perverse sheaves in Remark
\ref{refpss} and  Lemma \ref{ugus}, one shows that
$K_{|N}[-l]$ is perverse on $N.$
}
\end{rmk}

\begin{rmk}
    \label{fullyf}
    {\rm 
Let $\alpha : Y'\lorw Y$ be a closed immersion
of algebraic varieties.
We have that  $\alpha_{!}  \simeq \alpha_{*}$
are fully faithful so that  
for every  $K \in Ob(D(Y)),$ the composition
of the adjunction maps $\alpha_{!} \alpha^{!}K \to
K \to \alpha_{*} \alpha^{*} K$ yields a natural map:
$ \;
\alpha^{!}K \to \alpha^{*}K.
$
}
\end{rmk}

\begin{lm}
\label{adapt}
Let $Y$ be an algebraic variety,  
 ${\frak Y}$ be a stratification of $\,Y$,  $K$ be ${\frak Y}-$cc, 
 $\alpha:S\to Y$ be the embedding of a connected component
 of  a  stratum 
 $S_{l}$, $y \in S.$

 \n
The natural map  
 $\alpha_{!}\alpha^{!}K \to \alpha_{*}\alpha^{*}K$ 
 coincides, when restricted to a standard neighborhood $W$
 of $y$ in $Y$,
 with $c^{*}$ of the analogous map
 $a^{!} K_{|N} \lorw a^{*}K_{|N}.$ 
 
 \n
 The same is true for the induced  maps
 $\csix{j}{\alpha^{!} K} \to \csix{j}{\alpha^{*}K},$
 $\csix{j}{\alpha^{!} K}_{y} \to \csix{j}{\alpha^{*}K}_{y}$
 induced
 on the cohomology sheaves and on their stalks at $y.$
 
 \n
 Finally:
 $\;
 \csix{j}{ \alpha^{!} K}_{y}     
 \simeq \csix{-j-2l}{ i_{y}^{*} {\cal D}(K)}^{\vee}
 \simeq 
\csix{-j}{a^{*} {\cal D}( K_{|N}) }^{\vee}.
 $
 \end{lm}
{\em Proof.} The question being  local
 around $y\in Y,$ we may work on $W.$ 
We may assume that
$K= \pi^{*}K_{|N}$ and that $S$ is closed, so that,
since $\alpha_{!} \simeq \alpha_{*}$ are fully faithful,
it is enough to study the map $\alpha^{!}K \lorw \alpha^{*}K.$

\n
We have $\alpha^{!} \pi^{*}K_{|N} \simeq \alpha^{!} \pi^{!}K_{|N} 
[-2s]
\simeq c^{!} a^{!} K_{|N}[-2s] \simeq c^{*} (a^{!}K_{|N}),$
i.e. $(\alpha_{!} \alpha^{!} K)_{|S}$ is a pull-back from
$p$ and so are its cohomology sheaves.
The statement concerning $\alpha^{!} K \lorw 
\alpha^{*}K,$ the induced maps on the cohomology sheaves
and associated stalks at $y$ follow.  
The duality statements stem from   Poincar\'e-Verdier
Duality on $y$ and the isomorphism $\alpha^{*}_{y} 
\alpha^{!} K \simeq i_{y}^{!} K [2l],$ which holds in view
of the fact that $\alpha^{!}K$ has locally constant cohomology sheaves
on $S$.
\blacksquare

 \bigskip
A global counterpart of a normal slice is the notion of 
stratified normally 
nonsingular inclusion; see \ci{g-m}, Theorem I.1.11. 
The embedding  $Z \to Y$ 
of a subvariety 
is said to be a {\em normally nonsingular inclusion} if 
%$Z$ has a 
%tubular 
%neighborhood in $Y$, that is, 
there exists a neighborhood $W$ of $Z$
in $Y$
and a retraction $\pi: W \to Z$ 
which is locally homeomorphic to a  projection:
every point $z \in Z$ has a neighborhood $U\subseteq Z$ and a 
homeomorphism 
$\pi^{-1}(U) \simeq U \times \comp ^l $ 
compatible with the maps to $U.$
In addition, the homeomorphism
$\pi^{-1}(U) \simeq U \times \comp ^l $ can be chosen  
so that it is stratum-preserving
with respect to the induced stratification on $\pi^{-1}(U)$
and to  the  stratification
on  $U \times \comp ^l$ 
given by the product of the trivial 
stratification  
on  $\comp ^l $ with that induced on  $U$ by the transversality 
assumption.

A normally nonsingular inclusion can produced by  intersecting 
a projective variety $Y$ with a subvariety of the ambient projective 
space, 
e.g. a hypersurface, which intersects transversally every stratum of
a given stratification ${\frak Y}$ of $Y$
(cf. \ref{cutstr}).
The universal hyperplane section construction in
\ref{hyperplane} is an example.
A normally nonsingular inclusion carries a cohomology class.

\smallskip
The following fact is well-known and will be used often
in this paper.
\begin{lm}
\label{ugus}
Let $i:Z \to Y$ be a normally nonsingular inclusion of
complex codimension $d$ of
complex varieties, 
transversal to 
every stratum of a stratification ${\frak Y}$ of $Y,$ and 
$K$ be ${\frak Y}-$cc. Let $\pi:W \to Z$ be a retraction 
of a tubular neighborhood of $Z$ in $Y$ onto $Z.$ Then  we have
a) $K_{|W} \simeq \pi^*\pi_*(K_{|W}) \simeq \pi^*{K_{|Z}}$ and
b) $i^!K\simeq i^*K[-2d].$
\end{lm}
{\em Proof.}  We denote $K_{|W}$ simply by $K.$

\n
(a) By virtue of the local triviality assumption, 
the natural adjunction map
$ \pi^{*} \pi_{*} (K) \to K $ is an isomorphism
by \ci{borel}, Lemma V.10.14. 
The second isomorphism follows from the first one and the 
identification
$ i^{*}  \pi^{*}\simeq Id_{Z}^{*}:$ $\pi^{*} \pi_{*} K
\simeq \pi^{*} (i^{*} \pi^{*}) \pi_{*} K 
\simeq \pi^{*} i^{*} K =   \pi^{*} K_{|Z}.$ 

\n
(b) We use the natural identifications
${\cal D}^{2} \simeq Id,$ $i^{!} \simeq {\cal D}_{Z} i^{*}
{\cal D}_{W}$ and the fact that $\pi^{!}\simeq \pi^{*} [2d]$, for
$\pi$ is a locally trivial $\comp^{d}-$bundle. Denote
the dualizing complexes of $W$ and $Z$ by $\omega_{W}$ and
$\omega_{Z}.$ One has $\omega_{W} \simeq \pi^{!} \omega_{Z}.$
We have
$i^{!}K \simeq {\cal D}_{Z} i^{*} Rhom (K, \omega_{W})
\simeq {\cal D}_{Z} Rhom (i^{*}K, i^{*} \pi^{!} \omega_{Z})
\simeq
({\cal D}_{Z}  Rhom (i^{*}K, i^{*} \pi^{*} \omega_{Z}) ) [-2d]
\simeq
({\cal D}_{Z}  Rhom (i^{*}K,  \omega_{Z}) ) [-2d]
\simeq {\cal D}^{2} (i^{*}K) [-2d].$
\blacksquare

\subsection{Perverse sheaves}
\label{perverse}
Let $Y$ be an algebraic variety.
We consider the $t$-structure on $D(Y)$ associated with the middle 
perversity
see \ci{bbd}, \ci{k-s}, $\S$10.
The associated 
 heart is denoted by  $Perv (Y)$ and it  is a full 
 {\em abelian} sub-category of $D(Y).$  Its objects are 
 called
 {\em perverse sheaves},  despite the fact
 that they  are complexes.
 In short, we have the following structure.
 
 \begin{itemize}
     
     \item
    Two full sub-categories $D^{\leq 0}(Y)$ and $D^{\geq 0} (Y)$
     of $D(Y):$ 
     $$
     Ob ( D^{\leq 0}(Y) ) = \left\{ K \in D(Y) \, |  \;
     \dim{\,\mbox{supp} \, \csix{j}{K} } \leq -j, \; \forall j \right\},
     $$
     $$
     Ob ( D^{\geq 0}(Y) ) = \left\{ K \in D(Y) \, |  \;
     \dim{\, \mbox{supp} \, \csix{j}{{\cal D} {(K)}} } \leq -j, \, 
     \forall j \right\}.
     $$
     Set  $D^{\leq m} (Y)$ $:= D^{\leq 0}(Y) [-m],$
     $D^{\geq m} (Y):= D^{\geq 0}(Y) [-m]$.
\end{itemize}

\begin{rmk}
\label{refpss}
{\rm 
     These conditions can be 
     re-formulated using a stratification ${\frak Y}$
     as follows.
     Let  $K$ be ${\frak Y}-$cc and
     $\alpha_{l}: S_{l} \to Y$ be the corresponding embedding.
     We have:
     
\n
$K \in
Ob ( D^{\leq 0}(Y) ) $ if and only if
$\csix{j}{\alpha_{l}^{*}K} =0,$ $ \forall \; l$  and 
$j \; s.t. \; j > - l$.
This is known as the {\em condition of support}.

\n
$F \in   Ob ( D^{\geq 0} (Y) )$ 
if and only if
$\csix{j}{ \alpha_{l}^{!}F }   =0,$  $\forall \; l$ and  
     $j \; s.t. \; j < - l.$
This is known as the {\em condition of co-support.}
}
\end{rmk}     

\begin{itemize}  
 \item
 If $P \in Perv(Y),$ then ${\cal H}^i(P)=0$ for $i \notin [-\dim{Y}, 
0]$.
More precisely, if $P$ is ${\frak{Y}}-$cc and $ 0\leq s \leq d,$ then
${\cal H}^i(P_{|U_s})=0$ for $i \notin [-\dim{Y}, -s]$. 

\n If $P \in Perv(Y)$ is  ${\frak{Y}}-$cc  and is 
supported on a closed $s-$dimensional stratum $S_{s}$, then
$P \simeq {\cal H}^{-s}(P)[s]$.

     \item
     If $F \in Ob ( D^{\leq m}(Y) ) $ and $G \in Ob
     (   D^{\geq m+t}(Y) )$ for  $t >0$, then
     \begin{equation}
         \label{otts}
     Hom_{D(Y)} (F, G) =0. 
     \end{equation}

\item
There are 
the {\em perverse truncations}  functors,
defined up to unique isomorphism,
     $\ptd{m}: D(Y) \to D^{\leq m} (Y)$,
and     $\ptu{m}: D(Y) \to D^{\geq m} (Y)$, 
adjoint to the inclusion functors, that is
$$
Hom_{D(Y)}(F, G)\, =\,
Hom_{ D^{\leq m}(Y) }(F,\ptd{m}G )  \quad \hbox{ if } F \in 
Ob ( D^{\leq m}(Y) ) 
$$
and 
$$
Hom_{D(Y)}(G,F)\, =\, Hom_{ D^{\geq m}(Y) }(\ptu{m}G, F) \quad \hbox{ if } F \in 
Ob ( D^{\geq m}(Y) ).  
$$
There are adjunction maps   $  \ptd{m}K \to  K$ and
$K \to  \ptu{m}K.$
By the boundedness hypothesis on $D(Y),$ 
$ \ptd{m}K \simeq 0 $ if $m\ll 0$ and
$ \ptu{m}K \simeq 0 $ if $m\gg 0.$

\item
If 
$K$ is ${\frak Y}-$cc, then so are
  $\ptd{m}K$ and $\ptu{m}K$.

\item 
 There are canonical isomorphisms of functors
     $$
     \ptd{m} \circ [l] \,  \simeq \,  [l] \circ \ptd{m+l}  , \qquad
     \ptu{m} \circ [l] \,  \simeq \, [l] \circ \ptu{m+l} .
     $$
     
\item
     For every $K$ and $m$  there is a functorial
     triangle
     $$
      \ptd{m} K \lorw K \lorw \ptu{m+1} K \stackrel{[1]}\lorw.
   $$

\item
 The heart $Perv(Y):= D^{\leq 0} (Y) \cap D^{\geq 0} 
(Y)$  of the $t$-structure
   is an abelian category  
   which  objects are called {\em perverse sheaves.} 
 An object $K$ of $D(Y)$ is perverse if and only if
   the two natural maps coming from adjunction
$\ptd{0} K \to K$ and $  K \to  \ptu{0} K $
are isomorphisms.

\item
   The functor 
   $$
   \phix{0}{-}: D(Y) \lorw Perv(Y), \qquad \phix{0}{K}: =
   \ptd{0} \ptu{0}K \simeq \ptu{0} \ptd{0} K,
   $$
   is cohomological.
   Define
   $$
   \phix{i}{K}:= \phix{0}{K[i]}.
   $$
   These functors are called {\em the perverse cohomology functors}.
   
   \n
   Any  triangle $K' \to K \to K'' \stackrel{[1]}\to$
   in $D(Y)$
   gives a long exact sequence  in $Perv(Y):$
   $$
    \ldots \lorw \phix{i}{K'} \lorw  \phix{i}{K} \lorw
    \phix{i}{K''} \lorw \phix{i+1}{K'} \lorw \ldots.
   $$
 If $K$ is 
${\frak Y}-$cc, then so are  $\phix{i}{K},$ $\forall i  \in \zed.$

\n
By boundedness and the Five Lemma,
a map $\phi :K \to K'$  is an isomorphism
iff $\phix{i}{\phi}: \phix{i}{K}\to \phix{i}{K'}$ is 
an isomorphism for every $i.$

\item
     Poincar\'e-Verdier Duality exchanges  
    $D^{\leq 0} (Y)$ with
    $D^{\geq 0} (Y)$ and  fixes $Perv(Y).$ There are canonical
     isomorphisms of functors
    $$\ptd{0} \circ {\cal D} \,  \simeq \, 
    {\cal D}\circ \ptu{0}, \qquad \ptu{0} \circ {\cal 
D} \,  
\simeq \, {\cal D} \circ \ptd{0}
\qquad
    {\cal D}  \circ ^{p}\!{\cal H}^{j}   \, 
    \simeq \,  ^{p}\!{\cal H}^{-j}  \circ  {\cal D} .
    $$ 
\end{itemize}

%\begin{lm}
%\label{qits}
%Let ${\cal A}$ be a triangulated category with 
%$t-$structure, $H^{i}$ be the associated cohomology theory
%and $G \in Ob ({\cal A})$ be such that $H^{i}(G)=0,$ 
%for every $|i|\gg 0.$

%\n
%(a)
%Let  $\phi: \oplus_{i} H^{i} (G) [-i] \lorw G$
%be such that $H^{i}(\phi)$ is an isomorphism
%for every $i$. Then $\phi$ is an isomorphism in ${\cal A}.$

%\n
%(b) Let $\phi: \oplus_{i}{P^{i}[-i]} \lorw G$
%be an isomorphism, where $P^{i}$ is in the heart of the $t-$structure
%for every $i.$ Then $H^{i}(\phi): P^{i} \lorw H^{i}(G)$ is an 
%isomorphism that composed with the isomorphism
% $diag ( H^{i}(\phi^{-1}))$
%gives an isomorphism $\phi': \oplus_{i} H^{i}(G) [-i] \lorw G$
%inducing the identity on $H^{i}(G).$
%\end{lm}
%{\em Proof.} Elementary and left to the reader.
%\blacksquare

\subsection{t-exactness}
\label{texact}
A functor $T: D_{1} \to D_{2}$ of triangulated categories with 
$t$-structures is  said to be {\em left} ({\em right,} resp.) {\em 
$t-$exact}
if $T(D_{1}^{\geq 0}) \subseteq D_{2}^{\geq 0}$ 
($T(D_{1}^{\leq 0}) \subseteq D_{2}^{\leq 0},$ resp.) and it is said 
to be {\em $t-$exact} if it is left and right $t-$exact. 
If $T$ is $t-$exact,
then it preserves the hearts of the two categories. 
In particular, if $T$ is $t-$exact, then
there is a natural isomorphism $\phix{l}{T(K)} \simeq
T ( \phix{l}{K}   ),$ $K \in Ob(D_1).$
See \ci{bbd}, especially $\S$4, and \ci{k-s},
$\S$10.

\medskip
Let $f: X \to Y$ be an algebraic map of algebraic varieties. We 
consider the triangulated categories $D(X)$ and $D(Y)$ with their
middle-perversity $t$-structure. 
Verdier Duality is an auto-equivalence of categories. It 
exchanges $f^{!}$ with $f^{*}$, $f_{!}$ with $f_{*}$ and $D^{\leq 0}(-)$ 
 with $D^{\geq 0}(-).$ Consequently,  statements about
the left (right, resp.) $t-$exactness of
the four functors $(f_{!},f^{!}, f^{*}, f_{*})$ are equivalent
to the analogous statements of right (left, resp.) $t-$exactness
of the four functors
$(f_{*},f^{*}, f^{!}, f_{!}).$ Similarly,
 $f^{*} [j]$ is left $t-$exact if and only if
$f^{!}[-j]$ is right $t-$exact, etc.

\smallskip
If $f$ is affine, then $f_{*}$ is right $t-$exact and $f_{!}$ is left
$t-$exact (cf.  \ci{bbd}, 4.1.1). 
This is a convenient re-formulation
of the theorem on the cohomological dimension
of affine spaces and it implies the Weak Lefschetz Theorem
(cf. \ref{hyperplane}).

If $f$ is quasi-finite and affine, then $f_{*}$ and $f_{!}$ are
$t-$exact.
If $f$ is smooth of relative dimension $d$, then
$f^{*}[d] \simeq f^{!}[-d]$  and they are $t-$exact. If, in addition,
$f$ is surjective and with connected fibers, then
the induced functor $f^{*}[d]: Perv(Y) \to Perv(X)$
is fully faithful.

\subsection{Intersection cohomology complexes,
semisimple objects and intermediate extensions}
\label{iccso}
Recall that  $P \in Ob (Perv(Y))$ is said to be {\em simple}
   if it has no non-trivial sub-objects and hence no non-trivial
   quotients. $P$ is said to be {\em semisimple}
   if it is isomorphic to a direct sum of simple objects.
   The category  $Perv(Y)$ is  abelian,
{\em artinian,} i.e. every 
   $P \in Perv (Y)$ admits a finite filtration by sub-objects
whose successive quotients are simple, 
and {\em noetherian}, i.e. any increasing 
filtration of $P$
   by  sub-objects stabilizes.

   Let $\beta : U \to Y$ be  a Zariski-dense 
open subset
   of $Y.$ 
   Given $P \in Ob (Perv (U) )$,
   there is  an object $\beta_{!*}P \in Perv(Y)$
   with the property that it extends $P$ and it has no non-trivial 
sub-object
   and  quotient supported on a closed subvariety  of $Y \setminus U.$
  It is  unique up to 
isomorphism and is
   called {\em the intermediate extension of $P.$}
   See \ci{bbd}, 1.4.25, 2.1.9, 2.1.11.
   
\n 
   Given any stratification
   ${\frak Y}$ of $Y$ for which $Y\setminus U$ is a union of 
connected components of strata,
   $\beta_{!*}P$ is, up to isomorphism, the 
    unique extension
   $\widetilde{P}$ of $P$  in $D(Y)$ such that, 
   given any connected component, $S \stackrel{i}\to  Y \setminus 
U$,
 of a stratum contained in
$Y\setminus U,$ we have 
   $\csix{j}{i^{*}\widetilde{P}}=0$, $\forall j  \geq -\dim{S}$
   and $\csix{j}{i^{!}\widetilde{P}}=0$, $\forall j  \leq -\dim{S}.$

\begin{rmk}
\label{iee}
{\rm 
   The intermediate extension $\beta_{!*}P$ can be described 
   explicitly in terms of  stratifications 
and   successive push-forwards and truncations:
let  ${\frak Y}$ be a stratification of $Y$
inducing stratifications on $U$ and $Y\setminus U$ with respect to
which $P$ is ${\frak Y}_U-$cc.
The construction is by
induction on the strata: if $U= U_{l+1},$ then 
${\beta_{l}}_{!*}P \simeq  \td{-l-1}{\beta_l}_* P.$
}
\end{rmk}
Let ${\cal L}$ be a local system on an open set $U$ contained in 
   the regular part  of $Y.$
   The {\em intersection cohomology complex associated with
   $\cal L$ }  is   $IC_{Y}  ({\cal L} ):= 
   \beta_{!*} ({\cal L}[\dim{Y} ]) \in Ob (Perv(Y)).$
   
   \n
   The case of ${\cal L} = \rat_{U}$ is of particular interest
   and gives rise to the intersection cohomology
   complex  $IC_{Y}$ of $Y.$ The
    groups ${\Bbb H}^{j}(Y, IC_{Y})
   \simeq I\!H^{\dim{Y}+j }(Y, \rat)$ are the  rational 
   {\em intersection cohomology
   groups} of $Y$ (cf. \ci{go-ma2}). 
   
   \n
If $Y$ is smooth, or at least a rational homology manifold, then
$IC_{Y} \simeq \rat_{Y} [\dim{Y} ].$
   
   \n
   The  complex $IC_{Y}( {\cal L} )$ is characterized, up to 
isomorphism,
   by the following conditions:
   
   \smallskip
    \n
 -
     $\csix{j}{IC_{Y} ( {\cal L} )} =0;$ for all $j <-\dim{Y};$

 \n
 -
    $\csix{-\dim{Y}}{IC_{Y}  ( {\cal L}   )_{|U}} \simeq {\cal L} $;
     
 \n
 -
     $ \dim{ \, \mbox{supp} \, \csix{j}{ IC_{Y} ({\cal L} ) }  }< -j  $, 
    if $j > - \dim {Y};$
    
  \n
  -
        $ \dim{\, \mbox{supp} \, ( \csix{j}{ {\cal D} ( IC_{Y} ({\cal L} ) } )) 
}  < -j  $, 
    if $j > - \dim {Y}.$

    \n
The last two conditions   can be 
     re-formulated using  stratifications
     as follows.
     Let ${\frak Y}$ 
     be a stratification
     with respect to which $IC_{Y} ({\cal L}   )$ is ${\frak Y}-$cc,
     $\alpha_{l}: S_{l} \to Y$ be the embedding.
     We have: 
     
\n
$\csix{j}{\alpha_{l}^{*}  IC_{Y} ({\cal L} )  } =0,$ $ \forall \; l$  
and 
$j > \dim{Y} \; s.t. \; j \geq - l$;

\n
$\csix{j}{ \alpha_{l}^{!} IC_{Y} ({\cal L}   )    } =0,$  $\forall \; 
l$ and  
     $j > -\dim{Y} \; s.t. \; j \leq  - l.$

\n   
One has ${\cal D}  ( IC_{Y} ( {\cal L} )  )\simeq 
IC_{Y} (  {\cal L}^{\vee} ).$
Given a closed subvariety $i: Y' \lorw Y$ and a complex
of type $IC_{Y'} ({\cal L'}) \in Perv (Y')$, we denote
$i_{*} IC_{Y'} ({\cal L'})$ simply by $IC_{Y'} (  {\cal L}' )$.
It is an object of $Perv (Y)$ satisfying the conditions above, with
$Y'$ replacing $Y$.
   
\smallskip
   An object $P \in Perv(Y)$ is simple if and only if
   $P \simeq  IC_{Y'} ({ \cal L'} )$, for some
   closed subvariety $Y' \subseteq Y$ and some simple 
 local system
   ${\cal L'}$ defined on an open  subvariety of the regular part
   of $Y'$.
   A semisimple $P$ is a finite direct sum
   of such objects.

   \begin{rmk}
       \label{nomaps}
       {\rm 
 If $Y'$ and $Y''$ are distinct,
then the properties of intermediate extensions imply that  
 $Hom_{D(Y)}( IC_{Y'}(L'),  IC_{Y''}(L''))=0.$
 }
       \end{rmk}

\section{Preparatory results}
\label{prepres}   
In this section we collect  a series of results needed
in the sequel
for which we could not find an adequate reference.

\subsection{A splitting criterion in $Perv(Y)$ }
\label{vsl}

One of the key results of this paper is the
geometric proof of the Semisimplicity Theorem
for the perverse sheaf $\phix{0}{\fxn{f}{X}{n}}.$
Every perverse sheaf can be written as a finite extension
of intersection cohomology complexes. We prove, using induction
on a stratification,  that all
the extensions are trivial. 
The set-up  is as follows.
Let $Y$ be an algebraic variety and 
 $s \in \nat$  be such that there is a stratification
 $\frak Y$ with 
 $Y= U \amalg S$, $U = \amalg_{l>s}
Y_{l}$  and $S= S_{s}.$ 
Denote by  $
S \stackrel{\alpha}\lorw Y \stackrel{\beta}\longleftarrow U
$
the corresponding closed and open embeddings.
Let $P \in Ob ( Perv(Y) )$  be ${\frak Y}-$cc.

\smallskip
Consider the truncation  triangle 
$$
\tau_{\leq -s-1}P \lorw \tau_{\leq -s} P \stackrel{\tau}\lorw \tau_{\geq -s}P 
\stackrel{[1]}\lorw.
$$
The  conditions  of support in \ref{perverse}
imply that 
$ \td{-s-1} \beta^{*}P \simeq \beta^{*} P, $
$P \simeq \tau_{\leq -s} P$ and that  
$\tu{-s} P \, \simeq \,
\csix{-s}{P} [s]$ is a (shifted) local system supported on $S.$
Using   the long exact sequences associated with
the  truncation triangle and with the triangles obtained
from it by applying $\alpha^{*}$ and $\alpha^{!},$ one checks that 
the complex $\tau_{\leq -s-1}P$ is perverse. We get a short
exact sequence in $Perv(Y):$
\begin{equation}
    \label{sesip}
    0 \lorw \tau_{\leq -s-1}P \lorw P \stackrel{\tau}\lorw \csix{-s}{P} [s] 
    \lorw 0.
    \end{equation}
The deviation of 
$\td{-s-1}P$ from being the intermediate extension
$\tau_{\leq -s-1}\beta_{*} \beta^{*} P
\simeq \beta_{!*} \beta^{*}P$
is measured  by the   map 
$
 \tau_{\leq -s-1}P \lorw \tau_{\leq -s-1}\beta_{*} \beta^{*} P
$ 
which arises from truncating the adjunction map
$P \to \beta_{*} \beta^{*} P.$

Since $P \simeq \td{-s}P,$ the adjunction map 
 admits a canonical lifting
$$
l \, : \, P \lorw \td{-s} \beta_{*} \beta^{*} P.
$$

We are looking for a condition implying that $P$
is a direct sum of intersection cohomology complexes, or equivalently,
that $P  \simeq \beta_{!*} \beta^{*}P \oplus \csix{-s}{P}[s]$
holds for every stratum $S_{s}.$ In this context the following 
assumption is natural (cf. Remark \ref{condmet}).
\begin{ass}
\label{natass}
$\dim_{\rat}{  (\csix{-s}{ \alpha_{!}\alpha^{!}P} )_{y} } \, =\,  
\dim_{\rat}{ (\csix{-s}{ \alpha_{*}\alpha^{*}P   })_{y}} , \; \, y \in S.$
\end{ass}

\begin{rmk}
\label{condmet}
{\rm
Assumption \ref{natass}
is automatically satisfied if, for example,
$P$ is a direct sum of intersection cohomology complexes
(this is what we are aiming to  prove for 
$\phix{0}{\fxn{f}{X}{n}}$)
or if $P$ is  self-dual (in our case, this  is automatic 
by Verdier duality). 
In the former case, by \ref{iccso}, $P$ must be isomorphic
to $\beta_{!*} \beta^{*}P \oplus \csix{-s}{P} [s]$
and \ref{natass} follows from
 the natural isomorphisms $\alpha^{*} \alpha_{*} \simeq
Id \simeq \alpha^{!} \alpha_{*}.$
In the latter case, we  apply the duality statement
 of Lemma \ref{adapt}.  
}
\end{rmk}
The triangle
$
 \alpha_{!} \alpha^{!}P \to  P \to  \beta_{*} \beta^{*}P 
\to $
gives
\begin{equation}
    \label{abc}
\csix{-s-1}{P} \stackrel{a}\lorw \csix{-s-1}{\beta_{*}\beta^{*}P} \lorw 
\csix{-s}{\alpha_{!}\alpha^{!}P} \stackrel{\iota}\lorw
\csix{-s}{P} \stackrel{b}\lorw \csix{-s}{ \beta_{*} \beta^{*}P}.
\end{equation}

\begin{lm}
\label{splitp}
({\bf Splitting Criterion})
Assume \ref{natass}. The  
following  are equivalent:

\smallskip
\n
1) $P \simeq \beta_{!*}\beta^{*}P \oplus \csix{-s}{P}[s].$

\n
2) $
\iota : \, \csix{-s}{\alpha_{!} \alpha^{!} P } \, \lorw  \, \csix{-s}{P}
$
is an isomorphism.

\n
3) The map 
$l: P \to
\td{-s} \beta_{*} \beta^{*}P
$
has a  lifting  $\widetilde{l}: P \, \to  \,  
 \beta_{!*} \beta^{*}P.$

 \smallskip
\n
If 3) holds, then the lifting $\widetilde{l}$ is unique and
 gives the  natural isomorphism (cf. (\ref{sesip}))
$$
(\widetilde{l},\tau)\,:\, P \,\simeq \,  
\beta_{!*}  P \, 
\oplus\,  \csix{-s}{P}[s].
$$
\end{lm}
{\em Proof.} By the characterization of intermediate extensions
(cf. \ref{iccso}) and the fact that $\alpha^{!}\alpha_{*} \simeq 
\mbox{Id},$ 1) implies 2).

\n
Apply the functor $Hom_{D(Y)} (P, -)$ to the triangle
$$
 \td{-s-1} \beta_{*} \beta^{*}P \lorw \td{-s}  \beta_{*} \beta^{*}P  
\lorw
  \csix{-s}{ \beta_{*} \beta^{*}P} [s] \stackrel{[1]}\lorw.
 $$
By (\ref{otts}),  
 $Hom^{-1}_{D(Y)} (P, \csix{-s}{ \beta_{*} \beta^{*}P}[s]) =\{0 \}.$ 
The associated long exact sequence shows that
if $\widetilde{l}$  exists, then it is unique.
Furthermore,
$\widetilde{l} $ exists if and only if the image
of $l$ in
$Hom_{D(Y)} (P, \csix{-s}{  \beta_{*} \beta^{*}P } [s] )$
is zero. This is equivalent to 
$b=0$ (cf. (\ref{abc})) and hence to 
$\iota$ being surjective hence an isomorphism 
(cf. \ref{natass}). This shows that 2) and 3) are equivalent.

\n
Assume that  $\widetilde{l}$  exists. We have an 
 exact sequence in the abelian category  
$Perv(Y): $
\begin{equation}
    \label{esip}
0 \lorw K \lorw P \stackrel{\tilde{l}}\lorw   \beta_{!*} 
\beta^{*}P 
\lorw C \lorw 0.
\end{equation}
Since $\widetilde{l}$ is an  isomorphism over $U,$
the complexes $C$ and $K$ are supported on $S.$
Since intermediate extensions do not admit
quotients supported on $S,$ the complex $C=0.$
By the conditions of support
\ref{perverse},  $K \simeq 
\csix{-s}{K}[s]\simeq \csix{-s}{P}[s]$ is a shifted local system
on $S.$
The sequence (\ref{esip})
 reduces to the triangle 
 \begin{equation}
     \label{bdl}
  K \lorw  P  \stackrel{\widetilde{l}}\lorw \beta_{!*} 
  \beta^{*}P \stackrel{[1]}\lorw     
  \end{equation}
 whose long exact sequence contains
$$
0 \lorw 
\csix{-s-1}{P} \stackrel{\widetilde{a}} \lorw \csix{-s-1}{\beta_{*}\beta^{*}P}
\stackrel{c}\lorw \csix{-s}{K} \lorw \csix{-s}{P} \lorw 0.
$$
Since $\iota$ is injective and $\tilde{l}$ is a lifting 
of $l$,  we have that $a=\widetilde{a}$ is surjective, so that $c=0$.

\n
Since $Hom_{ D(Y) } ( \td{-s-1} \beta_{*}\beta^{*}P,
K[1]) 
\simeq Hom_{ Sh(Y) } ( \csix{-s-1}{ \beta_{*} \beta^{*} P },
\csix{-s}{P} )$, we see that $c=0$ implies that
the triangle (\ref{bdl}) splits, i.e. that there is some isomorphism
$P \simeq  
 \beta_{!*} \beta^{*} P  \oplus  \csix{-s}{P} [s].$
It follows that  $\widetilde{l} \oplus \tau$ is an isomorphism
on cohomology sheaves, hence an isomorphism.
 The fact that 3) implies  1) is now trivial.
\blacksquare

\subsection{The perverse filtration}
\label{tpf}
 The perverse truncation functors
define  increasing filtrations in hypercohomology.
\begin{defi}
    \label{defofpf}
    {\rm 
    Let $K \in Ob (D(Y))$ and $j \in \zed.$  The
{\em perverse filtration} on  $\ihixc{j}{Y}{K}$ is defined by setting:
\begin{equation}
    \label{defoffilth}
    \fihixcj{j}{i}{Y}{K} \,  :=\, \im{ \,\left\{ 
\ihixc{j}{Y}{\ptd{i} K} \to \ihixc{j}{Y}{K}  \right\}}.
\end{equation}
The graded pieces  are
\begin{equation}
    \label{defofpcg}
\gihixcj{j}{i}{Y}{K}\, : = \,
\fihixcj{j}{i}{Y}{K}/ \fihixcj{j}{i-1}{Y}{K}.
\end{equation}
}
\end{defi}

We have canonical maps, $a$ injective, $b$ surjective:
\begin{equation}
    \label{pccm}
\ihixc{j}{Y}{ \phix{i}{K}} \, \stackrel{a}\longleftarrow\, 
\mbox{Coker} \left\{ \ihixc{j}{Y}{\ptd{i-1} K}  
\to  \ihixc{j}{Y}{\ptd{i} K} \right\} \, \stackrel{b}\lorw\,  
\gihixcj{j}{i}{Y}{K}. 
\end{equation}
Given a map $\phi: K \to K'$ in $D(Y),$ the map
${\Bbb H}^{j}( \phi)$ is {\em filtered}, but not necessarily
 strict.

\smallskip
Let $f: X \to Y$ be a  map of algebraic varieties,
$n:=\dim{X}.$
There is a  canonical
isomorphism
$
H^{n+j}(X)\simeq \ihixc{j}{Y}{\fxn{f}{X}{n}}.
$
\begin{defi}
\label{pcg}
{\rm
({\bf Perverse filtration and cohomology groups})
Let $i,\, j \, \in \zed$ and set:
$$
H^{n+j}_{\leq i}(X) \, :=\,
\fihixcj{j}{i}{Y}{\fxn{f}{X}{n}}, 
\qquad
\qquad
H^{n+j}_{i}(X): = H^{n+j}_{\leq i}(X)/H^{n+j}_{\leq i-1}(X).
$$
We call the groups
$H^{l}_{b}(X)$ the {\em perverse cohomology groups of 
$X$} (relative to $f$).
}
\end{defi}

Note that,
in view of the natural equivalence $\ptd{i} \circ [l] \simeq 
[l]\circ \ptd{i+l},$ we have 
$$
H^{n+j}_{\leq i}(X) \, \neq \, \fihixcj{n+j}{i}{Y}{f_{*}\rat_{X}}
\, = \, 
 H^{n+j}_{\leq i+n}(X).
$$
The dimensional shift in Definition
\ref{pcg} is convenient for our purposes. 

\begin{rmk}
    \label{filtrfac}
    {\rm 
Let $r:Z^{m} \to  X^{n},$ $f:X \to Y$ be  maps of algebraic
    varieties of the indicated dimensions, $g:= f \circ r.$
    The natural map $\fxn{f}{X}{n}\to \fxn{g}{Z}{m}[n-m]$ and
    the rule
     $\ptd{i} \circ [l] \simeq 
[l]\circ \ptd{i+l} $ imply that the natural restriction map satisfies
(cf. \ref{joined}):  
   $$
    r^{*}\, : \, H^{l}_{\leq b}(X) \lorw H^{l}_{\leq b + n-m}(Z).
$$
}
\end{rmk}

\subsection{p-splitness and some consequences}
\label{psic}
In our inductive approach to the results of this paper,
the Decomposition Theorem \ref{tm1}.b is established rather early
(cf. \ref{rhl}). 
In this section, we collect some simple consequences of the
Decomposition Theorem which are used in the remainder
 of the proof by 
induction.

\begin{defi}
    \label{p-split}
    {\rm ({\bf p-split}) A complex $K \in D(Y)$ is said to be {\em p-split}
    if there is an isomorphism 
    $$
    \phi \, : \,  K\, \simeq\,  \dsdix{i}{K}.
    $$
    }
    \end{defi}

\begin{rmk}
\label{always}
{\rm 
We can and shall always consider isomorphisms satisfying
the  additional condition $\phix{i}{\f} = Id_{ \phix{i}{K}  }.$ 
Given $\phi,$ the map
$$
\f \, : \, \left(  \sum_{i}{ \phix{i}{\phi^{-1}}[-i]}    \right) 
\circ \phi  \, : \, K \lorw
\bigoplus_{i}{ \phix{i}{K} [-i] }
$$
is one such isomorphism.
}
\end{rmk}
\begin{rmk}
\label{qeiso}
{\rm 
If $\phi \, : \,  K\, \simeq\,  \oplus_{i}P_{i}[-i]$ is an 
isomorphism
with $P_{i} \in Perv(Y)$ for every $ i \in \zed,$ then 
$\phix{i}{\phi}: \phix{i}{K} \simeq P_{i}$ for every $i \in \zed$
and $K$ is p-split.
}
\end{rmk}

\medskip
In the remainder of this section, $f:X \to Y$ is a proper map of 
algebraic varieties, $X$ is nonsingular,
$n :=\dim{X},$ $A$ is an ample line bundle on $Y$ and 
$L:= f^{*}A.$

\begin{rmk}
    \label{comincia}
   {\rm  If 
   $ \fxn{f}{X}{n}$ is  p-split, then 
   for  $b, \, j \, \in \zed$  we have isomorphisms:
   $$
H^{n+j}_{\leq b}(X) \, \stackrel{\f}{\simeq} \, 
\bigoplus_{i \leq b}{\Bbb H}^{j-i}(Y, 
   \phix{i}{\fxn{f}{X}{n}} ), 
$$
$$
H^{n+j}(X) \stackrel{\f}{\simeq} \bigoplus_b  H^{n+j}_{b}(X).
$$
On the other hand, the maps $a$ and $b$ in (\ref{pccm})
are isomorphism and we get a {\em canonical} identification:
\begin{equation}
    \label{canid}
a \circ b^{-1} : \,   H^{n+j}_{b}(X)  \, =  \, {\Bbb H}^{j-b}(Y, 
   \phix{b}{\fxn{f}{X}{n}} ). 
\end{equation}
The isomorphism induced by  $\f $ coincides with the one above
(cf. \ref{always}).
%if one assumes $\phix{i}{\f}=Id_{\phix{i}{F}}.$
}
\end{rmk}

\smallskip
 Let $R: D(Y) \to D(Y)$ be a functor of triangulated categories, 
$\nu: R \to Id$ ($Id \to R$, resp.) 
be a natural transformation of functors of triangulated categories
compatible with coproducts 
 and $  K \in Ob(D(Y))$ be p-split. 
The  splitting of $K$
induces
a filtration on ${\Bbb H}^{*}(R(K)).$ This filtration is independent
of the splitting $\f$
and need not to coincide with the perverse 
filtration.
The maps ${\Bbb H}^{*}(Y, \nu)$ are
 strict.

\begin{defi}
    \label{defigf}
    {\rm Given $R$ and $K$ p-split as above,
   the filtration
    on ${\Bbb H}^{*}(Y, R(K))$ described above  is called
     the {\em induced   filtration}.
     }
     \end{defi}

\begin{lm}
\label{collectinfo} 
 Assume that
$\fxn{f}{X}{n} $ is  p-split.
Let $U \subseteq Y$ be an euclidean open neighborhood
of $y \in Y$ and  $U'=f^{-1}(U).$
The natural maps 
$$
H^{BM}_{n-l}(f^{-1}(y)) \to
H^{BM}_{n-l}(U') \to
H^{BM}_{n-l}(X) \simeq
H^{n+l}(X) \to H^{n+l}(U') \to
H^{n+l}(f^{-1}(y))
$$
are  strict with respect to the induced filtrations
for every $l \in \zed.$

\n
For every $b \in \zed$:
 $$
 H^{BM}_{n-b, \leq b} (f^{-1} (y) ) = H^{BM}_{n-b} (f^{-1} (y) ),
 \qquad 
 H^{n+b}_{\leq b-1} ( f^{-1} (y) ) = \{ 0 \}.
 $$
In particular, the map induced on the graded space
$$
H^{BM}_{n-b,a} (f^{-1} (y) ) \lorw
H^{n+b}_{a} (f^{-1} (y) )
$$
is the zero map for every $a \neq b.$
\end{lm}
{\em Proof.}
Let  $\alpha: y \to Y$ and
$\beta: U \to Y$ be the embeddings.
Note that $f_{*}\omega_{X}[-n] \simeq \fxn{f}{X}{n} $ is p-split
and that  the induced  filtration
on $H^{*}(f^{-1} (U) )$ coincides with the perverse filtration.
The strictness assertions follow from the discussion preceding
  \ref{defigf} applied to the duality functor and to  
 the two adjunction maps associated with $\alpha$ and $\beta.$

\n
The third  statement follows immediately from the second one
which in turn follows from the
 conditions of (co-)support  
  (cf. Remark
\ref{refpss}):
$$
\ihixc{t}{Y}{ \alpha_{!} \alpha^{!} \phix{l}{\fxn{f}{X}{n}      
[-l] }}
\,=\, \{0 \},
\qquad \forall  \, l > b,
$$
$$
\ihixc{t}{Y}{\alpha_{*} \alpha^{*} \phix{l}{\fxn{f}{X}{n}  } [-l] } 
\, =\, \{0 \},
\qquad \forall  \, l < b.
$$
\blacksquare

\begin{rmk}
 \label{isofib}
{\rm
Theorem \ref{nhrbr} states  that   the refined
intersection product 
induces  
isomorphisms  $H^{BM}_{n-b,b}(f^{-1}(y)) \simeq
H^{n+b}_{b}(f^{-1} (y) )$  (cf. \ref{alk}).
 }
\end{rmk}

\begin{lm}
\label{dtpb}
Let
$$
\xymatrix{
{X'}^{m} \ar[r]^{v} \ar[d]^{f'} & X^{n} \ar[d]^f \\
Y'  \ar[r]^{u}  & Y}
$$
be a Cartesian diagram of maps of algebraic varieties
of the indicated dimensions, $f$ proper.
Let $\frak X$ and $\frak Y$ be a stratification for $f.$
Assume   that $\fxn{f}{X}{n}$ is p-split
% there is an isomorphism
%$\phi: \dsdix{i}{ \fxn{f}{X}{n} } \simeq \fxn{f}{X}{n}$
and that
either $u$ is smooth, or it is a
normally nonsingular inclusion.

\n
Then $\fxn{f'}{X'}{m}$ is p-split, 
$
\phix{i}{ \fxn{f'}{X'}{m} } 
\simeq  u^{*} \phix{i}{\fxn{f}{X}{n}} [m-n],$
for every  $i \in \zed$ 
and the  natural map $v^{*}: H^{k}(X) \lorw H^{k}(X')$
is  compatible with the direct sum decomposition in 
perverse cohomology groups and it is strict.

\n
Let $y' \in Y'$ and $y=u(y').$ Then the  filtrations induced by 
$f$ and $f'$ on the (co)-homology of $f^{-1}(y)={f'}^{-1}(y')$
coincide.
\end{lm}
{\em Proof.} 
Let $\f: \dsdix{i}{ \fxn{f}{X}{n} } \simeq \fxn{f}{X}{n}$ be 
a splitting.
 By   base change, 
\begin{equation}
 \label{okko}
u^{*} \f \, :  \, u^{* } \left( \dsdix{i}{ \fxn{f}{X}{n} } )  [m-n] 
\right) \, 
\stackrel{\simeq }\lorw \,
u^{*} \fxn{f}{X}{n} [m-n] \simeq \fxn{f'}{X'}{m}.
\end{equation}
By \ref{texact} and
by  Remark
\ref{resshperv}, respectively,
the hypotheses imply that $u^{*}[m-n]$ is $t-$exact
so that the left-hand side is a direct sum of shifted
perverse sheaves.
The first two statements follow (cf. \ref{qeiso}).

\n
The  pull-back map $v^{*}$ on cohomology stems from
adjunction and base change,
$f_{*} \rat_{X} \to f_{*}  v_{*}v^{*}\rat_{X} \simeq
u_{*}f'_{*}  v^{*}\rat_{X} \simeq u_{*} u^{*}f_* \rat_{X},$
and
 preserves direct sum decompositions so that
 $v^{*}$ is  strict.
 
 \n
 The last statement follows from (\ref{okko}).
\blacksquare

\begin{rmk}
\label{joined}
{\rm The same conclusions hold
if $u=u'' \circ u',$ 
with $u'$ smooth and $u''$ a ${\frak Y}$-transverse embedding
into $Y.$ Note how the last statement improves on
Remark \ref{filtrfac}.
}
\end{rmk}

\subsection{The cup product with a  line bundle }
\label{cupetap}
Let $\mu$ be a line bundle on $X$ and 
denote by the same symbol its first Chern class
$\mu \in H^2(X)\simeq Hom_{D(X)}(\rat_X,\rat_X[2]).$
For every
$K \in Ob(D(X))$,
the isomorphism 
$K \simeq K \stackrel{\Bbb L}{\otimes} \rat_X$
defines
a map $\mu : K \to K[2].$

\begin{rmk}
\label{classicchern}
{\rm Suppose $s \in \Gamma(X,\mu)$ is a section  whose 
zero locus defines a normally nonsingular codimension 
one inclusion $i: \{s=0 \} \to X$; see \ref{locstr}.
By Lemma \ref{ugus}.
the map $\mu : K \to K[2]$ can be described geometrically as the composition:
$$
K \,\lorw \, i_*i^*K\,  \simeq \, i_!i^!K[2] \, \lorw K[2].
$$
If $K=\rat_X$, then  $\mu$ is the cohomology class associated with the 
normally nonsingular inclusion $i$ and we find one of the classical 
definitions of first Chern class.
}
\end{rmk}

The resulting map
${\Bbb H}^i(X, \mu):{\Bbb H}^i(X, K) \to {\Bbb H}^{i+2}(X,K)$ 
is the cup product with $\mu$.

\n
By functoriality, we get maps 
$$
f_*\mu: \, f_*K \lorw f_*K[2],
$$
$$
\ptd{i}f_*\mu: \, \ptd{i}f_*K 
\lorw (\ptd{i+2}f_*K)[2],
$$ 
$$
\phix{i}{f_*\mu}:\, \phix{i}{f_*K} \lorw \phix{i+2}{f_*K}.
$$
Applying the functor ${\Bbb H}^{*}(Y, -),$ we obtain the cup product with $\mu$
$$
{\Bbb H}(X, \mu)=
{\Bbb H}(Y, f_*\mu): \fihixcj{*}{a}{X}{K} \to \fihixcj{*+2}{a+2}{X}{K}
$$
which is filtered in an obvious sense 
and defines a cup product map
on the graded objects, still denoted 
\begin{equation}
\label{cupsh}
\mu : \, H^*_a(X,K) \, \lorw \, H^{*+2}_{a+2}(X,K).
\end{equation}

Let  $\f : K \simeq \dsdix{i}{K}$ be a p-splitting
as in \ref{always}. We have
$$
\f[2] \circ f_{*}\mu \circ \f^{-1} \, = : \, \widetilde{\mu}
 \, = \, \sum_{ij}{\widetilde{\mu}_{ij}} \, : \, 
  \dsdix{i}{K}  \lorw  \dsdix{j}{K}[2].
$$
By the choice of $\f,$ $\widetilde{\mu}_{i,i+2} = \phix{i}{f_{*} 
\mu}.$
By (\ref{otts}), $\widetilde{\mu}_{ij}= 0$ for $j > i+2.$ In general,
$\widetilde{\mu}_{ij} \neq 0$ for $ j < i+2;$ e.g. if $f =Id_{X},$ then
$\mu = \widetilde{\mu} = \widetilde{\mu}_{00}.$
It is immediate to verify 
the following important compatibility:
\begin{lm}
\label{compa}
Let $f_*K$ be $p$-split (cf. \ref{always}),  
then the isomorphism 
 $ {\Bbb H}^{*-i}(Y, \phix{i}{f_*K})\simeq $
 ${\Bbb H}^*_i(X,K)$ 
 of \ref{comincia}
is compatible with $\mu$, i.e. the 
following diagram is commutative
$$ 
\xymatrix{
{\Bbb H}^*_i(X,K) \ar[rrr]^{\mu}  \ar[d]^{\simeq} & 
&  & {\Bbb H}^{*+2}_{i+2}(X,K)  
\ar[d]^{\simeq} \\
{\Bbb H}^{*-i}(Y, \phix{i}{f_*K} )  
\ar[rrr]^{{\Bbb H}^{*-i}(Y, \phix{i}{f_* \mu})} 
& & &{\Bbb H}^{*-i}(Y, \phix{i+2}{f_*K}).
}
$$
\end{lm}

\begin{rmk}
\label{sepullback}
{\rm Let $\nu$ be a line bundle on $Y,$
denote the first Chern class in $H^{2}(Y)$ by $\nu$ 
and let $\mu= f^*\nu.$ 
Then $f_* \mu = f_*f^*\nu:f_*K \to f_*K[2]$ 
coincides with $\nu,$ as it can be seen by
considering a section $s \in \Gamma(Y,\nu)$ such that $ \{s=0 \} \to Y$ and 
$ f^{-1}(\{ s=0 \})=\{f^*s=0 \} \to X$ are normally nonsingular.
This also follows from the Change of Coefficients Formula
in \ref{stabi}.
In particular, if   
$f_*K\, \simeq\,  \dsdix{i}{f_*K}$
is p-split, then $f_*\mu=f_*f^*\nu=\nu$ is a direct 
sum map   sending 
$\phix{i}{f_*K}$ to $\phix{i}{f_*K}[2].$
In this case,
$\phix{i}{{(f_*\mu)}^{r}}:\phix{i}{f_*K} \to \phix{i+2r}{f_*K}$
is the zero map, for every $r >0$
and the cup product map with a pull-back is filtered strict
$$
f^* \nu: \, H^{*}_{i}(X) \lorw H^{*+2}_{i}(X) 
$$
with a compatibility  analogous to the one in  Lemma \ref{compa}.
}
\end{rmk}

\begin{tm}
\label{serhl}
Suppose that  
$$
\phix{-i}{f_*(\mu)^i}: \, \phix{-i}{\fxn{f}{X}{n}} \simeq  \phix{i}{
\fxn{f}{X}{n}}, \quad \forall \, i \geq 0.
$$
Then:

\n
a) $f_* \rat_X[n]$ is p-split.

\n
b)  Let $i \geq 0,$ $k \leq l$ and define 
$$
{\cal P}^{-i}_{\mu} \, :=\,  \ke \left\{
\phix{-i}{f_*(\mu)^{i+1}}:\, \phix{-i}{f_*\rat_X[n]} \to 
\phix{i+2}{f_* \rat_X[n]} \right\},
$$ 
$$
\mu^k {\cal P}^{-l}_{\mu}  \, :=\, \im \left\{ 
\phix{-l}{f_*\mu^{k}}:{\cal P}^{-l}_{\mu} \to 
\phix{-l+2k}{f_* \rat_X[n]} \right\}.
$$
(Note that the map above is a split monomorphism). 
There is a direct sum decomposition:
$$
\phix{-i}{f_*\rat_X[n]} \, \simeq \,
\bigoplus_{k \geq 0} \mu^k {\cal P}^{-i-2k}_{\mu}, \qquad 
\qquad
\phix{i}{f_*\rat_X[n]} \, \simeq \,
\bigoplus_{k \geq 0} \mu^{i+k} {\cal P}^{-i-2k}_{\mu}.
$$
c)  Let $i\geq 0,$ $j \in \zed.$   The isomorphisms 
$$
H^{n+j}_{-i}(X) \, \simeq \,  {\Bbb H}^{j}(Y, \phix{-i}{f_*\rat_X 
[n]} [i] ) \, 
\simeq \,
\bigoplus_{k\geq 0} {\Bbb H}^{j+i}(Y,\mu^k {\cal P}^{-i-2k}_{\mu} )
$$
identify: 

\n
1) $\ke{\, \mu^{i+1}} \subseteq H^{n+j}_{-i}(X)$ with the summand
${\Bbb H}^{j+i}(Y,{\cal P}^{-i}_{\mu})$ and

\n
2) the image of the injection $\mu^k: 
H^{n+j-2k}_{-i-2k}(X) \to H^{n+j}_{-i}(X)$
 restricted to 
${\Bbb H}^{j+i}(Y, {\cal P}^{-i-2k}_{\mu}),$
 with the  summand ${ \Bbb H}^{j+i}(Y,\mu^k {\cal P}^{-i-2k}_{\mu}).$ 
\end{tm}
{\em Proof.}
For   a) and b) see \ci{dess}. The rest  follows
from the constructions and Lemma \ref{compa}.
\blacksquare

\begin{rmk}
\label{nprim}
{\rm 
Given  $\f :\fxn{f}{X}{n} \simeq \dsdix{i}{\fxn{f}{X}{n}},$
the space
$\f^{-1} (\ihixc{-i}{Y}{ {\cal P}_{\eta}^{-i}[i]})$
$\subseteq H^{n-i}_{\leq -i}(X)$ is not
contained in  $\ke{\,\eta^{i+1}},$ i.e. the space of {\em
classical} primitive classes. What is true
is that $\ihixc{l}{Y}{ {\cal P}_{\eta}^{-i}[i]}$
is the kernel of $\eta^{{i+1}}: H^{n+l}_{-i}(X) 
\lorw
H^{n +l+2i+2}_{i+2}(X).$ 
}
\end{rmk}

\begin{rmk}
{\rm In the sequel of this paper, for simplicity of notation, 
we shall denote by the same symbol  all the cup 
product morphisms, without indicating which functor has been applied;
e.g. we will use $\mu$ instead of $f_*\mu,$   $\ptd{i}f_*\mu,$
$\phix{i}{f_*\mu}$ etc...
}
\end{rmk}

\subsection{Weight filtrations}
\label{fs}
For the notions introduced in this section see \ci{weil2} and \ci{stzu}.
An increasing filtration $W$ on a finite
dimensional vector space $H$ is
a collection of subspaces $W_{i}\subseteq W$ such that
$W_{i-1}\subseteq W_{i}$ for every $i \in \zed.$
The associated 
graded spaces are  $Gr^{W}_{i}H: = W_{i}/W_{i-1}$ 
The pair
$(H,W)$ is called a {\em filtered  space}.

\n
A {\em splitting} of  $(H,W)$ is an 
isomorphism $H \simeq \oplus \,Gr_i^W H.$

\n
Let $j \in \zed.$ The {\em shifted} filtration $W[j]$ is defined by
setting $W[j]_{i}:= W_{j+i}.$

\n
A {\em filtered} map
$\f:  (H,W) \to (H',W')$ is a linear map $\f: H \to H'$ such
that $\f(W_{i}) \subseteq  W'_{i}.$ It induces  linear maps
$Gr \f: Gr^{W}_{i}H \to Gr^{W'}_{i}H'$ which  we simply denote
by $\f.$

\n
A  filtered map $\f$ is called {\em strict} if
 $\f(W_{i}) = \f(H) \cap W'_{i}.$
 
 \n
Given a space $H'$ 
which is either a subspace or a quotient of $H,$
one easily defines 
an induced filtration $W(H')$ so that the maps in sight are filtered. 
In particular, given two filtrations $W$ and $W'$ on $H,$ the associated 
graded spaces $Gr^{W}_{i}H: = W_{i}/W_{i-1}$ are naturally filtered
by $W'.$

\smallskip
Given a finite dimensional vector space $H$ and a nilpotent 
endomorphism $N,$ there is a unique filtration $W$ with the 
properties that i) $N W_{i} \subseteq  W_{i-2}$ and, denoting
again by $N$ the induced map on graded spaces, ii) $N^{i}:
Gr^{W}_{i} H \simeq  Gr^{W}_{-i} H,$ for every $i\geq 0$
(cf. \ci{weil2}).

\n
Set $N^{i}= 0,$  $P^{-i}= 0 $ if $i <0,$ 
$P^{-i}= \ke{\, N^{i+1}} \subseteq Gr^{W}_{i},$ if 
$i\geq 0.$
There is the  Lefschetz decomposition
\begin{equation}
\label{ald}
Gr^{W}_{i} \, = \, \bigoplus_{l \in \zed}{ N^{-i+l} P^{i-2l} }, \quad
i \, \in \zed.
\end{equation}
This unique filtration is called the {\em weight filtration} of $N$ and is 
given by the ``convolution formula''
\begin{equation}
    \label{convform}
W_k\, =\, \sum_{i+j=k} \ke \, N^{i+1} \cap \im \, N^{-j}.
\end{equation}
 We denote the weight  filtration  of $N$ by $W^{N}$
 and
the graded spaces by $Gr^{N}_{i}.$
 
\n
Let $S$ be a  nondegenerate bilinear form on $H$ which is either
symmetric or skew-symmetric and satisfies
\begin{equation}
    \label{iaon}
S(Na,b)+S(a,Nb)=0.
\end{equation}
When (\ref{iaon}) holds, one says that
$N$ is an {\em infinitesimal automorphism} of $(H,S).$
In this case,  the weight filtration
is self-dual, i.e.
\begin{equation}
 ({W _{i}^{N}})^{\perp} \, = \, W_{-i-1}^{N}, \qquad i \in \zed  
    \end{equation}
and $S$ descends to nondegenerate
forms on  $Gr^{N}_{i}H$ for every $i \in \zed.$
More precisely,
\begin{equation}
    \label{sni}
    S^{N}_{i}( [a],[b])\, :\, =\,
S(a,N^{i}b), \quad \quad   \mbox{for} \; i \geq 0
\end{equation}
and one  requires
that $N^{i}: Gr^{N}_{i} H \simeq  Gr^{N}_{-i} H$
is an isometry for every  $i \geq 0.$ This is achieved  as follows:
if  $[a], [b] \in Gr_{-i}^{N}H$, we have
$[a]= N^{i}[a']$ and $[b]= N^{i}[b']$ for  
uniquely determined $[a'], [b'] \in  Gr_{i}^{N}H;$ set 
\begin{equation}
    \label{snin}
    S_{-i}^{N}([a],[b])\,:\,=\,
S_{i}^{N}( [a'], N^{i}[b']), \quad i >0.
\end{equation}

\n
The Lefschetz decomposition (\ref{ald})
is $S^{N}_{i}-$orthogonal for every $i \in \zed.$

\smallskip
Let $(H,W)$ be a filtered space. 
For every $i \in \zed,$ a  nilpotent endomorphism
$N$  of $(H,W)$  descends to
a nilpotent map $Gr_{i} N: Gr_{i}^{W}H \to Gr_{i}^{W}H$
and yields the   weight filtration $W^{Gr_{i}N}$ on $Gr_{i}^{W}H.$ 

\n
Let $k \in \zed $ be fixed. There is at most one filtration $W'$ of $H,$
called the {\em weight$-k$ filtration of $N$ relative to $W,$} such that
i) $N W'_{i} \subseteq W'_{i-2}$ and
ii) $W' ( Gr^{W}_{k} H ) = W^{Gr_{k}N}.$ 
See \ci{stzu}.

\smallskip
Let $(H,S)$ be as above, $N$ and $M$ be commuting nilpotent 
infinitesimal automorphisms of $(H,S).$
By the convolution formula  (\ref{convform}) for $W^{N},$  one has
$MW^{N}_{j} \subseteq W^{N}_{j}.$ 
For ease of notation we denote the map induced by
$M$ on $Gr^{N}_{j} H$ 
simply by $M.$

\n
Assume that $W^{M}[j]$ is the weight$-j$ filtration of $M$
relative to $W^{N}$ on $H$  for {\em every} $j \in \zed.$
In particular, this means that
\begin{equation}
    \label{morwf}
    M^{i}\, :\, Gr_{j+i}^{M} Gr_{j}^{N}H \, \simeq\,  Gr_{j-i}^{M} 
    Gr_{j}^{N}H,
    \qquad i \geq 0.
    \end{equation}
Set $P^{-j}_{-i}= \ke{ \, M}^{i+1} \cap \ke{\, N}^{j+1} \subseteq
Gr^M_{j+i}Gr^N_{j}H$
if $i, j \geq 0$ and zero otherwise.
We have the  double Lefschetz decomposition
\begin{equation}
    \label{adld}
    Gr_{j+i}^{M}Gr_{j}^{N}H \, = \, \bigoplus_{l,m \, \in \zed }{ M^{-i+l} 
    N^{-j+m} 
    P^{j-2m}_{i-2l}}, \quad i,j\, \in \zed.
    \end{equation}
The nondegenerate forms $S^{N}_{j}$ descend to
nondegenerate forms $S^{MN}_{ij}$ on $Gr^{M}_{j+i}Gr^{N}_{j}H.$
For $i, j \geq 0$
we have 
\begin{equation}
    \label{sijmn}
    S^{MN}_{ij}([a], [b]) \, = \, 
S(a, M^{i} N^{j} b),
\end{equation}
where  $a, b \in W^{M}_{j+i}H \, \cap \, W^{N}_{j}H$
are representatives of $[a],  [b] \in Gr_{j+i}^{M} Gr_{j}^{N}H.$
For the remaining  values of $i$ and $j,$ $S^{MN}_{ij}$ 
is defined by imposing that
(\ref{morwf}) is an isometry (cf. (\ref{snin})).

\n
The decomposition (\ref{adld}) is $S^{NM}_{ij}-$orthogonal.

\n
\begin{rmk}
\label{thenhodge}
{\rm
Let $n \in \zed$ be fixed 
and assume that the spaces $Gr_{j+i}^{M}Gr_{j}^{N}H$ are pure Hodge structures
of weight $(n-i-j)$ and that the induced maps
\begin{equation}
    \label{men1}
    N \, : \, Gr^{M}_{j+i} Gr^{N}_{j} H \lorw Gr^{M}_{j-2+i} 
    Gr^{N}_{j-2} H, 
    \end{equation}
    \begin{equation}
        \label{men2}
         M \, : \, 
     Gr^{M}_{j+i} Gr^{N}_{j} H \lorw Gr^{M}_{j+i-2} 
    Gr^{N}_{j} H.
    \end{equation}
are of pure type $(1,1).$ Then (\ref{adld}) is a direct sum
of pure Hodge substructures.
}
\end{rmk}

\begin{rmk}
    \label{uptosign}
    {\rm
If in addition  the form  $(-1)^{n-i-j}S^{MN}_{ij}$ is a polarization
of $P^{-j}_{-i}$ for  every pair of indices $(i,j) \neq (0,0)$ such that
$i, j \geq 0,$ then $(-1)^{i+j-m-l-1}S^{MN}_{ij}$ is a polarization
of  the summands
 $M^{-i+l} N^{-j+m} P^{j-2m}_{i-2l}$ in (\ref{adld}) except, possibly,
for $P^{0}_{0}.$ 
In this case, we simply say that the forms $S^{MN}_{ij}$ polarize
the summand spaces in question {\em up to sign}.
}
\end{rmk}

 \subsection{Filtrations  on $H^{*}(X)= \oplus_{j}{H^{j}(X)}$}
\label{asoh}
Let $f:X \to Y$ be a map of projective varieties, $X$ nonsingular,
$n =\dim{X},$ $\eta$ be an ample  line bundle on $X,$ $A$ be
an ample line bundle
on $Y$ and $L= f^{*}A.$ 
Let 
\begin{equation}
    \label{defofh}
    H^{*}(X)\,:\,=\, \oplus_{l}{ H^l(X)}
    \end{equation}
    and  define  the {\em twisted
Poincar\'e form} by setting:
 \begin{equation}
    \label{defofs}
 S  \left( \sum{\alpha_l}, \sum{\beta_l} \right)  \, :\,= \, \sum_{l}{ 
 (-1)^{\frac{l(l-1)}{2}} 
  \int_{X}{ \alpha_{l} \wedge \beta_{2n-l}} }.
  \end{equation}
  Note that the bilinear form $S$  on $H^{*}(X)$ is $(-1)^{n}-$symmetric, 
that $\eta$ and $L$ act via cup product
as nilpotent, commuting  operators on $H^{*}(X)$ with 
$\eta^{n+1}= L^{n+1}=0$ and that $\eta$ and $L$ are infinitesimal
automorphisms of $(H^{*}(X), S)$ (cf. \ref{fs}).

\n
The line bundles $\eta $ and $L$ act via cup product
on the cohomology of $X$  in a nilpotent fashion and  induce the weight
filtrations $W^{\eta}$ and $W^{L}$ on
$H^{*}(X).$

\smallskip
By the classical Hard Lefschetz Theorem \ref{chl} the
weight filtration $W^{\eta}$ for $\eta$ is the  filtration
by degree $W^{deg}:$
$$
W^{\eta}_i\,  =\, W^{deg} \, := \,  \bigoplus_{l \geq n-i} H^l(X).
$$
We also consider   
the {\em total} filtration on $H^{*}(X):$
$$
W_{i}^{tot} \, := \, \bigoplus_{b\in \zed}{  H^{n+b}_{ \leq b+i}(X)  }.
$$
Clearly,
\begin{equation}
\label{grgr}
Gr_{j+i}^{\eta} Gr_{j}^{tot} H^{*}(X) \, = \, H^{n-i-j}_{-i}(X), \quad 
i, \, j \, \in \zed
\end{equation}
and 
 (\ref{cupsh}) implies:
\begin{equation}
    \label{etaprtot}
  \eta \, W^{deg}_{i} \, \subseteq W^{deg}_{i-2}, \qquad  
  \eta \, W^{tot}_{j} \, \subseteq W^{tot}_{j}.
    \end{equation}
One of the main results of this paper is that $W^{L}=W^{tot},$
i.e. that, roughly speaking, the perverse filtration coincides with
the weight filtration induced by $L=f^{*}A$ 
(cf.  \ref{tiafiltr}).

\subsection{The Universal Hyperplane section and the 
defect of semismallness}
\label{hyperplane}
The universal hyperplane section plays 
a very important role in our inductive 
proof of the Relative Hard Lefschetz Theorem.
In this section we  prove that the defect of 
semismallness $r(f)$ of a map 
decreases when taking a map naturally associated with the universal 
hyperplane section. We also prove Weak Lefschetz-type results.

Let $X \subseteq {\Bbb P}$ be an embedded quasi-projective
variety and consider the 
 the universal hyperplane section diagram
$$
\xymatrix{
{\cal X}=\{ (x,s) \, | \; s(x) =0  \} \ar@{^{(}->}[r]^(.65){i} \ar[rd]^g & 
X \times \pd  \ar[d]^{f'} \\
     & {\cal Y}= Y \times \pd }.
$$
Let
$j: (X \times \pd ) \setminus{\cal X} \to  X \times \pd$ be
the open 
embedding. The morphism   $u:= f' \circ j$
is affine.
\medskip
In the special case $X= \p$ with $f=Id_{\p}$
we get
${\cal P} := \{ (p,s) \, | 
\:
s(p)=0 \} \subseteq \p \times \pd$  
which is  nonsingular  of codimension one and for 
which
the natural projection ${\cal P} \lorw \p$ is  smooth.

\n
A stratification ${\frak X}$  of $X$
with strata $S_{l}$ induces a 
stratification
on $X \times \pd$ with strata $S_{l} \times \pd$.
The following is elementary and left to the reader.
See \ref{iccso} and Lemma \ref{ugus}.b.

\begin{pr}
\label{uhssok}
The embedding $i: {\cal X} \lorw  X \times \pd$ is transversal
with respect to the stratification induced
by any stratification of $X,$ i.e.
the intersection  ${\cal P} \cap  (X \times \pd) = {\cal X}$ 
is transversal along every stratum $S_{l} \times \pd$ of $X \times 
\pd.$ In particular,
$$
IC_{\cal X} \,  \simeq  \,
i^{*}
IC_{X \times \pd} [-1] \,  \simeq \, i^{*} {p'}^{*} IC_{X}[d] [-1] .
$$
\end{pr}

Let $f: X \to Y$   be as in 
\ref{tpcss}.
We  recall the definition of the 
 defect of semismallness of the map $f$.  It plays a crucial role
in Goresky-MacPherson's version of the Weak Lefschetz Theorem in \ci{g-m}.
Set  $Y^{i} = \{ y \in Y \,  | \; 
\dim{f^{-1} (y)} = i \}.$

\begin{defi}
\label{defr}
{\rm The {\em defect of semismallness} of the map $f$
 is the integer
 $$
 r\, =\, r(f) :=\,  \max_{i| Y^i \neq \emptyset} \left\{ 2i+\dim{Y^i}-\dim{X} 
 \right\} .
 $$
 }
\end{defi}

Note that $r(f) \geq 0.$ If $r(f)=0$, then we say $f$ is {\em semismall}.
Note that this implies that $f$ is generically finite.
If $r(f)=0$ and the maximum is realized only for
$i=0$, then we say that $f$ is {\em small}.

\begin{rmk}
\label{sstri}
{\em
Let $f$ be as in \ref{tpcss}. If $f$ is semismall,
then $\fxn{f}{X}{n}
\simeq \phix{0}{\fxn{f}{X}{n}}$ and  theorems \ref{tm1}.a.b
hold trivially. In fact, 
 $r(f)=0$ implies that $\fxn{f}{X}{n}$
satisfies the conditions of support of Remark
\ref{refpss} (cf. \ci{b-m}). The conditions of co-support are automatic
since $\fxn{f}{X}{n}$ is self-dual.
}
\end{rmk}

The geometric quantity $r(f)$ plays a crucial role in our
proof by induction. The key point is that if it is not
zero, then it decreases by taking  hyperplane
sections.

\begin{tm}
\label{goesdown}({\bf $r(f)$ goes down})

\n
(a) If $r(f) >0$, then $r(g) < r(f).$ 

\n
(b) If $r(f) =0$, then $g$ is small.
\end{tm}
{\em Proof.}
For $s \in  {\Bbb P}^{\vee},$ let
$X_s:=  \{x \in X \, | \;  s(x)=0 \}$
be the  
corresponding  hyperplane section. If $(y,s) \in {\cal Y},$
then 
the projection 
$p: {\cal X} \to X$ identifies $g^{-1}(y,s)$ with $f^{-1}(y) \cap 
X_s$.
Set 
$$
{\cal Y}^{i \, '}= \left\{(y,s) : \dim{f^{-1}(y)}=i= 
\dim{f^{-1}(y)\cap X_s} \right\}.
$$
The point $(y,s) \in  {\cal Y}^{i \, '}$   
if and only if $X_s$ contains a top dimensional component of  
$f^{-1}(y)$.
It is a closed algebraic subset of ${\cal Y}^i.$
Set 
$$
{\cal Y}^{i \, ''}\,= \, \left\{(y,s)\, | \;  \dim{f^{-1}(y)}=i+1 \; \mbox{\rm 
and}
\; 
\dim{f^{-1}(y)\cap X_s} =i \right\}.
$$
It is an open algebraic subset of ${\cal Y}^i$.
We have that
$$
{\cal Y}^{i \,} ={\cal Y}^{i \, '} \amalg {\cal Y}^{i \, ''}.
$$
Since the set of hyperplanes in a projective space 
containing a given 
irreducible subvariety of dimension $d$ is  a linear space of 
codimension at least $d+1$, the definition of $r(f)$ implies that
$
\dim{ {\cal Y}^{i \, '}}  \leq $ $\dim{Y^i} + \dim{{\Bbb P}^{\vee}}
- (i-1) \leq $ $r(f) -2i + \dim{X} + \dim{{\Bbb P}^{\vee}}
- (i-1) =
r(f) -3i + \dim{\cal X}.$
It follows that
$$
2i + \dim{ {\cal Y}^{i \, '}} - \dim{\cal X} \leq r(f) -i, \quad
\forall i \geq 0.
$$
Since the general hyperplane section does not contain any
irreducible component of $f^{-1}(y),$ we have that
$\dim{{\cal Y}^{i \, ''}}=\dim{Y^{i+1}}+
\dim{{\Bbb P}^{\vee}}$ $\leq r(f) + -2(i+1) +\dim{X}+\dim{{\Bbb 
P}^{\vee}}=$
$r(f)-1+ -2i + \dim{{\cal X}}$. 
It follows that
$$
2i + \dim{ {\cal Y}^{i \, ''}} - \dim{\cal X} \leq r(f) -1, \quad
\forall i \geq 0.
$$
Suppose that either $Y^0$ is empty, or 
$ \dim{Y^0} - \dim{X} < r(f).$ Then the first inequality above
is  strict for $i=0.$ Combining it with the second inequality,
we get that $r(g) \leq r(f)-1.$

\n
Suppose that $Y^0$ is not empty and that $ \dim{Y^0} - \dim{X} = 
r(f).$ 
Then $r(f)=0$. The two inequalities above give
$r(g) \leq r(f)$, hence $r(g)=0$.
Moreover, $\dim{ {\cal Y}^i }  - 2i + \dim{\cal X} <0,$ $\forall
\, i >0$  so that $g$ is small.
\blacksquare

\medskip
A similar argument, 
based on ``Hironaka's principle of counting constants,'' 
as explained in \ci{so}, proves the following proposition, 
left to the reader:

\begin{pr}
\label{sohicc}
Let $X$ be nonsingular and $\eta$ be an ample line bundle on $X.$
There exists  $m_{0} \gg 0$ such that for
every $m \geq m_{0}$, having denoted by $X^{k}$ the transversal
intersection of $k$ general hyperplane sections
in  $|m \eta|,$ $k \geq 1$ and by 
$f_{k}: X^k \to Y$ the resulting morphism, we have :

\n
(a) If $r(f) \geq k$, then $r(f_{k}) \leq  r(f) -k.$ 

\n
(b) If $r(f) =0$, then $f_{1}: X^{1} \to Y$ is small.
\end{pr}

%\subsection{Weak-Lefschetz-type results}
%\label{wltpc}

The left $t-$exactness
of affine maps has important implications for the topology
of algebraic varieties.
\begin{lm}
 \label{wlmech}
 ({\bf Left $t-$exactness and Weak Lefschetz})
 Let 
% \begin{equation}
%     \label{wldiag}
%     \begin{array}{ccccc}
%X' & \stackrel{i}\lorw & X & \stackrel{j}\longleftarrow & X\setminus 
%X' \\
%& \searrow g & \downarrow f & \swarrow u & \\
%&& Y &&
%          \end{array}    
%      \end{equation}

 \begin{equation}
     \label{wldiag}
     \xymatrix{
X' \ar[r]^{i} \ar[dr]^g & X \ar[d]^f &  X\setminus X' \ar[dl]^u \ar[l]_j \\
& Y & &}
              \end{equation}
be a commutative diagram of algebraic varieties with $i$ a closed 
embedding, $f$ proper, $u$ affine and let $P \in Perv(X).$ Then

\smallskip
\n
(i) the natural map $\phix{l}{ f_{*}P} \to \phix{l}{ g_{*}i^{*} P}$
is iso for  $ l \leq -2$ and  mono for $l =-1;$

\smallskip
\n
(ii) the natural map $  \phix{l}{ g_{*}i^{!} P} \to \phix{l}{ f_{*}P}$
is  iso  for $ l \geq 2$ and  epi for $l =1.$
\end{lm}
{\em Proof.} By  applying $f_{*}\simeq f_{!}$ to the  triangle
$ j_{!}j^{!}P \to P \to i_{*}i^{*}P \to$ and by using
the isomorphisms
$u_{!} \simeq f_{!}j_{!},$ $j^{!} \simeq j^{*},$ $f_{*}i_{*} \simeq 
g_{*}$
one  gets the  triangle 
$ u_{!} j^{*}P \to f_{*}P \to g_{*}i^{*}P \to.$
Since $j^{*}P \in Perv (X \setminus X')$ and $u_{!}$ is left 
$t-$exact, (i) follows by taking the long exact sequence of 
perverse cohomology on $Y.$

\n
(ii) is obtained by first applying  (i) to ${\cal D}(P) \in Perv(X)$
and then by applying the rules
${\cal D} \phix{l}{-} \simeq \phix{-l}{ {\cal D}(-)},$
${\cal D} f_{*} {\cal D} \simeq f_{!} \simeq f_{*}$
and ${\cal D} g_{*} i^{*} {\cal D}  \simeq g_{!}i^{!}{\cal D} {\cal D}
\simeq g_{*} i^{!}.$
\blacksquare

\medskip
We shall need the following immediate consequence of 
Lemma  \ref{wlmech}.
\begin{pr}
\label{pwl} 
Let $Y$ be a projective variety, $i: Y_{1} \to Y$ be a
hyperplane section
and $P \in Perv (Y)$. Then the natural  maps
$$
i^*\, : \, \ihixc{j}{Y}{P} \lorw \ihixc{j}{Y_{1}}{  i^{*}P},
\qquad \qquad
i_* \, : \, \ihixc{l}{Y_{1}}{i^{!}P} \lorw \ihixc{l}{Y}{ P}
$$
are  isomorphisms for $j\leq -2$ and $l \geq 2,$
 injective for $j=-1$ and surjective for $l=1.$
\end{pr}
{\em Proof.} Apply Lemma \ref{wlmech} and take hypercohomology.
\blacksquare

\medskip
Consider  the universal hyperplane section 
in \ref{hyperplane}.
There is the  commutative diagram
$$
\xymatrix{
X \ar[dd]_f & & X \times \pd \ar[ll]_{p'} \ar[dd]^{f'}  \\
 & {\cal X} \ar@{^{(}->}[ur]^i \ar[dr]^{g} & & X \times \pd 
 \setminus {\cal X} \ar[ul]^j \ar[dl]^u \\
Y && {\cal Y}=Y \times \pd \ar[ll]_p 
}
$$
%$$
%\begin{diagram}
%{\cal X}  &       \rInto^{i}     & X \times \pd     & \lInto^j & X \times \pd \setminus {\cal X} \\
%          &  \rdTo               &   \dTo^{f'}           & \ldTo  \rdTo^u       &              \\ 
%          &                      &{\cal Y}=Y \times \pd &           &        X   \\
%         &                      &                  &       \rdTo^p  & \dTo    \\
%         &                      &                  &               &         Y
%\end{diagram}.
%$$
%Setting
%$j: (X \times \pd ) \setminus{\cal X} \to  X \times \pd$ be
%the open 
%embedding,  
where $u:= f' \circ j$
is an 
affine morphism.

\begin{pr}
\label{wltphc} ({\rm {\bf The Relative Weak Lefschetz Theorem}})

\n
Let $K \in Perv(X),$ 
$K':= p'^{*} K [d],$ $M:= i^{*} K'[-1].$
Then

\n
(i)
$p^*\phix{l}{f_* K }[d]
\lorw 
\phix{l+1}{g_*M }
$
is  iso for $ l \leq -2$
and mono for $l=-1;$

\n
(ii)
$\phix{l-1}{g_*M }\lorw 
p^*\phix{l}{f_* K }[d]$
is  iso for $ l\geq 2$ and epi
for $ l=1.$
\end{pr}
{\em Proof.} Since $p'^{*}[d]$ is $t-$exact,
(i) follows from Lemma \ref{wlmech}.i
applied to $K'$ in  the set-up
of \ref{hyperplane},
keeping in mind
that
$f'_{*}{p'}^{*} 
\simeq p^{*} f_{*}$  
and
$p^{*}[d] (\phix{l}{f_* K } )   \simeq ) \, 
\phix{l}{p^{*}[d] (f_* K) }.$

\n
By transversality, \ref{uhssok}
and  Lemma \ref{ugus}.b, $i^{!}K' \simeq i^{*}K'[2]$
and
(ii) follows from Lemma \ref{wlmech}.ii.
\blacksquare

\medskip
The following complements Proposition
\ref{wltphc}; see \ci{bbd}, 5.4.11.
\begin{pr}
\label{maxi}
$p^*[d] \phix{-1}{f_* K }$ can be identified with the biggest 
perverse subsheaf of 
$\phix{0}{g_*M }$ coming from $Y,$ and  
$p^*[d] \phix{1}{f_* K }$ with the biggest quotient perverse sheaf
 of $\phix{0}{g_*M }$ coming from $Y.$
\end{pr}

\section{Proof of the main theorems in \ref{tpcss}}
\label{pomt}
In this section, we prove Theorem \ref{tm1}, i.e.
the Relative Hard Lefschetz Theorem, the Decomposition Theorem
and the Semisimplicity Theorem (except for $\phix{0}{\fxn{f}{X}{n}}$),
the Hard Lefschetz Theorem for Perverse Cohomology
\ref{tm3}, the Hodge Structure Theorem \ref{uf}, the 
$(\eta, L)-$Decomposition Theorem \ref{etaldecompo}
and the Generalized Hodge-Riemann Bilinear Relations \ref{tmboh}. 

\medskip
The set-up is the one of \ref{tpcss}.
We assume \ref{basicass} and remind the reader of Remark 
\ref{ovviamento}.

\subsection{Relative Hard Lefschetz, Decomposition
and  Semisimplicity ($i \neq 0$) }
\label{prhlt}
We use the notation and results in \ref{hyperplane}.
Let  $\eta':=  i^{*}{ p'}^{*} \eta;$ it is $g-$ample.

\begin{lm}
\label{rhl1}
Suppose that 
${\eta'}^r:
\phix{-r}{g_*M }\stackrel{ \simeq }\lorw
\phix{r}{g_*M }$ for all $r \geq 0$ and that 
$\phix{0}{g_*M }$ is semisimple. Then 
$\eta^r:
\phix{-r}{f_*K }\stackrel{ \simeq }\lorw
\phix{r}{f_*K }$ for $r \geq 0.$
\end{lm}
{\em Proof.}
We have  $\eta^{r} = i_{*} \circ {\eta'}^{r-1} \circ  i^{*}.$
If $r \neq 1,$ then we conclude by
the Weak Lefschetz Theorem  \ref{wltphc}. 

\n
Let $r=1.$
Since $p^*[d]$ is fully faithful, $\eta$ is an isomorphism if and 
only if 
$p'^*\eta[d]:
p^*\phix{-1}{f_*K }[d]\lorw
p^*\phix{1}{f_*K }[d]$ is an isomorphism. This map 
is the composition of 
the monomorphism  $i^{*}$
%$
%p^*\phix{-1}{f_* K }[d]
%\stackrel {i^*}\lorw 
%\phix{0}{g_*M } 
%$ 
with the epimorphism $i_{*}.$
%$
%\phix{0}{g_*M }\stackrel{i_*}\lorw 
%p^*\phix{1}{f_* K }[d].
%$
%Suppose $\ke \, p'^*\eta[d] \neq 0$. 
By the semisimplicity of $\phix{0}{g_*M },$
the sequence of perverse  subsheaves
$
i^* \ke \, p'^*\eta[d] \subseteq \ke\, i_* \subseteq \phix{0}{g_*M }
$
splits and we get  a direct sum decomposition
$$
%\ke i_*=i^*\ke \,p'^*\eta[d] \oplus R,
%\qquad
\phix{0}{g_*M }\, = \, i^*\ke \, p'^*\eta[d] \oplus R \oplus S,
$$
where the restriction of $i_*$ to $S$ is an isomorphism with 
$p^*\phix{1}{f_* K }[d]$. 

\smallskip
\n
The projection 
$
\phix{0}{g_*M } \lorw i^*\ke \, p'^*\eta[d]\oplus S
\simeq i^*\ke \, p'^*\eta[d]\oplus p^*\phix{1}{f_* K }[d]
$
is an epimorphism  and both summands come from $Y.$
By the maximality statement of Proposition \ref{maxi},
$ i^*\ke \, p'^*\eta[d]=0,$ i.e
. 
$\eta:
\phix{-1}{f_*K }\lorw
\phix{1}{f_*K }$
is a monomorphism.

\n
Since the complex  $K$ and the map $ \eta$  are
self-dual, $\eta$ is also an epimorphism, hence an isomorphism.
\blacksquare

\begin{rmk}
\label{mple}
{\rm 
The relative Hard Lefschetz Theorem \ref{tm1} 
holds  when  $\eta$  is $f-$ample, i.e.
ample when restricted to the fibers of $f.$
In fact, if $\eta$ is $f-$ample, then 
$\widetilde{\eta}:= \eta + mL$ is ample for every $m \gg 0$
and,  by Remark \ref{sepullback},
$\phix{j}{\eta} = \phix{j}{\widetilde{\eta}}$, for every $j \in \zed.$
Note that   ample implies $f-$ample.
We need Theorem
\ref{tm1} for $f-$ample line bundles in the next proposition.
}
\end{rmk}

\medskip
We can now prove Theorem \ref{tm1} parts (a), (b) and (c),
except for $i=0.$

\begin{pr}
\label{rhl}
Let 
$
f : X \to  Y
$
and $\eta$ 
be as in  \ref{tpcss} and
assume that
\ref{basicass} holds.

\n
Then  the Relative Hard Lefschetz Theorem \ref{tm1}.a 
and the Decomposition Theorem \ref{tm1}.b 
hold for $f.$ The Semisimplicity Theorem \ref{tm1}.c   holds for  
$\phix{i}{ \fxn{f}{X}{n} }$ with $i\neq0$.
\end{pr}
{\em Proof.} 
We apply the inductive hypothesis \ref{basicass} to $g: {\cal X} \lorw 
{\cal Y}$, 
which satisfies $r(g)<r(f)$ (cf. \ref{goesdown}). 
Setting $K=\rat_X[n],$ we have $M=\rat_{\cal X}[n+d-1].$
By the inductive hypothesis and Remark
\ref{mple}, we have:
1) ${\eta'}^{r}$ is an isomorphism for every $r$  and 2)
$\phix{0}{g_{*}M}$ is semisimple.
By Lemma \ref{rhl1}, 
$\eta^{r}$ is an isomorphism for $r \geq 0.$
This proves that Theorem \ref{tm1}.a holds for $f.$
As it has  already  been observed in 
\ref{serhl},
the well-known Deligne-Lefschetz Criterion for $E_{2}-$degeneration
\ci{dess} yields  Theorem \ref{tm1}.b for $f.$ 
The semisimplicity statement b) for $i\neq 0$ follows  from the
Weak Lefschetz Proposition \ref{wltphc}
and the semisimplicity of $\phix{0}{g_{*}M}.$
\blacksquare

\subsection{Hard Lefschetz  for perverse cohomology groups}
\label{phltpcq}
Note that  $r(f) \leq  \dim{X}.$                            
If $r(f) >0,$ consider $ 1 \leq k \leq r(f),$  let
$X^k,$ 
be the transversal intersection of $k$ general hyperplane sections
of the linear system $\eta$ and  $f_k: X^k \lorw Y$ be the resulting map.
By  Proposition
\ref{sohicc}, we may assume that
$\eta$ is such that if $r(f) >0,$ then
$r(f_k) \leq r(f) -k,$
for every $  1 \leq k \leq r(f)$.

\medskip
The following  contains  Weak-Lefschetz-type results
for the hyperplane sections of $X.$

\begin{pr}
\label{wltspl}
Assumptions as in \ref{phltpcq}. Let $ 1 \leq k \leq r.$

\n
(i) The 
natural  restriction map
$
\phix{l-k}{\fxn{f}{X}{n}} \lorw \phix{l}{ \fxn{f_{k}}{X^{k}}{n-k}} 
$
is iso for $l < 0$ and a splitting mono for $l=0.$

\n
(ii) The
natural Gysin map
$
\phix{l}{ \fxn{f_{k}}{X^{k}}{n-k}} \lorw \phix{k + l}{\fxn{f}{X}{n}}
$
is an iso for $l>0$ and a
 splitting epi for $l=0.$

 \n
 (iii)
The maps induced by $L$ in hypercohomology are compatible with
the splittings (i) and (ii). 
In particular,  if for given $h$ and $j$ the map 
$$
L^{j} : 
H^{h}_{0}(X^{k})
\lorw H^{h+2j}_{0}(X^{k}) 
$$
is injective (resp. surjective, resp. bijective),
then the maps
$$
L^{j}: H^{h}_{-k}(X)
\lorw H^{h+2j}_{-k}(X), \qquad  L^{j}: H^{h+2k}_{k}(X)
\lorw H^{h+2k+2j}_{k}(X), 
$$
are injective (resp. surjective, resp. bijective).
\end{pr}
{\em Proof.} Let $k=1.$ Lemma \ref{wlmech}
applied to $X^{1} \subseteq X$ and $\rat_{X}[n]$ implies (i) and (ii). The
case $k>1$ follows by induction.

\n
Since cupping with $L=f^{*}A$ commutes with any direct sum decomposition 
on $Y,$ (iii) follows.
 \blacksquare

\medskip
We need the following easy consequence of Proposition
\ref{wltspl} to prove Theorem \ref{tm3}.
The more precise statement $\ke{\,L_{n-1}} \subseteq
H^{n-1}_{\leq -1}(X)$ is true, but its proof would
require
that we prove 
Theorem \ref{tm3} first.

\begin{lm}
\label{kln-1} Assumptions as in \ref{phltpcq}.
$$
\ke{\,L_{n-1}} = \ke{\, L_{n-1, \leq 0}} \subseteq H^{n-1}_{\leq 0}(X).
$$
\end{lm}
{\em Proof.}
Since $L$ acts compatibly with the  p-splitting,
 it suffices to prove that
$\ke{L_{n-1}} \subseteq H^{n-1}_{k} (X)$ is trivial,
for every $k >0.$
By Proposition \ref{wltspl}, $H^{n-1}_{k} (X)$ is isomorphic to 
a direct summand
of $H^{(n-k) - (k+1)}_{0}(X^{k})$ with $L$ acting
as the restriction of  $L_{|X^{k}}$ to the direct summand. 
Inductively, this map
is injective. We conclude by  Proposition \ref{wltspl}.iii.
\blacksquare

\begin{pr}
\label{prtm3}
Under the assumption \ref{basicass}
the Hard Lefschetz theorem 
for perverse cohomology groups  
\ref{tm3} holds for $f,$ i.e.  
$$
\eta^{k}\; : \;  H^{j}_{-k}  (X)
\, \simeq \,
H^{j+2k}_{k}(X),  \quad
L^{k}\; : \;  H^{n+b-k}_{b}  (X)
\,\simeq \,
H^{n+b+k}_{b}(X), \quad k\geq 0,\;\; b, \, j \, \in \zed.
$$
\end{pr}
{\em Proof.} 
Since $f_* \rat_X[n]$ p-splits by \ref{rhl}, there is a decomposition
$$
H^{n+j}_{\leq b}(X) \,
\stackrel{\f}{\simeq}
\, 
\bigoplus_{i \leq b}{\Bbb H}^{j-i}(Y, 
   \phix{i}{\fxn{f}{X}{n}} ). 
$$
The statement for $\eta^{k}$ follows from the
previously established Relative Hard Lefschetz
Theorem for $f$ (cf. \ref{rhl}) and from the compatibility 
\ref{compa}. 

\n
The rest of the proof is concerned with $L^{k}.$
The cases $b \neq 0$ follow from the inductive hypotheses:
apply Theorem  \ref{tm3} and 
Proposition \ref{wltspl}.iii to $f_{|X^{k}} : X^{k} \to Y.$

\n
Let $b=0$.  Note that 
the statement is trivial for $k=0.$
Choose a sufficiently general hyperplane section $Y_{1}$ of $Y$
(cf. \ref{cutstr}) and let $X_{1}$ be the nonsingular
$f^{-1}(Y_{1}).$

\n
Recalling the canonical identification
$H^{n+k}_{0}(X)= \ihixc{-k}{Y}{\phix{0}{\fxn{f}{X}{n}}}$
and the compatibility \ref{sepullback},
we have a commutative diagram
$$
\xymatrix{
H^{n-k}_{0}(X) \ar[d]^{i^*} \ar[r]^{L^{k}}
& H^{n+k}_{0}(X)  \\
H^{(n-1)-(k-1)}_{0}(X_{1})
\ar[r]^{L_{|X_{1}}^{k-1}}  &
H^{(n-1)+(k-1)}_{0}(X) \ar[u]^{i_*}
}
$$
where $i^*$ is restriction 
%given the fact that $\phix{0}{\fxn{f}{X}{n}}$ is self-dual, 
and $i_*$ is the (dual) Gysin map.

\n
Let $k \geq 2.$ By
 Proposition \ref{pwl}, $i^*$ and $i_*$ are isomorphisms.
 By 
 Theorem \ref{tm3} applied to   $f_1: X_1 \lorw Y_1,$ the map
$L^{k-1}_{|X_{1}}$ is an isomorphism and so is the map 
 $L^{k}.$

\n
Let $k=1.$ 
We must check that  
$L: H^{n-1}_0(X)  \to  H^{n+1}_0(X)$
is an isomorphism. 

\n
By the self-duality of $\phix{0}{\fxn{f}{X}{n}},$
  the two spaces have the same dimension so that
  it is enough to check that $L$ is injective.

  \smallskip
  \n
Let $\alpha \in \ke \, L  \subseteq H^{n-1}_{0}(X).$ According
to \ref{serhl}, $\alpha = \sum_{j\geq 0}{ \eta^{j} \alpha_{j}}.$
By Remark \ref{sepullback}, $L$ commutes with any direct sum 
decomposition
so that $L \eta^{j} \alpha_{j}=0, $ for every $j\geq 0.$
Since $L$ commutes with $\eta$ and, by what we have already
proved for $\eta^{k},$ the map $\eta^{j}$ is injective
on $H^{n-2j -1}_{-2j}(X) \ni \alpha_{j}$ for every $j \geq 0,$
we have that $L\alpha_{j}=0$ for every $j \geq 0.$
The case $b \neq 0,$ which was dealt-with above, implies that
if $j >0,$ then $L$ is an isomorphism from
$H^{n-2j -1}_{-2j}(X).$ This implies that $\alpha_{j}=0,$ for every
$j >0, $ i.e. $\alpha = \alpha_0.$

\n
Since $L$ acting on $H^{*}(X)$ is strict and $\ke\, L \subseteq 
H^{n-1}_{ \leq 0}(X)$ 
by Lemma \ref{kln-1}, we have
$\ke\, L / ( \ke\, L \cap  H^{n-1}_{\leq -1} ) = 
\ke{ 
\{ 
\,
H^{n-1}_{0}(X) \stackrel{L}\to H^{n+1}_0(X) \,
\} 
}.$
It follows that there exists $a \in \ke \, L \subseteq
H^{n-1}_{\leq 0}(X)$ such that its class in $H^{n-1}_{0}(X)$ satisifes
$[a] = \alpha.$

\n
Let $a = \sum{a^{pq}}$ be the $(p,q)-$decomposition for the natural
Hodge structure on $H^{n-1}(X).$ Clearly $a^{pq} \in \ke \, L.$
It follows that we may assume that $\alpha = [a],$
with $a \in \ke \, L $ of  some pure type $(p,q).$

\n
 By way of contradiction, assume that $\alpha \neq 0.$ By Lemma 
 \ref{dtpb}, we have
 $i^{*}\alpha = [ a_{ | X_{1} } ].$ Since $i^{*}$ is injective 
 by Proposition \ref{pwl}, $0 \neq i^{*} \alpha \in H^{n-1}_{0}(X_{1}).$
 This restricted class
 is of pure type $(p,q),$ for the pure Hodge structure coming
 inductively from 
 Theorem \ref{uf} applied to $f_{|X_{1}} : X_{1} \to Y_{1}.$
 By Lemma \ref{dtpb}: i) $L_{|X_{1}} i^{*} \alpha= i^{*}(L\alpha)$ and
 ii) since $({\cal P}_{\eta})_{|X_{1}} \simeq {{\cal 
 P}_{\eta}}_{|X_{1}} [1],$ we have $\eta_{|X_{1} } i^{*} \alpha =0.$
Since $\alpha = \alpha_0,$   $i^{*} \alpha \in P^{0}_{0}(X_{1}).$
 
 \n
 By the inductive Generalized Hodge-Riemann Relations
 \ref{tmboh} for $f_{ |X_1 }: X_1 \to Y_1$
 $$
 0 \, \neq  \,  \pm S^{\eta_| L_| }_{00}( i^*\alpha, 
 \overline{ i^*\alpha } ) \, = \,  \int_{X_1}{  a_{|X_{1}} \wedge
 \overline{a}_{ |X_1 } } \, = \, \int_{X}{ L \wedge a \wedge 
 \overline{a} } \, = \, 0,
 $$
a contradiction. 
\blacksquare

\bigskip
Proposition \ref{prtm3} allows  to complete the identification
of the  filtrations defined in \ref{fs}:

\begin{pr}
    \label{tiafiltr}
For every  $i,j \in \zed:$ 
\begin{equation}
    \label{repetitaL}
 W^{\eta} \, = \, W^{deg}, \quad 
W^{L}  \,= \,  W^{tot},
\end{equation}
\begin{equation}
\label{givesht}
W^L_i  \, = \, \bigoplus_{l \in \zed}H^{n+l}_{\leq l+i}(X),
\quad
Gr^L_i  \, = \, \bigoplus_{l \in \zed}H^{n+l}_{ l+i}(X),
\quad
W^{L}_{i} \, \cap \, H^{n+j}(X) \, = \, H^{n+j}_{ \leq i+j}(X),
\end{equation}
\begin{equation}
    \label{shofpf}
Gr^{\eta}_{j+i}Gr^{L}_{j} H^{*}(X) \, = \, H^{n-i-j}_{-i}(X)
\end{equation}
and the forms $S^{\eta L}_{ij}$ of (\ref{sijmn})
are therefore defined on $H^{n-i-j}_{-i}(X).$

\n
The filtration $W^{\eta}[j]$ is the weight-j filtration
of $\eta$ relative to $W^{L}$ (cf. (\ref{morwf})).
\end{pr}
{\em Proof.}  
Since  $L$
is compatible with the splitting $H^{*}(X)= \oplus_{l}{ H^{l}(X)},$
the convolution formula (\ref{convform}) implies
that
$$
W_i^{L}\,=\,  \bigoplus_l{ (W_i^L \cap H^l(X))}.
$$
By the characterization of  weight filtrations,
in order to prove that  $W^{deg}=W^{\eta}$ and $W^{tot}=W^{L}$, 
it is enough to show that
i) $\eta \, W^{deg}_{k} \subseteq W^{\deg}_{k-2},$ 
ii) $\eta^{k}: H^{n-k}(X) \simeq H^{n+k}(X),$ 
iii) $L \, W^{tot}_{k} \subseteq W^{tot}_{k-2},$
and
iv) $L^{k}: Gr^{tot}_{k} \simeq Gr^{tot}_{-k},$ for all $k \geq 0.$

\n
While i) is obvious  and iii) follows from  \ref{sepullback},
ii) and iv) are mere re-formulations
of  the
 Hard Lefschetz Theorems \ref{chl} for $\eta^{k}$
 acting on ordinary cohomology
 and of  the just-established Hard Lefschetz for Perverse Cohomology
 Theorem \ref{tm3} for the action of  $L^{k},$
 respectively. 
 This proves (\ref{repetitaL}), (\ref{givesht}) and (\ref{shofpf}).
 The assertion on $S^{\eta L}$  follows from   (\ref{shofpf}) and 
 the definition (\ref{sijmn}).

\n
The second assertion follows from $W^{L}=W^{deg}$, i) above and 
the Hard Lefschetz Theorem \ref{tm3}
for $\eta.$ 
\blacksquare

\medskip
\n
{\bf {\em Proof of the Hodge Structure Theorem \ref{uf}.}}
Since $L$ is of pure type $(1,1)$ acting on $ H^{*}(X),$
the convolution formula
(\ref{convform}) for $W^{L}$ implies that
$W^{L}_{i} \cap H^{n+j}(X)$ is a Hodge sub-structure
of $H^{n+j}(X)$ for every $i,\,j \in \zed.$
We conclude by (\ref{givesht}).
\blacksquare.

\medskip
\n
{\bf {\em Proof of the $(\eta,L)-$Decomposition 
Corollary \ref{etaldecompo}.}}
By Proposition \ref{tiafiltr},
$W^{\eta}[j]$ is the weight$-j$ filtration
of $\eta$ relative to $W^{L}$ and
\ref{fs} applies.
\blacksquare

\medskip

\subsection{The
Generalized Hodge-Riemann Bilinear Relations ($P^{i}_{j} \neq P^0_0$)}
\label{eccettopoo}
In this section we are going to prove
the Generalized Hodge-Riemann Bilinear Relations
 \ref{tmboh} for $f,$ except for $P^{0}_{0}.$
 As in \ref{phltpcq}, we may assume that the $\eta-$sections
 are general enough.

\begin{lm}
\label{chepalle}
(a) 
Let $r>0$  and $X^{r}$ be the complete transversal intersection of
$r$ general sections of $ \eta.$
Then the natural restriction map 
$
i^{* } :  H^{n -r-j}_{-r}(X)  \lorw  H^{n-r-j}_0(X^{r}),
$
is an injective map of pure Hodge structures for every $j \in \zed,$
and $i^{*} (P^{-j}_{-r}(X) ) \subseteq P^{-j}_0( X^{r} )$ for every
$ j \geq 0.$

\n 
(b)
Let $ X_{1} =f^{-1}(Y_{1}),$
where $Y_{1}$ is a
general section of $A,$ 
transversal to the strata of $Y.$
Then, for every $j >0,$  the natural restriction map
$
i^{*}: H^{n-j}_0(X)  \lorw H^{n-j}_0(X_{1})
$
is an  injection of  pure  Hodge structures
and $i^{*}( P^{-j}_0(X)) \subseteq P^{-j+1}_0( X_{1} ).$
\end{lm}
{\em Proof.} 
The inductive hypotheses apply to
$X^{r} \to Y$ and to $X_{1} \to Y_{1}$
so that all perverse cohomology groups
have natural Hodge structures.

\n
(a) We have proved Theorem \ref{uf} for $f.$ 
The map $H^{*}(X) \to H^{*}(X^{r})$ 
is a map of Hodge structures  and 
Remark \ref{filtrfac} implies that so is the map in question.
 Lemma \ref{wltspl}.i. implies the injectivity 
statement.
The fact that $P^{-j}_{-r}(X)$ maps
to $P^{-j}_0(X^{r})$ can be shown as follows.
Let $l > 0$ and $l' \geq 0$ and let $X^l$
be the transversal complete intersection
of $l$ general sections of $ \eta.$
The map $f_* ( \eta^{l+l'}): \fxn{f}{X}{n} \lorw \fxn{f}{X}{n+2l+2l'}$
factorizes as
$$
\fxn{f}{X}{n} \lorw \fxn{g}{X^l}{n-l}[l] 
\stackrel{g_* ( \eta^{l'}_{|X^l} ) }\lorw
\fxn{g}{X^l}{n-l} [l] [2l'] \lorw 
\fxn{f}{X}{n+2l+2l'}.
$$
The statement follows from applying
the functors ${\Bbb H}^*( \phix{-r}{-}   )$ to the factorization above
when $l= r$ and $l'=1.$

\n
(b) The compatibility with the perverse decomposition
of Lemma \ref{dtpb} implies that
the map in question has the indicated range.
The map $H^{*}(X) \to H^{*}(X_{1})$ 
is a map of Hodge structures and  so is the map in question.
The injectivity statement 
follows from Proposition \ref{pwl}. The fact that $P^{-j}_0(X)$ maps
to $P^{-j+1}_0(X_{1})$ follows from the fact
that $L^{j} =   
i_* \circ L^{j-1}_{| X_{1}} \circ i^*$, where $i_*$ is an 
isomorphism
in the range we are using it.
\blacksquare

\begin{pr}
 \label{ufcind}
 The 
Generalized Hodge-Riemann Bilinear Relations
 \ref{tmboh} hold for the direct summands
 $\,
 \eta^{-i+l} \, L^{-j+m} \, P^{j-2m}_{i-2l} \neq P^{0}_{0}.
 $
 \end{pr}
{\em Proof.}
By Remark \ref{uptosign}, it is enough to
consider the case $ i, j \geq 0,$ $(i,j) \neq (0,0).$
Also, we assume that $l=m=0,$ and leave the easy necessary modifications
to deal with the other cases to the reader.

\n
 Let $X^r_{s}$ be the complete intersection
 of $r$ general sections of $\eta$ and $s$
 general sections of $L.$ If $r=0,$ then we consider only the
 sections of $L.$ Similarly, if $s=0.$
 
 \n
 Since $\int_{X}{ \eta^{r} \wedge L^{s}\wedge a \wedge b} = 
 \int_{X^{r}_{s}}{ a_{|X^{r}_{s}} \wedge b_{|X^{r}_{s}}},$ we have
 that 
 $$
 S_{rs}^{\eta L}(X)(a,b) \, = \, S_{00}^{\eta_{|} L_{|}}
 ( X^{r}_{s})  ( a_{|X^{r}_{s}}, 
 b_{|X^{r}_{s}}).
 $$
 
 \n
 The statement follows from a repeated application 
 of Lemma \ref{chepalle}, of 
 the inductive hypothesis
 Theorem  \ref{tmboh} applied to $X^{r}_{s} \to Y_{s}$ and 
 from  Remark \ref{important}.
 \blacksquare

\subsection {The space $\Lambda \subseteq H^n(X),$ its approximability 
 and  the polarization of $P^0_0$}
 \label{tapwat}
In this section we are going to complete
the proof of the Generalized Hodge-Riemann Bilinear Relations
 \ref{tmboh} for $f,$ by polarizing  $P^{0}_{0}.$
For convenience, we consider cohomology with real coefficients.

\medskip
Let $\e >0$ be a real number. Define
$$
\Lambda_{\e} := \ke \, (\e\, \eta +L) \subseteq H^{n}(X).
$$
The spaces $\Lambda_{\e}$ are Hodge sub-structures.
By the classical Hard Lefschetz Theorem, 
$\dim{\Lambda_{\e}}= b_n - b_{n-2},$ where $b_{i}$ are the Betti numbers of 
$X.$
Define
$$
\Lambda\,:= \, \lim_{\e \to 0} \Lambda_{\e},
$$
where the limit is taken in the Grassmannian $
G(b_n - b_{n-2}, H^n(X)).$
The space $\Lambda \subseteq H^n(X)$ is a real Hodge sub-structure
for that is a closed condition.
Clearly, $
\Lambda \subseteq \ke \{H^n(X) \to H^{n+2}(X) \},
$
but in general there is no equality, since 
$$
\dim  \ke\,  \left\{ H^n(X) \stackrel{L}\to H^{n+2}(X) \right\}
\, = \, b_n-b_{n-2}+
\dim  \ke\,\left\{H^{n-2}(X) \stackrel{L} \to H^{n}(X) \right\}.
$$

The main goal of this section is to characterize the subspace 
$\Lambda$ in terms of
$\eta$ and $L.$

Since $L^{k}: Gr_{k}^{L} \simeq Gr_{-k}^{L},$
we have 
\begin{equation}
   \label{clqeivl614}
   \ke\,  L^k \subseteq W^L_{k-1}.
   \end{equation}
In order to  keep track of cohomological  degrees, we 
 set 
$$
L^k_r \,:\, Gr_r^{\eta} \,  = \, H^{n-r}(X)  \lorw  H^{n-r+2k}(X)
\, =  \,
Gr_{r-2k}^{\eta}.
$$

\medskip
The following two lemmata will allow to identify $\Lambda.$

\begin{lm}
\label{sliceke} 
$\;
\eta \,  \ke\,  L_2\,  \cap\,  (\eta\,  \ke\, L_{2})^{\perp} \,
\cap\,  \cdots\,  \cap\,  (\eta^i\,  
\ke\,  L^i_{2i})^{\perp}\, =\, \{0\} \, \in H^{n}(X), \;  i \gg 0.  
$
\end{lm}
{\em Proof. } It is enough to show   that 
$$
\eta \, \ke \, L_2\,  \cap \, (\eta\, \ke\,  L_{2})^{\perp}\,
\cap\,  \cdots\,  \cap\,  (\eta^i \, 
\ke\,  L^i_{2i})^{\perp} \, = \,
\eta \, \ke\,  L_2 \, \cap\,  W_{-i}^L , \qquad \forall i \geq 0,
$$
for then the  lemma follows by taking $i = n+1,$ for example.

\n
The claim above can be  proved by induction as follows.
The starting step $i=0$ of the induction follows from 
 (\ref{clqeivl614}): 
 \begin{equation}
\label{cesvp}
\ke \, L_2\, \subseteq \, W^L_{0}.
\end{equation}
Suppose the claim proved for $i.$ 
We are left with showing that
\begin{equation}
\label{fave}
\eta \, \ke\,  L_2 \, \cap \, W_{-i-1}^L\, 
=\, \eta\,  \ke\,  L_2\,  \cap\,  W_{-i}^L\, \cap \,
(\eta^{i+1} \,  \ke\,  L^{i+1}_{2i+2})^{\perp}.
\end{equation}
The inclusion ``$\subseteq$'' follows at once
from (\ref{clqeivl614}) and self-duality:
$
\eta^{i+1} \ke\, L^{i+1}_{2i+2} \subseteq W_i^L=(W_{-i-1}^L)^{\perp}.
$
The other inclusion follows from the  nondegeneracy of 
the forms 
$S^{\eta L}_{i+2,i},$ as we now show.

\n
Let $\alpha  = \eta \lambda \in \eta \, \ke\, L_2 \cap W_{-i}^L.$ 

\smallskip
\n
{CLAIM 1:} $\lambda \in W^{L}_{-i}.$ 

\n
We have $\lambda \in W^{L}_{0}.$ By way of contradiction,
assume that $\lambda \in W^{L}_{-i'},$ for some $ -i < -i' \leq 0.$
Since $\eta \lambda \in W^{L}_{-i},$ we have $\eta \lambda
\in H^{n}_{\leq -i}(X)$ and we would have that the map
$$
\eta\, :\, Gr^{\eta}_{2} \, Gr^{L}_{-i'}\,  =\, H^{n-2}_{-i'-2}(X)
\lorw Gr^{\eta}_{0}\, Gr^{L}_{-i'} \, = \, H^{n}_{-i'} 
$$
is not injective, contradicting Proposition \ref{tiafiltr},
i.e. the injectivity of  $\eta$  for $i \geq -1.$

\smallskip
\n
CLAIM 2: there exists $\lambda' \in W^{\eta}_{2i+2} \cap W^{L}_{i}=
H^{n-2-2i}_{ \leq -2-i}(X)$ such that $\lambda = L^{i} \lambda'.$

\n
We have that $L^{k}: Gr^{L}_{k} \simeq Gr^{L}_{-k}$ for every $k \geq 0.$
Using the case $k=i$ we may write
$$
 \lambda \, = \, L^{i} \lambda_{1} + \tau_{1}, \qquad
 \lambda_{1} \in H^{n-2-2i}_{ \leq -2-i}(X), \; \; \tau_{1} \in 
 H^{n-2}_{\leq -2-i-1}(X).
$$
Replacing $\lambda$ with $\tau_{1},$ $k=i$ with $k=i+1$ and iterating
we get
$$
\lambda \, = \, L^{i}\, \sum_{t=1}^{j}{ L^{t-1}\lambda_{t} } + \tau_{j},
\qquad  L^{t-1}\lambda_{t} \in H^{n-2-2i}_{ \leq -2-i}(X), 
\; \; \tau_{j} \in H^{n-2}_{ \leq -2-i-j}(X).
$$
CLAIM 2 follows by taking $j \gg 0.$

\smallskip
\n
Since $L_2 \lambda=0$, we have $L^{i+1}_{2i+2} \lambda' =0.$
So far, we have proved that
$$
\hbox{ if $\alpha \in \eta \, \ke \, L  \cap W^{L}_{-i}, \;\; $ then
$\alpha = \eta L^{i}_{2i+2} \lambda' $ 
with $L^{i+1}_{2i+2} \lambda' =0.$  }
$$
Let $\beta \in \eta^{i+1} \ke \, L^{i+1}_{2i+2}.$ By 
(\ref{clqeivl614}), $\beta = \eta^{i+1} \beta',$ for some
$\beta' \in \ke \, L^{i+1}_{2i+2} \subseteq W^{L}_{i}.$
By the  very definition (\ref{sijmn}) of the forms $S^{\eta L}$
it follows that,   denoting  by the same symbol
an element in some $H^{*}_{\leq *}$
and the corresponding  class in $H^{*}_{*}(X),$
$$
S^{\eta L}_{ij} ( \lambda', \beta') \, = \, 
S( \eta L^{i} \lambda' , \eta^{i+1} \beta') \, = \, 
S(\alpha, \beta).
$$
Finally, let us assume that 
$\alpha \in \eta \, \ke \, L_2 \, 
\cap W^{L}_{-i} \cap (\eta^{i+1} \ke \, L^{i+1}_{2i+2}  )^{\perp}.$ 
This implies $S(\alpha , \beta)=0.$

\n
By the Generalized Hodge-Riemann Bilinear Relations \ref{tmboh},
the restriction of the form $S^{\eta L}_{i+2,i} $ to 
$\ke\, L^{i+1} \subseteq H^{n-2i-2}_{-i-2}(X)$ is nondegenerate.
It follows that the class
of $\lambda'$ in $H^{n-2i-2}_{-i-2}(X)$ is zero, i.e. $\lambda' \in
W^{L}_{i-1}.$ Since $\eta$ respects the filtration  $W^{L}$ while
$L$ shifts it by $-2,$ we conclude that
$\alpha = \eta L^{i} \lambda' \in W^{L}_{-i-1}.$ We have proved the 
 remaining inclusion. 
\blacksquare

\begin{lm}
\label{w}
$\;
\Lambda \, = \,
\ke\, L_0 \cap \left(  \cap_{i\geq 1} (\eta^i \ke \, L_{2i}^i )^{\perp}   
\right) 
\subseteq H^n_{\leq 0}(X)
$

\n
and there is the orthogonal direct sum decomposition:
$$
\ke{\,  L_0} \, = \,  \Lambda\, \bigoplus \, \eta\, \ke{\, L_2}.
$$
\end{lm}
{\em Proof.}
We show that 
$$
\Lambda_{\e} \subseteq \bigcap_{i\geq 1} ( \eta^i \ke \, L_{2i}^i     )^{\perp}.
$$
It is enough to show that
$\Lambda_{\e} \subseteq ( \eta^i \ke \, L_{2i}^i )^{\perp}, $ for every $i\geq 1.$
Let $u_{\e} \in \Lambda_{\e}$, i.e. $\eta u_{\e}= - \e^{-1} L u_{\e}$. We 
have
$\eta^i u_{\e} = (- \e^{-1} )^i L^i u_{\e}.$
Let $\lambda \in H^{n-2i}(X)$  be such that $L^i \lambda =0$.
We have  $\int_X u_{\e}\wedge \eta^i  \lambda= 
(-\frac{1}{\epsilon} )^{i}\int_X u_{\e}\wedge L^i \lambda =0.$
The wanted inclusion follows and $\Lambda \subseteq 
 \cap_{i\geq 1} ( \eta^i \ke L_{2i}^i     )^{\perp}$.

\n
We show that 
$
\Lambda \subseteq \ke{\, L_0}:
$
if $\Lambda \ni u = \lim_{\e \to 0} u_{\e},$ with $u_{\e} \in \Lambda_{\e},$
then
 $Lu= \lim_{\e \to 0} L u_{\e}=
\lim_{\e \to 0} ( -\e \eta u_{\e} )=0$.

\n
It follows  that
$
\Lambda \subseteq 
\ke{\, L_0} \cap (  \cap_{i\geq 1} ( \eta^i \ke \, L_{2i}^i     )^{\perp}  ).
$

\n
By Lemma \ref{sliceke},
$\cap_{i\geq 1} ( \eta^i \ke L_{2i}^i     )^{\perp} \cap \eta{\ke{\, L_{2}}}
= \{0 \}$ and therefore
$
 \Lambda \cap \eta \ke L_2 = \{ 0 \}.
$
By counting dimensions, the internal direct sum
$\Lambda \oplus \eta \,\ke{\, L_2} = \ke{\, L_0}.$
On the other hand, we also have an internal direct sum
$(\ke{\, L_0} \cap (\cap_{i\geq 1} ( \eta^i \ke\, L_{2i}^i     )^{\perp}))
\oplus 
\eta \ke \, L_2 \subseteq \ke \, L_0$
and this implies that the inclusion $\Lambda \subseteq 
\ke{\, L_0} \cap (\cap_{i\geq 1} ( \eta^i \ke{\, L_{2i}^i}     )^{\perp})$
is in fact an equality. The orthogonality of the decomposition is immediate.
\blacksquare

\medskip
The form $S_{00}^{\eta L}$  is nondegenerate
on each direct summand of the $(\eta,L)-$decomposition
for $H^{n}_{0}(X).$ In particular, it is so on 
$P^0_0.$

\begin{lm}
\label{pwo}
The form  $(-1)^{n}S^{\eta L}_{00}$ defines a polarization
of
$\Lambda_0 := \Lambda / 
(\Lambda \cap \lij{H}{n}{\leq -1} (X))$.
\end{lm}
{\em Proof.}
By the classical Hard Lefschetz Theorem,
the Poincar\'e pairing multiplied by $(-1)^{\frac{n(n+1)}{2}}$
is a polarization of $\Lambda_{\e}$ for every $\e >0.$
In particular, the form  $(-1)^{n}S( -, C(-))$
is semipositive definite when restricted to
$\Lambda.$ 

\n
It follows that $(-1)^{n}S_{00}^{\eta L}(-,C(-)),$ being
 semipositive definite and nondegenerate  on $P^{0}_{0},$
 is in fact 
 positive 
definite,
i.e. 
$(-1)^{n}S^{\eta L}_{00}$ is a polarization of $P^{0}_{0}.$
\blacksquare

\begin{lm}
\label{finita}
We have an orthogonal direct sum decomposition:
$$
\ke L_0 /  (\ke L_0 \cap H^n_{\leq-1}(X))=
\Lambda_0 \bigoplus \left(\eta \ke L_2 / \eta \ke L_2 \cap 
H^n_{\leq-1}(X)\right).
$$
\end{lm}
{\em Proof.} The statement follows from the following elementary fact: 
if $V$ is a vector space with a bilinear form 
and $V_1 \subseteq V$ is its radical, an orthogonal direct sum decomposition
$$
V = U_1  \oplus U_2
$$
induces an orthogonal direct sum decomposition
$$
V/V_1=U_1 / (V_1 \cap U_1) \oplus U_2 / (V_1 \cap U_2)
$$
and the bilinear form
is nondegenerate on the two summands. 
We apply this to $V= \ke{\, L_0},$ $V_1= \ke{ L_0} \cap H^n_{\leq-1}(X),$
$U_1=\Lambda$ and $U_2= \eta \, \ke{\, L_2}.$
\blacksquare

\bigskip
We  now conclude the proof of the Polarization 
Theorem \ref{tmboh} for $P^{0}_{0}:$

\medskip
\n
{\bf {\em Proof of Theorem \ref{tmboh}.}}
Since $P^{0}_{0} \subseteq \ke \,\eta,$ we have an inclusion
of Hodge structures
$$
P^{0}_{0} \subseteq (\eta \ke{\, L_2})^{\perp}/(\eta \ke{\, L_2})^{\perp}
\cap H^n_{\leq -1}(X) =\Lambda_0
$$ 
 which, in view of Lemma \ref{pwo} and Remark \ref{important}, are 
 polarized by $(-1)^{n}S_{00}^{\eta L}(-,C(-)).$
\blacksquare

  \section{The Semisimplicity Theorem \ref{tm1}.c
  for $\phix{0}{\fxn{f}{X}{n}}.$}
  \label{sios}
  The set-up is as in \ref{tpcss} and we assume
  that
  \ref{basicass} holds. Recall Remark \ref{ovviamento}.

\bigskip
In this section we prove 
that  $\phix{0}{\fxn{f}{X}{n}}$
is semisimple, i.e. we establish
the remaining case $i=0$ of 
Theorem \ref{tm1}.c 
for $f.$ Along the way we also prove the Generalized Grauert
Contractibility Criterion \ref{nhrbr} and
the Refined Intersection Form Theorem  \ref{rcffv}, thus proving 
all the results
in \ref{tpcss}.

\subsection{The induction on the strata: reduction to $S_{0}.$}
\label{defws}
We  introduce the stratification with which we  work.
 Let $K \in Ob (D(X))$. The typical example will be
  $K=\rat_{X} [n]$.  
  We fix once and for all
  ${\frak X}$ and $\frak Y$ finite algebraic Whitney stratifications
  for $f$ such that $K$ is ${\frak X}-$cc.
  By \ref{stabi} and \ref{perverse},
  $f_{*}K$ and all of its perverse cohomology
  complexes 
  $\phix{j}{ f_{*}K }  $, $\forall \,j \in \zed,$
  are ${\frak Y}$-cc.
  
\medskip
 We employ the notation in  \ref{wsoav}.
Let $ 0 \leq s \leq d.$ Denote by 
$$
 \; S_{s} \; \stackrel{\alpha_{s}}\lorw U_{s} 
 \stackrel{ \beta_{s}}\longleftarrow \; U_{s+1}  
$$
the corresponding closed and open embeddings.

\n
The stratification ${\frak Y}$ induces 
a  stratification ${\frak Y}_{U_{s}}$ on $U_{s}$ and 
the trivial one, ${\frak Y}_{S_{s}},$ on $S_{s}.$ 
The  maps $\alpha_{s}$ and $\beta_{s}$ are 
 stratified with respect to these stratifications.

 \n
 Let $K'$ be ${\frak Y}-$cc, e.g. $K'=f_{*}K$
 or $K'= \phix{j}{f_{*}K},$ $ l \, \in \zed.$
 Then $\alpha_{s}^{*}K'$ is ${\frak Y}_{S_{s}}-$cc
 and $\beta_{s}^{*}K'$ is ${\frak Y}_{U_{s+1}}-$cc.
 
\n 
Let $K'  \in Ob (D( U_{s}  ) )$ be $\frak{Y}_{U_{s}}-$cc.
By \ref{stabi},
all terms of the   triangle
$
{\alpha_{s}}_{!} \alpha_{s}^{!} K' \to K' \to {\beta_{s}}_{*} 
\beta_{s}^{*}K' \stackrel{ [1] }\to
$
are ${\frak Y}_{U_{s}}-$cc and the  maps
induced at the level of cohomology sheaves
are, when restricted to the strata  $S_{l},$ $l\geq s,$  
maps of local systems. 

\n
Let $n: =\dim{X},$  $m:=\dim{f(X)}.$
The stratum $S_{m}$ has a unique connected component
$S_{f}$
contained in the open subset of $f(X)$ over which $f$
is smooth.

\n
Clearly,  all the complexes we shall be interested in  have support 
contained in $f(X).$ In addition,  depending on whether
they are defined on $Y$, $U_{s}$ or $S_{s}$, they are either ${\frak 
Y}$-cc,
${\frak Y}_{U_{s}}-$cc, or ${\frak Y}_{S_{s}}-$cc.

\begin{rmk}
\label{h-s}
{\rm 
By the condition of (co)support in \ref{perverse}, we have 
$$
\phix{j}{\fxn{f}{X}{n}}_{|U_{s}}\, \simeq \,  \tu{-m} \, \td{-s} \,
\phix{j}{\fxn{f}{X}{n}}_{|U_{s}}, \qquad 
 \forall \;
0  \leq s \leq m.
$$
The sheaf
$ \csix{-s}{  \phix{j}{ \fxn{f}{X}{n}  }_{| U_{s}}    }$
is a local system on $S_s.$
}
\end{rmk}
Let  $f_{s} : U'_{s}:= f^{-1} (U_{s}) \lorw U_{s}$ be the 
corresponding
maps.
Note that $U'_{s} =\emptyset$, $\forall s > m.$
We have natural  restriction isomorphisms
$$
\phix{j}{ \fxn{f}{X}{n}_{| U_{s}   } } \simeq
\phix{j}{ \fxn{f_{s}}{ U'_{s} }{n} } . 
$$
Recall that,  if $P \in Perv (U_{s+1})$, then ${\beta_{s}}_{!*}P
\simeq \td{-s-1}P \in Perv (U_{s} )$ (cf. \ref{iee}).

\n
In this set-up, 
Deligne's Theorem \ci{dess}
 can be re-formulated in terms of the existence
of an isomorphism
$$
\fxn{f_{m}}{ U'_{m} }{n} \simeq \dsdjix{j}{j}{\fxn{f_{m}}{ U'_{m} 
}{n} }
$$
where  $\phix{j}{ \fxn{f_{m}}{U_{m}'}{n} }$ is supported,
as a complex on $U_{m}$,
precisely on $S_{f}$  and is  there isomorphic to
$(R^{n-m+j} {f_{m}}_{*} \rat_{ f^{-1} ( S_{f} ) } )
[m]$.

\begin{rmk}
\label{tassok}
{\rm
The local systems 
$R^{n-m+j} {f_{m}}_{*} \rat_{ U_{m}'} $
on $S_f$
are semisimple  by Deligne Semi\-sim\-pli\-city Theorem \ref{dss}. In 
particular,
the complexes $\phix{j}{ \fxn{f_{m}}{U_{m}'}{n} }$
are semisimple in $Perv(U_{m})$.
}
\end{rmk}

\medskip
The following Lemma essentially 
reduces the proof of the missing part of the Decomposition Theorem
to the local criterion of Lemma \ref{redtopt}.b.
To prove that the local  criterion is met 
we reduce it to  a global property of projective maps, Theorem 
\ref{pwo}.

\begin{lm}
\label{redtopt}

\n
(a) For every $(j,s) \neq (0,0)$ we have a
canonical isomorphism
in $Perv(U_{s})$:
$$
\phix{j}{ \fxn{f}{X}{n}  }_{| U_{s}}
\simeq 
{\beta_{s}}_{*!} ( \phix{j}{ \fxn{f}{X}{n}  }_{| U_{s+1}}  )
\oplus \csix{-s}{  \phix{j}{ \fxn{f}{X}{n}  }_{| U_{s}}    } [s],
$$
where the projection to the first summand
is the (unique) lifting of truncation and the projection
to the second  stems from truncation
(cf. \ref{h-s}).

\n
(b)
For $(j,s)=(0,0)$
$$
\phix{0}{ \fxn{f}{X}{n}  }
\simeq 
{\beta_{0}}_{*!} ( \phix{0}{ \fxn{f}{X}{n}  }_{| U_{1}}  )
\oplus \csix{0}{  \phix{0}{ \fxn{f}{X}{n}  }    } [0]
$$
if and only if the natural map of dual skyscraper sheaves (cf.
\ref{fullyf})
$$
\csix{0}{ { \alpha_{0}  }_{!} 
{ \alpha_{0}  }^{!} \,\phix{0}{ \fxn{f}{X}{n}  } }
\lorw
{ \alpha_{0}  }_{*} { \alpha_{0}  }^{*}
\csix{0}{ \phix{0}{ \fxn{f}{X}{n}  } }
$$
is an isomorphism.    
\end{lm}
{\em Proof.} 
(a)  The perverse
sheaf $\phix{j}{ \fxn{f}{X}{n}}$ is semisimple
for  $j \neq 0$
by Theorem \ref{tm1}.c for $f.$
We  apply the Splitting Criterion
 \ref{splitp} whose hypotheses are met
in view of Remark \ref{condmet}.

\n
Let $j=0.$ 
$\phix{0}{ \fxn{f}{X}{n}}$ is self-dual
 by Poincar\'e-Verdier duality.
By Remark \ref{condmet}, it is enough
to check that the Splitting Criterion  \ref{splitp}
holds for $ 1 \leq s \leq m.$
In the case $s=m,$ $\phix{0}{\fxn{f}{X}{n}}_{|U_{m}}$ is a shifted 
local system and there is nothing to prove.
Let $1 \leq s \leq m-1.$
Let $Y_{s} \subseteq Y$  be the complete intersection
of $s$  hyperplane sections chosen so that 
1) it meets {\em every} connected component of
the pure and positive dimensional $S$ transversally at
a finite set $T$ and  2) $X_{s}: = f^{-1} (Y_{s})$ is a  nonsingular
variety (cf. \ref{cutstr}).
We obtain a projective morphism
$f_{s}: X_{s}  \lorw Y_{s}$. 
We have $\dim{X} > \dim{X_{s}}$ and we can apply our inductive 
hypotheses: Theorem \ref{tm1}.b and c hold and
$\phix{0}{ \fxn{f_{s}}{X_{s}}{n-s}} $ is semisimple.
By Lemma \ref{dtpb}:
$\phix{0}{ \fxn{f_{s}}{X_{s}}{n-s} } \simeq
\phix{0}{  \fxn{f}{X}{n}}_{|Y_{s}}[-s].$

\n
The semisimplicity of
$\phix{0}{\fxn{f_{s}}{X_{s}}{n-s}}$ implies, via Remark
\ref{condmet}, that the conditions for the splitting criterion
for $\phix{0}{\fxn{f_{s}}{X_{s}}{n-s}}$
of Lemma \ref{splitp} are met at every point
of $T$ which is a subset of the set of  zero-di\-men\-sio\-nal strata
for $f_{s}$.

\n
By the second part of Lemma \ref{adapt}, we have that
the splitting condition for
$\phix{0}{\fxn{f}{X}{n}}_{|U_{s}}$ is met as well.

\n
(b) Since we have the result for $U_{1},$ 
the statement is a mere re-formulation of Lemma \ref{splitp}.
\blacksquare

\subsection{The local system $\csix{-s}{\alpha^{!}_{s} \phix{0 }{ 
\fxn{f}{X}{n}  } }$ on $S_{s}$}
\label{tlsa}
 \begin{lm}
\label{surjinpcds}
Let $Z$ be  an affine algebraic variety,  $Q \in  D^{\leq 0}(Z)$, 
i.e. $\dim{ \hbox{supp} ( {\cal H}^{i}(Q))} 
\leq -i.$ Let $\alpha: \Sigma \lorw Y$ be the closed
embedding of the possibly empty support of ${\cal H}^{0}(Q).$

\n
Then 
the natural restriction map below is surjective 
$$\ihixc{0}{Z}{Q} \lorw \ihixc{0}{Z}{\alpha_{*}
\alpha^{*}Q}.
$$ 
\end{lm}
{\em Proof.}
We have the two spectral sequences
$E^{pq}_2 (Q) = \ihixc{p}{Z}{ {\cal H}^q (Q) }
 \Longrightarrow \ihixc{p+q}{Z}{ Q }$,
$E^{pq}_2 ( \alpha_{*} \alpha^{*} Q) = 
\ihixc{p}{Z}{ {\cal H}^q ( \alpha_{*} \alpha^{*} Q) }
 \Longrightarrow \ihixc{p+q}{Z}{ \alpha_{*} \alpha^{*}  Q }$.
The natural adjunction map
$a : Q \lorw \alpha_* \alpha^* Q$ induces a map of spectral sequences
$E_r(Q) \lorw E_r ( \alpha_* \alpha^* Q   ).$

\n
Note that $\dim{ \Sigma} \leq 0.$ The assumptions
on $Q$ imply that 
$ E^{pq}_2 ( \alpha_{*} \alpha^{*} Q)= 0$ if either
$p \neq  0,$ or $q>0$
so that
$E_2 (\alpha_{*} \alpha^{*} Q    )
= E_{\infty} ( \alpha_{*} \alpha^{*} Q ).$  
In particular,
 $\ihixc{0}{Z}{\alpha_{*}\alpha^{*}Q}=
\lijx{E}{00}{\infty}{\alpha_{*}\alpha^{*}Q}
=\lijx{E}{00}{2}{\alpha_{*}\alpha^{*}Q} .$

\n
Since $Q \in D^{\leq 0}(Z),$
$Q$ is ${\frak Z}-$cc with respect to some
stratification ${\frak Z}$ of $Z$, $\mbox{supp} {\cal H}^{q}(Q)$
is a closed affine subset of $Z$
of dimension at most $-q.$
The theorem on the cohomological dimension of affine sets
with respect to constructible sheaves,
\ci{k-s}, Theorem 10.3.8, implies that
$E^{pq}_2(Q)=0$ for every $p+q >0.$

\n
We have $E^{pq}_2(Q) = E^{pq}_{\infty}(Q) =0$
if either $p <0$ or $p+q >0.$

\n
It follows that we have the surjection
$$
\ihixc{0}{Z}{Q} \lorw E^{00}_{\infty}(Q)
= E^{00}_{\infty} ( \alpha_{*}
\alpha^{*}Q   )= \ihixc{0}{Z}{\alpha_{*}
\alpha^{*}Q}. 
$$
\blacksquare

\medskip
\n
In what follows, by the conditions of support for perverse sheaves,
$\mbox{supp} {\cal H}^0 (\phix{b}{\fxn{f}{X}{n} })$ is either empty,
or a finite set of points. In the first case, Proposition
\ref{pcdsurj} is trivial.

\begin{pr}
\label{pcdsurj}
Let $b \in \zed$ and  $\alpha$ be the closed embedding  into $X$
of the
zero-di\-men\-sio\-nal set $\mbox{supp} {\cal H}^0 (\phix{b}{\fxn{f}{X}{n} })=
\{y_1, \cdots,y_r\}.$
The  restriction  map
$$
 H^{n+b}_b(X) =
\ihixc{0}{Y}{   \phix{b}{\fxn{f}{X}{n}}}\lorw 
\ihixc{0}{Y}{ \alpha_*\alpha^*\phix{b}{\fxn{f}{X}{n}} }=
\bigoplus_{i=1,\cdots,r}H^{n+b}_b(f^{-1}(y_i)) 
$$
is surjective.
Dually, the cycle map below is injective
$$
\bigoplus_{i=1,\cdots,r} H^{BM}_{n+b,-b}(f^{-1}(y_{i}) ) 
\lorw 
H^{n-b}_{-b}(X). 
$$ 
\end{pr}
{\em Proof.} 
Let $U \subseteq Y$ be an affine open set such that  $\mbox{supp}
\csix{0}{ \phix{b}{\fxn{f}{X}{n} } } \subseteq U'$ and $U':=
f^{-1}(U).$
Consider the commutative diagram

$$
\xymatrix{
H^{n+b}(X) \ar[r]^A & H^{n+b}(U')  \ar[r]^{B}
& H^{n+b} (f^{-1}(y) ) \\
\lij{H}{n+b}{\leq b}(X)\ar[u] \ar[d]^{p_1} \ar[r]^{A_{\leq b}}
& \lij{H}{n+b}{\leq b}(U') \ar[u] \ar[d]^{p_2} \ar[r]^{B_{\leq b}}
& \lij{H}{n+b}{\leq b}(f^{-1}(y)) \ar[u] \ar[d]^{p_3}\\
\lij{H}{n+b}{ b}(X) \ar[r]^{A_{ b}} \ar[d]^{\stackrel{\f}{\simeq}}
& \lij{H}{n+b}{ b}(U') \ar[r]^{B_{ b}} \ar[d]^{\stackrel{\f}{\simeq}} 
& \lij{H}{n+b}{b}(f^{-1}(y))  \ar[d]^{\stackrel{\f}{\simeq}}      \\
\ihixc{0}{Y}{ \phix{b}{\fxn{f}{X}{n} }}  \ar[r] & \ihixc{0}{U}{ 
\phix{b}{\fxn{f}{X}{n} }} \ar[r] & {\cal H}^0
(\alpha_* \alpha^*  \phix{b}{\fxn{f}{X}{n} })
}
$$
%$$
%\begin{array}{ccccc}
%H^{n+b}(X) & \stackrel{A}\lorw & H^{n+b}(U') & \stackrel{B}\lorw
%& H^{n+b} (f^{-1}(y) ) \\
%\uparrow & & \uparrow & & \uparrow \\
%\lij{H}{n+b}{\leq b}(X) & \stackrel{A_{\leq b}}\lorw 
%& \lij{H}{n+b}{\leq b}(U') & \stackrel{B_{\leq b}}\lorw
%& \lij{H}{n+b}{\leq b}(f^{-1}(y)) \\
%\downarrow p_1 & & \downarrow p_2 & & \downarrow p_3 \\
%\lij{H}{n+b}{ b}(X) & \stackrel{A_{ b}}\lorw 
%%& \lij{H}{n+b}{ b}(U') & \stackrel{B_{ b}}\lorw
%& \lij{H}{n+b}{b}(f^{-1}(y)) \\
%\downarrow ||\f & & \downarrow ||\f & & || \f \\
%\ihixc{0}{Y}{ \phix{b}{\fxn{f}{X}{n} }} & \lorw & \ihixc{0}{U}{ 
%\phix{b}{\fxn{f}{X}{n} }}
%& \lorw & {\cal H}^0(\alpha_* \alpha^*  \phix{b}{\fxn{f}{X}{n} }).
%\end{array}
%$$
where the vertical maps pointing up are the natural injections,
the quotient maps 
 $p_i$ are surjective and the vertical maps on the bottom row
 are the identifications of Remark \ref{comincia}. 
In view of the existing splitting $\f,$ the maps $A$ and $B$ 
are strict with respect to  the perverse  filtrations
 on $H^{n+b}(X),$ $ H^{n+b}(U)$
 and the induced  filtration on $H^{n+b}( f^{-1} (y))$
 (cf.  \ref{collectinfo}).

\n
By Lemma \ref{surjinpcds}, $B_b$ is surjective. 
This  implies that $B_b \circ p_2 = p_3 \circ B_{\leq b}$ is surjective.

\n
By Deligne's Theory of Mixed Hodge Structures,
\ci{ho3}, Proposition 8.2.6, $\im B \circ A = \im B$. 
 By  the strictness 
with respect to the perverse and to the induced   filtration,
we infer that
$\im B_{\leq b} \circ A_{\leq b} = \im B_{\leq b}$

\n
It follows that
$p_3 \circ B_{\leq b} \circ A_{\leq b}= (B_b \circ A_b) \circ
p_1$ is surjective and so is $ B_b \circ A_b$, i.e.
we have proved the wanted surjectivity.
\blacksquare

\subsection{$\phix{0}{\fxn{f}{X}{n}}$ is a direct sum of intersection 
cohomology complexes}
\label{iadsoicc}

In this section we prove Proposition
\ref{itissum}, i.e. an important step towards the semisimplicity
of $\phix{0}{\fxn{f}{X}{n}}.$ 
A key ingredient is the the Generalized Grauert
Contractibility Criterion \ref{nhrbr}, which 
is concerned with the Hodge-theoretic properties
of the refined  intersection form
$H^{BM}_{n-b,b}(f^{-1}(y)) \to H^{n+b}_{b}(f^{-1}(y))$
introduced in \ref{sperem}.
Together with
Lemma \ref{collectinfo},  Theorem \ref{nhrbr}  gives complete information
on the structure of the refined intersection form 
on the fibers of $f.$

\medskip
\n
{\bf {\em Proof of the  Generalized Grauert
Contractibility Criterion \ref{nhrbr}.} }
The injectivity follows from the second statement
in Proposition \ref{pcdsurj}. 
If
$y \notin \hbox{supp} {\,\cal H}^0( \phix{b}{\fxn{f}{X}{n}} ),$
then $H^{BM}_{n-b,b}( f^{-1}(y) )=0$ and the injectivity statement
is trivial.

\n
The inclusion in $\ke \, L$ follows from the fact that one can find
a hyperplane section of $Y$ avoiding $y.$

\n
The class  map  $cl: H^{BM}_{n-b}( f^{-1}(y)) \lorw H^{n+b}(X)$ is a map
of mixed Hodge structures so that the image $\im (cl) \subseteq
H^{n+b}(X)$ is a pure Hodge sub-structure.

\n
By Lemma \ref{collectinfo}, the class map is filtered
and  $H^{BM}_{n-b}(f^{-1}(y)) =
H^{BM}_{n-b,\leq b}( f^{-1}(y) ).$

\n
It follows that the projection, 
$\im{ (cl_{b}) },$ of $\im{(cl)}$ to $H^{n+b}_{b}(X)$ is a pure
Hodge sub-structure.

\n
The compatibility with the direct sum decomposition given by 
Theorem \ref{tm1}.a and  \ref{serhl} follows from the
additivity of $\alpha_{!} \alpha^{!},$ where $\alpha: y \to Y.$

\n
By the Generalized Hodge-Riemann Bilinear Relations
\ref{tmboh}, the direct summands of $\im{(cl_{b}) }$ are
$S^{\eta L}_{-b0}-$orthogonal. Since $\im{(cl_{b})}$ is a 
Hodge
sub-structure of $H^{n+b}_{b}(X),$
the form
$S^{\eta L}_{-b0}$ induces a polarization on each direct summand
(cf. \ref{important}).
\blacksquare

\bigskip
The proof of Proposition \ref{itissum} requires only the case $b=0$ 
of Theorem \ref{nhrbr}.  
Consider the natural adjunction map
$$
A \; : \;   \alpha_! \alpha^! \phix{0}{ \fxn{f}{X}{n} }
\lorw  \phix{0}{ \fxn{f}{X}{n} }.
$$

\begin{pr}
\label{cacciucco}
The map
$$
{\cal H}^0(A)_y: {\cal H}^0 ( \alpha_! \alpha^! \phix{0}{ 
\fxn{f}{X}{n} }  )_y
\lorw 
 {\cal H}^0 ( \phix{0}{ \fxn{f}{X}{n} }  )_y
$$
is an isomorphism.
\end{pr}
{\em Proof.}
Since the domain and the  target have the same rank, it is enough
to show injectivity.
Let
$$
A' \; : \;   \alpha_! \alpha^! \phix{0}{ f_{*}\omega_{X}[-n]}
\lorw \phix{0}{ f_{*}\omega_{X}[-n]}
$$
be the natural adjunction map.
In view of Remark \ref{thepoint}, the statement
to  be proved is equivalent
to  the analogous statement for the
map
${\cal H}^0(A')_y.$
Consider the composition
$$
 I\,:\; \alpha_! \alpha^! \phix{0}{ f_{*}\omega_{X}[-n]}
\lorw \phix{0}{ f_{*}\omega_{X}[-n]} \simeq
\phix{0}{ f_{*}\rat_{X}[n]} \lorw
\alpha_{*} \alpha^{*} \phix{0}{ f_{*}\rat_{X}[n]}.
$$
By the self-duality of $I,$
the domain and target of  ${\cal H}^0(A')_y$ have the same rank.
The linear map ${\cal H}^0(I)_y$ is the 
refined intersection form
$H^{BM}_{n,0}(f^{-1}(y) \lorw H^{n}_{0}(f^{-1}(y) )$
which is an isomorphism by Theorem \ref{nhrbr}.
This implies that
${\cal H}^0(A')_y$ is injective
and hence an isomorphism.
\blacksquare

%Note that ${\cal D} ( \phix{0}{\e} ) = \phix{0}{\e}.$ 

%\n
%Consider the composition $ \phix{0}{\e} \circ A$ and its dual
%${\cal D}(A) \circ   \phix{0}{\e}.$

%\n
%The composition
%$$
%\begin{array}{ccc}
%\phix{0}{ \fxn{f}{X}{n} } & \stackrel{ \phix{0}{\e}}\lorw &
%{\cal D} ( \phix{0}{ \fxn{f}{X}{n} } ) \\
%A \uparrow &  &  
%{\cal D}(A) \downarrow  \\
%\alpha_! \alpha^! \phix{0}{ \fxn{f}{X}{n} }
%&
%\stackrel{B}\lorw  &
%{\cal D} ( \alpha_! \alpha^! \phix{0}{ \fxn{f}{X}{n} } )
%\end{array}
%$$
%defines  pairings
%$$
%\begin{array}{ccc}
%H^n_0(X) \times H^n_0(X) & \stackrel{\Phi^n_0}\lorw  & \rat \\
%\uparrow & & || \\
%H_{n,0}(f^{-1}(y) \times H_{n,0}(f^{-1}(y) ) 
%& \stackrel{B^n_0}\lorw  &\rat \\
%\end{array}
%$$
%where we have identified
%$\ihixc{0}{Y}{\alpha_! \alpha^! \phix{0}{\fxn{f}{X}{n}}}
%$
%with $H_{n,0}(f^{-1}(y),$ and ${\Bbb H}^0 (Y, \phix{0}{\e} )$ with
%$\Phi^n_0$ in view of Lemma \ref{dmatrix}.

%\n
%We have that
%${\cal H}^0(A)_y$ is an isomorphism
%if and only if ${\cal H}^0(B)_y =
%{\Bbb H}^0(Y, B) $ 
%is an isomorphism if and only if $B^n_0$ is non-degenerate.
%We conclude by Theorem \ref{nhrbr}.
%\blacksquare

\begin{pr}
\label{itissum}
There are  canonical isomorphisms in $Perv (Y)$ for every $b:$
$$
\phix{b}{\fxn{f}{X}{n}} \simeq
\bigoplus_{l=0}^{\dim{Y}}{ 
IC_{\overline{S_{l}}}( \alpha_{l}^{*} {\cal H}^{-l} (  
\phix{b}{\fxn{f}{X}{n}}  )).}
$$
\end{pr}
{\em Proof.} It follows from
Lemma \ref{redtopt} and Proposition \ref{cacciucco}.
\blacksquare

\medskip
\n
{\bf {\em Proof of the Refined Intersection 
Form Theorem \ref{rcffv}.}}
By Lemma \ref{collectinfo} we only need to deal with the case
$a=b.$ In this case, by
\ref{splitp}, the nondegeneracy of the refined intersection form
in question is precisely the obstruction to the splitting 
of $\phix{b}{\fxn{f}{X}{n}}$ so that the statement follows from
Proposition \ref{itissum}.

\blacksquare

\subsection{ The semisimplicity of
$\phix{0}{\fxn{f}{X}{n}}$}
\label{86}
The goal of this section is to prove Theorem
\ref{sss}. 
 The local systems in question do not seem to arise as the ones associated with 
the cohomology of the fibers of a smooth map, so that
Deligne's semisimplicity result  \ref{dss} does not apply directly.
The idea of the proof is to use hyperplane
sections on $Y$ to find a smooth projective
family  ${\cal X}_{T} \lorw T$  of $(n-s)-$dimensional
varieties over a Zariski-dense open  subset
$T$ 
of $S_s$ in a way that allows to use
Proposition \ref{pcdsurj} to infer that, over $T,$
the local system $H^{n-s}_{0}( {\cal X}_{t}) $
maps surjectively onto the local system
$H^{n-s}_{0}(f^{-1}(t)).$
The left-hand side is semisimple by Deligne's
Semisimplicity Theorem \ref{dss}.
It  follows that so is the right-hand side.
On the other hand, the latter  is the restriction
of ${\cal H}^{-s}( 
\alpha^{*}_{s}\phix{0}{\fxn{f}{X}{n}} )$ to $T$
and the semisimplicty over $S_{s}$ follows (cf. \ref{sszos}).

\medskip
We need a relative version of Proposition \ref{pcdsurj}.

\begin{lm}
\label{pcdsrel}
Let 
$$
\xymatrix{
{\cal X} \ar[r]^{\Phi} \ar[rd]_F & {\cal Y} \ar[d]_{\pi} \\
 & T \ar@/_/[u]_{\theta}
}
$$
%$$
%{\cal X} \stackrel{\Phi}\lorw {\cal Y} \stackrel{\pi}\lorw T 
%\stackrel{\theta}\lorw {\cal Y}
%$$
be projective maps of quasi-projective varieties such that:

\n
1)
${\cal X}$ is nonsingular of dimension $n$, $T$ is nonsingular
of dimension $s;$ 

\n
2) $F:= \pi \circ \Phi$ is surjective and smooth
of relative dimension $n-s;$

\n
3) the map $\Phi$ is stratified in the sense
of Theorem \ref{tila} and the strata of ${\cal Y}$
map  smoothly and surjectively onto $T;$

\n
4) $\theta $ is a section of $\pi,$ i.e. $\pi \circ \theta = Id_{T}$
 and
$\theta (T)$ is a stratum of ${\cal Y};$

\n
5) there is an isomorphism $ \fxn{\Phi}{\cal X}{n} \simeq
\dsdix{l}{\fxn{\Phi}{\cal X}{n} }.$

\n
Then there is a surjective map of local systems on $T:$
$$
\sigma \; : \; R^{n-s}F_{*} \rat_{\cal X} \lorw 
{\cal H}^{-s}(\theta^* \phix{0}{\Phi_*\rat_{\cal X}[n]}).
$$
In particular, the local system  
${\cal H}^{-s}(\theta^* \phix{0}{\Phi_*\rat_{\cal X}[n]})$
on $T$ is semisimple. 
\end{lm}
{\em Proof.}
By assumption 3),  the sheaves in question are indeed
local systems on $T.$
The closed embedding  $i: t \lorw  T$ of a point in $T$ induces
the following  diagram with Cartesian squares:
$$
\xymatrix{
{\cal X}_{t}  \ar[r]^{\Phi_{t}} \ar[d]^J & {\cal Y}_t \ar[r]^{\pi_t} \ar[d]^j & t \ar@/^/[l]^{\theta_t} \ar[d]^i \\
X  \ar[r]^{\Phi} & {\cal Y} \ar[r]^{\pi} & T \ar@/^/[l]^{\theta}
}
$$
%$$
%\begin{array}{ccccccc}
%{\cal X}_{t} & \stackrel{\Phi_{t}}\lorw & {\cal Y}_{t} &
%\stackrel{\pi_{t}}\lorw & t & \stackrel{\theta_{t}}\lorw  & {\cal 
%Y}_{t} \\
%\downarrow J &  & \downarrow j &  & \downarrow i  &  &
%\downarrow j\\
%{\cal X} & \stackrel{\Phi}\lorw & {\cal Y} & \stackrel{\pi}\lorw & T &
%\stackrel{\theta}\lorw & {\cal Y}. 
%\end{array}    
%$$
There is the commutative diagram
%$$
%\begin{array}{ccccc}
%    \pi_{*} \fxn{\Phi}{\cal X}{n} & \stackrel{a}\to  & 
%    \bigoplus_{l}{ \theta^{*} \phix{l}{ \fxn{\Phi}{\cal X}{n} } [-l] }
%    & \stackrel{b}\to  & \theta^{*} \phix{0}{\fxn{\Phi}{\cal X}{n}} \\
%    \downarrow &  & \downarrow && \downarrow
%    \\
%  i_*i^{*} \pi_{*} \fxn{\Phi}{\cal X}{n} & \to & 
%   i_*i^{*}  \bigoplus_{l}{ \theta^{*} \phix{l}{ \fxn{\Phi}{\cal X}{n} } [-l] }
%    & \to  & 
%    i_* i^{*}\theta^{*} \phix{0}{\fxn{\Phi}{\cal X}{n}}  \\
%    \simeq &  & \simeq && \simeq
%    \\
%    i_* {\pi_{t}}_{*} \fxn{\Phi_{t}}{{\cal X}_{t}}{n-s} [s] & 
%    \to & 
%   i_*  \bigoplus_{l}{ \theta^{*}_{t} j^{*} \phix{l}{ \fxn{\Phi}{\cal X}{n} }
%   [-l] }
%    & \to  & 
%    i_*\theta^{*}_{t} j^{*} \phix{0}{\fxn{\Phi}{\cal X}{n}} \\
%     \simeq &  & \simeq && \simeq
%    \\
%    i_* {\pi_{t}}_{*} \fxn{\Phi_{t}}{{\cal X}_{t}}{n-s} [s] & 
%    \to & 
%   i_*  \bigoplus_{l}{ \theta^{*}_{t}  \phix{l}{ \fxn{\Phi_{t}}{{\cal 
%   X}_{t}}{n-s}[s] }
%   [-l] }
%    & \to  & 
%    i_*\theta^{*}_{t}  \phix{0}{\fxn{{\Phi}_{t}}{{\cal 
%    X}_{t}}{n-s}} [s]
%    \end{array}
%$$
$$
\xymatrix{
\pi_{*} \fxn{\Phi}{\cal X}{n} \ar[r]^a \ar[d] & \bigoplus_{l}{ \theta^{*} \phix{l}{ \fxn{\Phi}{\cal X}{n} } [-l] }
    \ar[r]^{b}  \ar[d] & \theta^{*} \phix{0}{\fxn{\Phi}{\cal X}{n}} \ar[d] \\
i_*i^{*} \pi_{*} \fxn{\Phi}{\cal X}{n} \ar[r] \ar[d]^{\simeq} & 
   i_*i^{*}  \bigoplus_{l}{ \theta^{*} \phix{l}{ \fxn{\Phi}{\cal X}{n} } [-l] }
    \ar[r]\ar[d]^{\simeq} & 
    i_* i^{*}\theta^{*} \phix{0}{\fxn{\Phi}{\cal X}{n}} \ar[d]^{\simeq} \\
i_* {\pi_{t}}_{*} \fxn{\Phi_{t}}{{\cal X}_{t}}{n-s} [s]  
    \ar[r] \ar[d]^{=} & 
   i_*  \bigoplus_{l}{ \theta^{*}_{t} j^{*} \phix{l}{ \fxn{\Phi}{\cal X}{n} }
   [-l] }
    \ar[r]  \ar[d]^{\simeq} & 
    i_*\theta^{*}_{t} j^{*} \phix{0}{\fxn{\Phi}{\cal X}{n}} \ar[d]^{\simeq} \\
i_* {\pi_{t}}_{*} \fxn{\Phi_{t}}{{\cal X}_{t}}{n-s} [s] \ar[r] &
    i_*  \bigoplus_{l}{ \theta^{*}_{t}  \phix{l}{ \fxn{\Phi_{t}}{{\cal 
   X}_{t}}{n-s}[s] }
   [-l] }
   \ar[r]  & 
    i_*\theta^{*}_{t}  \phix{0}{\fxn{{\Phi}_{t}}{{\cal 
    X}_{t}}{n-s}} [s]
    }
$$
which is obtained as follows.
The first row: the first map is obtained by applying
$\pi_{*}$
 to the adjunction map for $\theta$ and by using 4) and 5); the second map
 is the natural projection.
 The first column of maps, is the adjunction relative to $i.$
The third  row is obtained from the second one
 using the usual base change relations:
 $i^{*}\pi_{*} \simeq {\pi_{t}}_{*}j^{*},$ 
 $j^{*}\Phi_{*} \simeq {\Phi_{t}}_{*}J^{*}$ and the equality
 $i^{*} \theta^{*} = \theta_{t}^{*} j^{*}.$
 The commutativity of the bottom follows from
 Lemma \ref{dtpb} in view of the fact that
  the  codimension $s$ embedding  $j: {\cal Y}_{t} \lorw {\cal Y}$ 
 is transverse to all 
 the strata of ${\frak Y}$ by 3).
 
\n
Keeping in mind 2), define $\sigma$ to be ${\cal H}^{-s}( b\circ a ).$

\n 
The map $\sigma_{t}$ at the level of stalks is identified
to the analogous map on the bottom row which reads as
$
H^{n-s}( {\cal X}_{t} ) \lorw H^{n-s}_{0}( \Phi_{t}^{-1}(t) )
$
and is surjective by Proposition \ref{pcdsurj}.

\n
We conclude by the Semisimplicity Theorem \ref{dss}.
 \blacksquare

\medskip
The following theorem concludes the proof of the semisimplicity
of $\phix{0}{\fxn{f}{X}{n}},$ for it shows that every
direct summand of it is an intersection cohomology complex associated 
with a semisimple local system on some locally closed
smooth subvariety.

\begin{tm}
\label{sss}
%Let $f: X \to Y,$
%$s \geq 0$ and $S_{s}$ be as in
%\ref{defws}.
%\n
The local systems ${\cal H}^{-s}( 
\alpha^{*}_{s}\phix{0}{\fxn{f}{X}{n}} )$
are semisimple. 
\end{tm}
{\em Proof.}
The statement is trivial
for $s=0$ and $s > m = \dim{f(X)}.$
 The case $s=m$  follows by Remark
\ref{tassok}. 
We may  assume that $1 \leq s  \leq m-1.$
We shall  reduce this  case to Lemma  \ref{pcdsrel}.

\n
In view of Remark \ref{sszos}, it is enough to show
semisimplicity over a Zariski-dense open subset
$T$ of every 
  connected component of $S_{s}.$
  
\n
Let 
$\pd = | A | \simeq {\Bbb P}^{d'}$ be the very ample linear system
on $Y$  associated with $A,$  
$\Pi: = ( \pd  )^{s},$ $d:= sd' = \dim{\Pi}.$
A point $p \in \Pi$ corresponds to
an $s-$tuple  $(H_{1}, \ldots , H_{s})$ of hyperplanes in ${\Bbb P}.$

\n
Consider the universal $s-$fold complete
intersection families ${\cal Y}:
= \{ (y,p) \, | \; y \in  \cap_{j=1}^{s}{ H_{j}} \} \subseteq Y \times \Pi$
and  ${\cal X}: = {\cal Y} \times_{Y \times \Pi} (X \times
\Pi )  \subseteq X \times \Pi.$
Note that ${\cal X}$ is nonsingular and the general member of the 
family ${\cal X}$ over $\Pi$ is nonsingular and connected by the 
Bertini 
Theorems;
in fact the assumption  $1 \leq s  \leq m-1$
implies $\dim f(X) \geq 2.$ 
However,
the connectedness plays no essential role.

\n
For every map $W \to Y$ there is the commutative diagram with cartesian squares
$$
\xymatrix{
{\cal X}_{W}  \ar[r] \ar[d]& {\cal X}  \ar[r]\ar[d] & X \times \Pi  
\ar[r]^{p}\ar[d]^{f'} & X \ar[d]^f \\
{\cal Y}_{W} \ar[r] \ar[d] & {\cal Y}  \ar[r] \ar[d]  & Y \times \Pi  
\ar[r]^{q} & Y \\
W    \ar[r] & Y. }
$$
%$$
%\begin{array}{ccccccc}
%{\cal X}_{W}  & \lorw & {\cal X} & \lorw & X \times \Pi & 
%\stackrel{p}\lorw & X \\
% \downarrow &  & \downarrow g &  & \downarrow f'&  &
%\downarrow f  \\
%{\cal Y}_{W} & \lorw & {\cal Y} & \lorw  & Y \times \Pi & 
%\stackrel{q}\lorw & Y \\
%  \downarrow &  & \downarrow  &&&& \\
%W &   \lorw & Y. &&
%\end{array}
%$$
The base-point-freeness of $|A|$ implies that
${\cal Y}_{W} \lorw W$ is Zariski-locally trivial.
Since
the general complete intersection of $s$ hyperplanes
meets the $s-$dimensional $S$ in a non-empty and finite set,  
the natural map 
$$
b\; : \;  {\cal Y}_{S} \lorw \Pi
$$
is dominant.

%Let ${\frak X}$ and ${\frak Y}$ be a stratification for
%$g$. In particular, $\phix{l}{\fxn{g}{\cal X}{\dim{\cal X}}}$ 
%is ${\frak Y}-$cc for every $l.$

%Let $S' \subseteq  S$ be a Zariski-dense open subset
%over which the family, ${\cal Y}_{S'}:= S'\times_{Y} {\cal Y},$ trivializes. 

\n
By  Bertini Theorem for $|L|$ (cf. \ref{cutstr}),
generic smoothness for ${\cal X} \lorw \Pi$, 
the algebraic Thom Isotopy Lemmas
(Theorem \ref{tila}) and 
generic smoothness for ${\cal Y} \lorw \Pi$,
there is   a Zariski-dense open subset
$\Pi^{0} \subseteq \Pi$ such that:

\n
1) the surjective map 
${\cal X} \lorw  \Pi$ is smooth over $\Pi^{0};$

\n
2) 
the complete intersections $Y_{s}$ of $s$ elements
associated with  the points of $ \Pi^{0}$ meet all strata
of $Y$ transversally;

\n
3) the restriction of  $h: {\cal Y} \lorw \Pi$ over $\Pi^{0}$ is  
stratified
so that every stratum maps surjectively and smoothly to
$\Pi^{0}.$

%\n
%Conditions 1) and 2) can be realized on some
%Zariski-open dense $V \subseteq \Pi$ by Bertini Theorem.
%Condition  3) is realized as follows.  
%Let $U''$ be the complement of the images of the closures
%of the strata of ${\frak Y}$ that do not dominate $\Pi.$
%By
%the Thom Isotopy Lemmas
%$h$ is topologically locally trivial
%on an  Zariski-open   dense subset $U'$ of $U''$.
%All strata of ${\frak Y}_{h^{-1} ( U' )}$ dominate 
%$U'$. By generic smoothness, there is a Zariski-dense
%open subset $U \subseteq U'$ over which all the strata
%of ${\frak Y}_{f^{-1}U}$ are smooth and map surjectively to $U.$

%\n
%Set $\Pi^{0} = U \cap V.$

\n 
Since $b$ is dominant, $b^{-1} \Pi^{0}$ is Zariski-dense
and open
in $ {\cal Y}_{S}.$ 

\n
Since ${\cal Y}_{S} 
\lorw S$ is  Zariski-locally trivial, 
there exists a Zariski-dense open subset $T\subseteq S$
 such that ${\cal Y}_{T} \lorw T$ 
 admits a section $\mu: T 
\lorw 
{\cal Y}_{T}$
 with the property that $\mu (T) \subseteq b^{-1} \Pi^{0}.$
 
 \n
 By shrinking  $T$, we may assume that
 the quasi-finite  map
  $b \circ \mu: $ 
 $T \lorw b(\mu (T))  \subseteq 
  \Pi^{0} \subseteq \Pi$ is smooth, of relative dimension zero.

 \n
We have a commutative diagram with cartesian squares
%$$
%\begin{array}{ccccc}
%{\cal X}_{T} & \stackrel{p'}\lorw & 
%{\cal X} & \stackrel{p''}\lorw & X \\
% \downarrow \Phi &  & \downarrow g &  &
%\downarrow f  \\
% {\cal Y}_{T} & \stackrel{q'}\lorw & 
%{\cal Y} & \stackrel{q''}\lorw & Y \\
% \downarrow \pi &  & \downarrow  h && \\
% T & \stackrel{b\circ \mu}\lorw & \Pi &&
%\end{array}
%$$

$$
\xymatrix{
{\cal X}_{T} \ar[r]^{p'} \ar[d]^{\Phi}& 
{\cal X}  \ar[r]^{p''}\ar[d]^{g} & X \ar[d]^f\\
 {\cal Y}_{T}  \ar[r]^{q'} \ar[d]^{\pi} & 
{\cal Y} \ar[r]^{q''} \ar[d]^h & Y \\
 T  \ar[r]^{b\circ \mu} & \Pi }
$$
The map $\Phi$ inherits a stratification
from the one on $g$ by pull-back and all strata
on ${\cal Y}_{T}$ map surjectively and smoothly onto $T.$

\n
For every $t \in T$, $Y_{t}:= \pi^{-1}(t)$ is a complete 
intersection of $s$ hyperplanes passing through
$t \in T \subseteq Y$, meeting all the strata
of $Y$ transversally and such that $ X_{t}:= (\pi \circ 
\Phi)^{-1}(t)$ is a smooth projective variety of dimension $n-s.$
Note that the map $\pi$ has a tautological section
$\theta: T \to  {\cal Y}_{T}$ assigning
to $t \in T$ the same point $t \in {\cal Y}_{t}.$

\n
We have a commutative diagram where the upper square   is cartesian:
%$$
%\begin{array}{ccc}
%{\cal X}_{T} & \stackrel{p_{X}}\lorw & X \\
%\downarrow \Phi &  & \downarrow f  \\
%{\cal Y}_{T} & \stackrel{p_{Y}}\lorw & Y \\
%\downarrow \pi & \nwarrow\theta & \uparrow \alpha_{T} \\
%T & = & T.
%\end{array}
%$$

$$
\xymatrix{
{\cal X}_{T} \ar[rr]^{p_X} \ar[d]^{\Phi} & & X \ar[d] ^f \\
{\cal Y}_{T} \ar[rr]^{p_Y} \ar[rd]^{\pi} & & Y \\
 &T \ar@/^/[ul]^{\theta} \ar[ur]^{\alpha_T} &
}
$$

\n
We have proved Theorem \ref{tm1}.b for $f$  so that $\fxn{f}{X}{n} \simeq 
\dsdix{l}{\fxn{f}{X}{n}}.$
By Lemma \ref{dtpb} and Remark \ref{joined}, 
we have, via pull-back, analogous  decompositions 
for $g$ and for $\Phi$ and an isomorphism
$p_{Y}^{*} \phix{0}{\fxn{f}{X}{n}} ) \simeq
\phix{0}{\fxn{\Phi}{{\cal X}_{T}}{n}} ).$

\n
We are now in the position to apply Lemma
\ref{pcdsrel} to the diagram ${\cal X}_{T} \stackrel{\Phi}\lorw
{\cal Y}_{T} \stackrel{\pi}\lorw T \stackrel{\theta}\lorw 
{\cal Y}_{T}$ and deduce the semisimplicity of
$$
{\cal H}^{-s}(\theta^{*}  \phix{0}{\fxn{\Phi}{{\cal X}_{T}}{n}} )
\simeq
{\cal H}^{-s}(\theta^{*} p_{Y}^{*} \phix{0}{\fxn{f}{X}{n}} )
\simeq {\cal H}^{-s} ( \alpha^{*}_{T}   \phix{0}{\fxn{f}{X}{n}}   ).
$$
We conclude by \ref{sszos} applied to $T\subseteq S.$
\blacksquare

\section{The pure Hodge structure on Intersection Cohomology }
\label{tphs}
In this section we prove the Purity Theorem  \ref{tmp1}
and the Hodge-Lefschetz Theorem for Intersection Cohomology
\ref{pic}.

\subsection{The Purity Theorem}
\label{rtoho}
Note that if $\dim{X} =1,$ then Theorem
\ref{tmp1} holds trivially.

\begin{lm}
\label{redtoo}
 If Theorem \ref{tmp1} holds for every  map $g:Z \lorw Z'$
 of projective varieties, $Z$ nonsingular, $\dim{Z} < \dim{X},$
 then it holds for every group 
 $H^{j}_{i}(X)$  $(i,j)\neq(0,n).$
\end{lm}
{\em Proof.}
Fix  $i<0.$  Let
$r: X^{1} \to   X$ be a nonsingular hyperplane section. 
Choose stratifications for $f$ and $g$ which have in common
the stratification of $Y.$

\n
By the Weak-Lefschetz-type 
Proposition \ref{wlmech}.i, the Semisimplicity Theorem
\ref{tm1}.c and   by the Hodge Structure Theorem \ref{uf} coupled with
Remark \ref{filtrfac}, the restriction 
map $r^{*}:  H^{j}_{i}(X)  \to   H^{j}_{i+1}(X^{1})$
is an injective  map of pure Hodge structures.

\n
By Remark \ref{nomaps},
the restriction map $r^{*}$ is a direct sum map
$$
r^{*} \, = \,
\sum r^{*}_{l,S}\, : \, \bigoplus_{l,S}{ H^{j}_{i,l,S}(X) } \lorw 
     \bigoplus_{l,S}{ H^{j}_{i+1, l,S}(X^{1})}.
$$
Let  $\widetilde{S}$ be a connected component of a stratum
$S_{\tilde{l}}.$
The inductive hypothesis holds for the  map
$g:= f \circ r: X^{1} \lorw Y.$ This implies that the natural
projection map
$$
\pi \, : \, 
\bigoplus_{l,S}{ H^{j}_{i+1, l,S}(X^{1})}
\lorw 
\bigoplus_{l,S\neq \tilde{S}}{ H^{j}_{i+1, l,S}(X^{1})}
$$
is a map of pure Hodge structures and so is
the composition $ \pi \circ r^{*}.$

\n
Clearly, $H^{j}_{i,\tilde{l}, \tilde{S}}(X) = \ke{\, (\pi \circ 
r^{*}) }$ is then a Hodge sub-structure for every $j \geq 0$
and every  $i <0.$

\n
The same argument as above, using  Proposition \ref{wlmech}.ii and 
a cokernel instead of a kernel, shows that any 
$H^{j}_{i,\tilde{l}, \tilde{S}}(X)$ is a 
Hodge sub-structure for every $j \geq 0$
and every  $i> 0.$

\n 
Let $i=0.$  There are two cases left: $j< n$ and $j>n.$
They are handled in the same way as above, by first replacing
$X^{1}$ with $X_{1}= f^{-1}(Y_{1})$
 the pre-image of a general hyperplane
section on $Y,$ and then by  using 
Proposition  \ref{pwl} instead of Proposition 
\ref{wlmech}.
\blacksquare

\medskip
Recall that the bilinear form $S^{\eta L}_{00}$
on $H^{n}_{0}(X)$ is induced by the Poincar\'e pairing
$\int_{X}{- \wedge -}$ on $X;$ see \ref{tiafiltr}, (\ref{sijmn}) and
(\ref{defofs}).
\begin{lm}
    \label{lm1} The direct sum decomposition
    $
    H^{n}_{0}(X) \, = \, \bigoplus_{S}{
     H^{n}_{0,l,S}(X)      }
    $ is $S^{\eta L}_{00}$-orthogonal.
    \end{lm}
 {\em Proof.} The duality isomorphism $\e: \rat_{X}[n] \simeq 
 \omega_{X}[-n]$ induces the isomorphism
 $$
 \phix{0}{\fxn{f}{X}{n}} \, \stackrel{  \phix{0}{\e}   }\simeq \,
 {\cal D}(\phix{0}{\fxn{f}{X}{n}})
 $$
 giving $S^{\eta L}_{00}$ in hypercohomology.
Setting $L_{0,l,S}:= {L_{0,l}}_{|S},$  this gives rise to an isomorphism 
 $$
 \bigoplus_{0,l,S}{ IC_{\overline{S}} (L_{0,l,S) }}  \simeq
 \bigoplus_{0,l,S}{ IC_{\overline{S}} (L_{0,l,S}^{\vee}) }
 $$
 which is a direct sum map by Remark \ref{nomaps}. 
 \blacksquare
    
    \begin{lm}
        \label{lm2}
        Let $V$ be a pure Hodge structure of weight $n.$
        $\Psi: V \otimes V \lorw \rat (-n)$
        be a map of pure Hodge structures which is nondegenerate
        as a bilinear form.
        %be of PHS and non-degenerate (sign convention on hodge II, page 26 
        %bottom and 27).
        Assume that $V= V_{1}\oplus V_{2}$
        with $V_{1} \subseteq V$ a  pure Hodge sub-structure and that
         $V_{1} \perp_{\Psi} V_{2}.$
        
        \n
        Then $V_{2} \subseteq V$ is a  pure Hodge sub-structure.
        \end{lm}
        {\em Proof.} The space  $V_{2}$ is the kernel of the 
        composition 
        $ V \to V^{\vee}  \to V^{\vee}_{1}.$
        \blacksquare

        \medskip
        \n
        {\bf {\em Proof of the Purity Theorem \ref{tmp1}.} }
The proof  is by induction on $\dim{X}.$

\n
The statement is trivial when $\dim{X}=1,$ for $f(X)$ is either a 
point or another curve and in either case there is only one  direct summand. 

\n
Assume that we have proved the statement for every map
$g: Z \to Z'$ of projective varieties with 
$Z$ nonsingular and $\dim{Z} < \dim{X}.$

\n
By Lemma \ref{redtoo} we are left with the case of $H^{n}_{0}(X).$

\n
Fix a connected component $S$ of a {\em non-dense} stratum $S_{l}.$
Let $L_{S}:= {L_{0,l}}_{|S}$ and $y \in S.$

\n
CLAIM: {\em $ \;
L_{S,y} = H^{n-l}_{0}( f^{-1}(y)) = H^{n-l}_{\leq 0} ( f^{-1}(y))
$
is a  weight-$(n-l)$  pure Hodge sub-structure of the mixed Hodge
structure $H^{n-l}(f^{-1}(y)).$}

\n
Proof of CLAIM. The first equality is the definition
of $L_{S}$ (see Proposition \ref{itissum}).  
Let $Y_{l} \subseteq $ be the intersection of $l$
sufficiently general
hyperplane sections of $Y$ through $y \in S$
and $f_{l} : X_{l}:=  f^{-1}(Y_{l})  \to Y_{l}$ be the resulting map.
Clearly, $f_{l}^{-1}(y)= f^{-1}(y).$ The filtrations
on the cohomology group $H^{n-l}(f^{-1}(y))$ induced by the two maps
$f_{l}$ and $f$ coincide by Lemma \ref{dtpb}.
The second equality follows from Lemma \ref{collectinfo}
applied to $f_{l}.$
Theorem \ref{nhrbr}, applied to $f_{l},$ gives the last statement
of the CLAIM.

\n
Let 
$$
\rho:\,  Z_{S} \lorw Z'_{S} := \,  \overline{f^{-1}(S)}
$$
be a proper surjective map, with $Z_{S}$ nonsingular and projective
and of dimension $\dim{Z_{S}} = \dim{Z'_{S}} <n.$
For example,  a resolution of the singularities of the irreducible components
of $Z'_{S}.$
Note that $Z_{S}$ is not necessarily pure-dimensional.

\n
Let $S_{sm} \subseteq S$ be the Zariski-dense open
set over which $f \circ \rho$ is smooth.

\n
By refining the stratification, we may assume that
$S_{sm} =S.$ In fact, the new strata in $S \setminus
S_{sm}$ will  not contribute any new direct summand to $H^{n}_{0}(X).$ 

\n
Setting  $g:= f \circ \rho: Z_{S} \to Y,$
the map
$$
\rho: g^{-1}(y) \lorw f^{-1}(y)
$$
is proper and  surjective from a nonsingular space.

\n
By Theorem \ref{propofmhs}, the map
$
H^{n-l}(  f^{-1}(y)  ) \stackrel{\rho^{*}}\lorw H^{n-l}( g^{-1}(y) )
$
is such that 
$$
\ke{ \rho^{*} } = W_{n-l-1} H^{n-l} (  f^{-1}(y)  ).
$$
Since $H^{n-l}_{\leq 0}( f^{-1}(y) )$ is of  pure weight
$n-l,$
$$
\ke{\, \rho^{*}} \cap H^{n-l}_{\leq 0} (  f^{-1}(y)  ) = \{0 \}
$$
so that
$$
\rho^{*}_{|} \, : \,
H^{n-l}_{\leq 0} (  f^{-1}(y)  ) \lorw H^{n-l} (  g^{-1}(y)  ) 
$$
is injective.  It follows that  so is the map of local systems
\begin{equation}
    \label{eq2}
L_{S} \lorw ( R^{n-l}g_{*} \rat_{Z_{S}})_{|S}.
\end{equation}
By Deligne's Semisimplicity Theorem \ref{dss}, this injection
splits.
Let
$$
Z_{S} = \coprod_{t \geq 0} Z^{t}_{S}
$$
be the decomposition into pure-dimensional
``components.''

\n
By Remark \ref{filtrfac} and the Hodge Structure 
Theorem  \ref{uf}, we get that the natural
pull-back  is a map of pure Hodge structures:
$$
\rho^{*}: \, H^{n}_{0}(X) \lorw \bigoplus_{0 \leq t \leq n-1 }{
H^{n}_{n-t}(Z^{t}_{S})}.
$$
Denote by $\pi: \oplus_{0 \leq t \leq n-1 }{
H^{n}_{n-t}(Z^{t}_{S})} \to V$ the projection
corresponding to the direct summand
associated with $(R^{n-s}g_{*} \rat_{Z_{S}})_{|S}.$

\n
Since $t<n,$ we can apply the inductive hypothesis and
$\pi$ is a map of  pure Hodge structures.

\n
By  (\ref{eq2}),
$$
\ke{\,( \pi \circ \rho^{*} ) } \,  = \,
\bigoplus_{l', \,S' \neq S}{H^{n}_{0,l',S'}(X)} \subseteq   H^{n}_{0}(X) 
$$
is a Hodge sub-structure.

\n
By Lemma \ref{lm1} and Lemma \ref{lm2}, we have that
$$
H^{n}_{0,l,S} (X) \subseteq H^{n}_{0}(X)
$$
is a pure Hodge sub-structure for every {\em non-dense} stratum $S.$
This implies immediately that so is any direct sum over non-dense 
strata.

\n
Applying the two lemmata again, we conclude that 
the contribution from the dense stratum is also
a pure Hodge sub-structure.
\blacksquare

\subsection{The Hodge-Lefschetz Theorem}
\label{poftmpic}
We need the following

\begin{lm}
    \label{nt}
Let 
$$
\xymatrix{
    \hat{X} \ar[r]^{r} \ar[dr]^{\hat{f}} & X \ar[d]^f \\
       & Y }
$$
 be such that $\hat{f}$ and $f$ are resolutions
 and $r$ is proper and surjective.
 Let $Y= \amalg{S_{l}}$ be  a stratification of $Y$
 part of stratifications for  $\hat{f}$ and for $f.$
 Let $y \in S$ be the choice of a point on a 
 connected component of a stratum.
 We have the diagram
 $$
\xymatrix{
    H^{n-l}(   f^{-1}(y) ) \ar[r]^{r^{*}}  & 
    H^{n-l}( \hat{f}^{-1}(y)   )\\
    H^{n-l}_{ \leq 0}(  f^{-1}(y)  )  \ar@{^{(}->}[u] \ar[r]^{r^{*}_{0}} 
    &  H^{n-l}_{\leq 0}(  \hat{f}^{-1}(y) ) \ar@{^{(}->}[u] .
    }
    $$
%$$
%\begin{array}{ccc}
%    H^{n-l}(   f^{-1}(y) )  & \stackrel{r^{*}}\lorw  & 
%    H^{n-l}( \hat{f}^{-1}(y)   )\\
%     UI & & UI   \\
%    H^{n-l}_{ \leq 0}(  f^{-1}(y)  )   & \stackrel{r^{*}_{0}}\lorw 
%    &  H^{n-l}_{\leq 0}(  \hat{f}^{-1}(y)     ) .
%    \end{array}
%    $$
Then $r^{*}_{0}$ is an injection of pure Hodge structures
and we have a splitting 
injection of local systems
$$
L^{f}_{0,l,S}  \lorw L^{\hat{f}}_{0,l,S}.
$$
\end{lm}
{\em Proof.} Let 
$$
\hat{Z}_{S} \stackrel{\rho}\lorw  \overline{ \hat{f}^{-1}(S)  }
$$
be any projective and  surjective map from a nonsingular space and 
 $\theta := \hat{Z}_{S} \to \overline{S}$ be the resulting map.

 \n
 We may assume, by refining the stratification if necessary, that
 $\theta$ is smooth over $S.$
 
\n
We have the commutative diagram
 
$$
\xymatrix{
  H^{n-s}_{ \leq 0}(  f^{-1}(y)  )    \ar[r]^{r^{*}_{0}} \ar[dr]_{\rho^{*} 
  \circ r^{*}} & 
  H^{n-s}_{\leq 0}(  \hat{f}^{-1}(y)     ) \ar[d]^{\rho^*} \\
      & H^{n-s}_{\leq 0}(  \theta^{-1}(y)     )}
$$
where $\rho^{*} \circ r^{*}$ is  injective by the same 
argument as in the proof of the Theorem \ref{tmp1}.

\n
It follows that $r^{*}_{0}$ is injective.
The existence of a splitting comes from semisimplicity.
\blacksquare

\bigskip
\n
{\bf {\em 
Proof of Theorem \ref{pic}.a.}}
There is a  commutative diagram  
$$
\xymatrix{
  H^{n+j}_{  0}( X )  \,\, \simeq  
  I\!H^{n+j}(Y)^{h}\; \bigoplus \bigoplus_{l\neq \dim{Y},S}{ 
  I\!H^{j}(Y, IC_{\overline{S}}( 
  L_{0,l,S}^{h}   ) )} \\
    H^{n+j}_{  0}( X'  ) \ar@<23ex>[u]^{g'^*} \simeq
  I\!H^{n+j}(Y)^{f'} \, \bigoplus \bigoplus_{l\neq \dim{Y},S}{ 
  I\!H^{j}(Y, IC_{\overline{S}}( 
  L_{0,l,S}^{f'}  ) )} \ar@<5.2ex>[u]
}
  $$
where, by \ref{nomaps},  the vertical map on the right is direct sum map of maps
of pure Hodge structures identifying the pure Hodge structure
on  $I\!H(Y)^{h}$ with 
the one on $I\!H(Y)^{\hat{f'}}(Y).$ 
\blacksquare

%\begin{rmk}
%    {\rm 
%    What about $r_{*} \rat_{\hat{X}} = \rat_{X} \oplus \ldots$? 
%    Shouldn't this give something helpful?
%    Also should try $r$ proper surjective, i.e. not necessarily 
%    resolution and see about ``functoriality'' for the pieces 
%    $H^{l}_{S,b}(X)$ not only $IH(Y)$ for resolution.
%    }
%    \end{rmk}
 
\bigskip
\n
{\bf {\em Proof of Theorem \ref{pic}.b.}} By transversality,
the complex   $r^{*}IC_{Y}[-1]$ satisfies the
conditions characterizing $ IC_{Y_{1}}$ (cf. \ref{iccso}). The result 
follows from Proposition \ref{pwl}.
\blacksquare

\bigskip
\n
{\bf {\em Proof of Theorem \ref{pic}.c.}} Let $f: X \to Y$
be a projective  resolution of the singularities of $Y.$
By the Semisimplicity
Theorem \ref{tm1}.c,
the complex $IC_{Y}$ is a direct summand of 
$\phix{0}{\fxn{f}{X}{\dim{X}}}.$ By Remark \ref{sepullback},
the cup product map $A^{j}$ is a direct sum map.
By Theorem \ref{tm3}, recalling that we are identifying $L^{j}$
with $A^{j},$ the map $A^{j}$ is an isomorphism on every direct summand,
whence the Hard Lefschetz-type statement and its 
 standard algebraic 
consequence, i.e. 
the primitive Lefschetz Decomposition. 

\blacksquare

\bigskip
\n
{\bf {\em Proof of Theorem \ref{pic}.d.}}
Since $f$ is birational,
 i)  the complexes $\phix{i}{\fxn{f}{X}{n}}$ 
 are supported on  proper closed algebraic  subsets of $Y$
 for every $i \neq 0$ and ii) 
  $IC_{Y}$ is a direct summand
of $\ke{ \, \eta } = {\cal P}_{\eta}^{0} \subseteq  
\phix{0}{\fxn{f}{X}{\dim{X}}}$ and is the only summand
supported on $Y.$

\n
Let $\eta$ be any ample line bundle on $X.$
The result follows from the 
Generalized Hodge-Riemann Bilinear Relations \ref{tmboh},
  \ref{pic}.a, Remark \ref{important}
and from the fact that $L$ acts compatibly with any direct sum 
decomposition by Remark \ref{sepullback}.
    \blacksquare

Authors' addresses:

\smallskip
\n
Mark Andrea A. de Cataldo,
Department of Mathematics,
SUNY at Stony Brook,
Stony Brook,  NY 11794, USA. \quad 
e-mail: {\em mde@math.sunysb.edu}

\smallskip
\n
Luca Migliorini,
Dipartimento di Matematica, Universit\`a di Bologna,
Piazza di Porta S. Donato 5,
40126 Bologna,  ITALY. \quad
e-mail: {\em migliori@dm.unibo.it}


\begin{thebibliography}{99}



\bibitem{bbd}{A.A. Beilinson, J.N. Bernstein, P. Deligne,
{\em Faisceaux pervers}, Ast\'erisque {\bf 100}, Pa\-ri\-s, Soc. Math. 
Fr. 1982.}


\bibitem{borel} A. Borel et al., {\em Intersection Cohomology}, 
Progress in Mathematics Vol. {\bf 50}, Birkh\"auser, Boston Basel 
Stuttgart 1984.

\bibitem{b-m}{W. Borho, R. MacPherson, ``Partial resolutions of 
nilpotent 
varieties," Ast\'erisque {\bf 101-102} (1983), 23-74.}

%\bibitem{cg} N. Chriss, V. Ginzburg, {\em Representation theory and 
%complex geometry}. 
%Bi\-r\-kh\"auser Boston, 1997.

\bibitem{clemens} C.H. Clemens,
``Degeneration of K\"ahler manifolds,''
Duke Math. J. {\bf 44} (1977), no. 2, 215--290.


\bibitem{demigsemi}{ M. de Cataldo, L. Migliorini, ``The 
Hard Lefschetz Theorem and the Topology of semismall maps,"      
Ann. Sci. \'Ecole Norm. Sup. (4) {\bf 35} (2002), no. 5, 759--772.}

\bibitem{decmigmot}{ M. de Cataldo, L. Migliorini, ``The Chow Motive 
of semismall resolutions,''
arXiv:math.AG/0204067 (2002), to appear in Math.Res.Lett.}

\bibitem{decmightamv1}{M. de Cataldo, L. Migliorini, ``The 
Hodge Theory of Algebraic maps,''
arXiv:math.AG/0306030 v1 (2003).}


\bibitem{dess} P. Deligne, ``Th\'eor\`eme de Lefschetz et crit\`eres 
de d\'eg\'en\'erescence
de suites spectrales,'' Publ.Math. IHES {\bf 35} (1969), 107-126.

\bibitem{shockwave} P. Deligne, ``D\'ecompositions dans la cat\'egorie 
D\'eriv\'ee'', 
Motives (Seattle, WA, 1991), 115--128, Proc. Sympos. Pure Math., {\bf 
55},Part 1, Amer. Math. Soc., Providence, RI, 1994.
 
\bibitem{ho2}  P. Deligne, ``Th\'eorie de Hodge, II,'' Publ.Math. 
IHES {\bf 40} (1971), 5-57.

\bibitem{ho3}  P. Deligne, ``Th\'eorie de Hodge, III,'' Publ.Math. 
IHES {\bf 44} (1974), 5-78.

\bibitem{weil2} P. Deligne, ``La conjecture de Weil, II,'' Publ.Math. 
IHES {\bf 52} (1980), 138-252.

\bibitem{elzein} F. El Zein, ``Th\'eorie de Hodge des cycles 
\'evanescents''
Ann. Sci. \'Ecole Norm. Sup. (4) {\bf 19} (1986), 107--184.


%\bibitem{e-v} H. Esnault, E. Viehweg, ``Vanishing and Non-Vanishing 
%Theorems,"
%Actes du Colloque de Th\'eorie de Hodge, 
%Asterisque {\bf 179-180} 1989, 97-112.



\bibitem{go-ma2} M. Goresky, R. MacPherson, ``Intersection homology 
II,''
Inv. Math. {\bf 71} (1983), 77-129. 


\bibitem{g-m} M. Goresky, R. MacPherson, {\em Stratified Morse 
Theory},
 Ergebnisse der Mathematik,
3.folge. Band 2, Springer-Verlag, Berlin Heidelberg 1988.

\bibitem{navarro} F. Guillen, V. Navarro Aznar, ''Sur le th\'eor\`eme 
local des cycles invariants,'' 
Duke Math. J. {\bf 61} (1990), no. 1, 133--155.

\bibitem{iv} B. Iversen, {\em Cohomology of Sheaves}, Universitext, 
Springer-Verlag, Berlin Heidelberg 1986.

\bibitem{k-s} M. Kashiwara, P. Schapira {\em Sheaves on manifolds},
Grundlehren der mathematischen Wissenschaften. Vol. 292, 
Springer-Verlag, Berlin Heidelberg 1990.


%\bibitem{laufer} H. Laufer, {\em Normal Two-Dimensinal Singularities},
%Annals of Mathematics Studies {\bf 71}, Princeton University Press 
%1971.



%\bibitem{looj} E. Looijenga, ``Cohomology
%and Intersection Hohomology
%of Algebraic Varieties,'' in {\em Complex Algebraic Geometry}, 
%IAS/Park City 
%Mathematics Series, Vol.3 J. Koll\'ar Editor, American Mathematical 
%Society, 1997, 221-264.


\bibitem{mac83}{ R. MacPherson}{ ``Global Questions in the Topology
of Singular Spaces,'' Proc. of the I.C.M.,  213--235 (1983), Warszawa.}

\bibitem{mo1}{T. Mochizuki, ``Asymptotic behaviour of tame nilpotent
harmonic bundles with trivial parabolic structure,'' J.Diff.Geom.
{\bf 62} (2002), 351-559.}

\bibitem{mo2}{T. Mochizuki, ``Asymptotic behaviour of tame harmonic
 bundles and an application to pure twistor D-modules,''
arXiv:math.DG/0312230 (2003).}

\bibitem{sabbah}{C. Sabbah, ``Polarizable twistor D-modules,''
preprint.}

 \bibitem{samhp}{    
 M. Saito, `` Modules de Hodge polarisables,''
 Publ. Res. Inst. Math. Sci. 
{\bf 24} (1988), no.~6, 849--995 (1989).}

\bibitem{samhm}
{ M. Saito,
`` Mixed Hodge modules,''
 Publ. Res. Inst. Math. Sci. {\bf 26} (1990), no.~2, 221--333.}

\bibitem{sa}{M. Saito, ``Decomposition theorem for proper K\"ahler 
morphisms,"
Tohoku Math. J. (2) {\bf 42}, no. 2,  (1990),  127--147.}

\bibitem{so} A.J. Sommese, ``Submanifolds of Abelian Varieties,'' 
Math. Ann.{\bf  233} (1978), no. 3, 229--256. 

\bibitem{steenbrink} J. Steenbrink, ``Limits of Hodge Structures,''
Inv. Math. {\bf 31} (1975-76), 229-257. 

\bibitem{stzu}{ J. Steenbrink, S. Zucker, 
`` Variation of mixed Hodge structure. I''.  Invent. Math.  {\bf 80}  
(1985),  no. 3, 489--542.}

\end{thebibliography}
\end{document}